\documentstyle[a4,11pt,bj,etcetera,xy]{article}
\newcommand{\comment}[1]{ }

\newcommand{\ups}{^{\ast}}

\newcommand{\iso}{\cong}

\newcommand{\adj}[6]{\mbox{$#1,#2:#3\dashv#4:#5\rightarrow#6$}}

\newcommand{\Set}{\mbox{${\cal S}\!et$}}

\newcommand{\morph}[3]{\mbox{$#1:#2\rightarrow #3$}}

\newcommand{\isomorph}[3]{\mbox{$#1:#2\stackrel{\sim}{\rightarrow}#3$}}

\newcommand{\cell}[3]{\mbox{$#1:#2\Rightarrow #3$}}

\newcommand{\CommaObj}[2]{\mbox{$#1\!\downarrow\! #2$}}

\newcommand{\twocatdef}[4]{
\setbox153 = \hbox{#1{\bf\hspace*{5mm}morphisms\hspace*{5mm}}}
\catdefsize\textwidth
\advance\catdefsize-\@totalleftmargin
\advance\catdefsize-\wd153
\begin{tabbing}
#1\={\bf\hspace*{5mm}morphisms\hspace*{5mm}}\=\kill
#1\>{\bf\hspace*{5mm}objects}\>\parbox[t]{\catdefsize}{#2} \\
\>{\bf\hspace*{5mm}morphisms}\>\parbox[t]{\catdefsize}{#3} \\
\>{\bf\hspace*{5mm}2-cells}\>\parbox[t]{\catdefsize}{#4}
\end{tabbing}}

\newtheorem{assumption}[definition]{Assumption}

\def\sslice{/{\kern-.7ex}/}

\newcommand{\alg}[1]{#1
         \mathchoice{\mbox{-}}{\mbox{-}}{\mbox{\scriptsize -}}{\mbox{\tiny -}}
{\cal A}lg}

\def\wrt.{w.r.t.}

\newcommand{\Mod}{{\cal M}\!od}

\def\ucorner#1{\raisebox{.5ex}{$\ulcorner$}{\kern-.4ex}#1
         {\kern-.4ex}\raisebox{.5ex}{$\urcorner$}}
\def\lcorner#1{\raisebox{-.5ex}{$\llcorner$}{\kern-.4ex}#1
         {\kern-.4ex}\raisebox{-.5ex}{$\lrcorner$}}

\newcommand{\cat}[1]{\mbox{$\Bbb #1$}}

\newcommand{\Cat}{\mbox{${\cal C\/}\!at$}}

\renewcommand{\frac}[2]{
\prooftree #1
\justifies
#2
\thickness=0.08em
\endprooftree}

\newcommand{\tuple}[1]
   {\langle #1 \rangle}

\def\Multicat{\mbox{${\cal M}\!\mathit{ulticat}$}}
\def\Repmulticat{\mbox{${\cal R}\!\mathit{ep}\Multicat$}}
\def\MonCat{\mbox{${\cal M}\!\mathit{on}\Cat$}}

\def\Pbk{\bold{Pbk}}

\def\wgph{\mbox{$\omega\mbox{\textsf{-gph}}$}}
\def\wcat{\mbox{$\omega\mbox{\textsf{-cat}}$}}

\newcommand{\ncat}[1]{\mbox{$#1\mbox{\textsf{-cat}}$}}
\newcommand{\ngph}[1]{\mbox{$#1\mbox{\textsf{-gph}}$}}
\newcommand{\intcat}[1]{\mbox{$\mathbf{#1}$}}
\newcommand{\pseudo}[1]{\mbox{${\mbox{\textsf{Ps-}}}{#1}$}}
\newcommand{\lax}[1]{\mbox{${\mbox{\textsf{Lax-}}}{#1}$}}
\newcommand{\laxrep}[1]{\mbox{${\mathsf{Lax}_{\mathit{rep}}}{#1}$}}
\newcommand{\Alg}[1]{\mbox{$#1\mbox{\textsf{-alg}}$}}

\newcommand{\Mnd}[1]{\mbox{${\mathsf{Mnd}}({#1})$}}
\renewcommand{\Mod}[1]{#1\mbox{\sf -mod}}

\newcommand{\ignore}[1]{}

\newcommand{\Spn}[1]{\bold{Spn}(#1)}
\newcommand{\SpnT}[2]{\bold{Spn}_{#1}(#2)}

\newcommand{\bimod}[1]{\mbox{$#1\!{\mbox{\textsf{-bimod}}}$}}
\newcommand{\Bimod}[1]{\mbox{${\mathsf{Bimod}}({#1})$}}
\newcommand{\BimodT}[2]{\mbox{${\mathsf{Bimod}}_{{#1}}({#2})$}}
\newcommand{\tensor}{\mathrel{\otimes}}
\newcommand{\lineardual}{o}
\newcommand{\Hom}[1]{\mbox{${\mathsf{Hom}}_{{#1}}$}}
\newcommand{\im}[1]{#1_{\#}}
\newcommand{\subadj}[1]{#1_{\!\scriptscriptstyle\dashv}}

\newcommand{\twocat}[1]{\mbox{$\mathbf{\mathcal{#1}}$}}

\newcommand{\bold}[1]{\mathbf{#1}}
\newcommand{\qed}{\hspace*{\fill}$\Box$}
\newcommand{\rel}[3]{{{#1}\colon{#2}\not\rightarrow{#3}}}

\input amssym.def
\input{amssym}
\newdimen\proofrulebreadth \proofrulebreadth=.05em
\newdimen\proofdotseparation \proofdotseparation=1.25ex
\newdimen\proofrulebaseline \proofrulebaseline=2ex
\newcount\proofdotnumber \proofdotnumber=3
\let\then\relax
\def\hfi{\hskip0pt plus.0001fil}
\mathchardef\squigto="3A3B
\newif\ifinsideprooftree\insideprooftreefalse
\newif\ifonleftofproofrule\onleftofproofrulefalse
\newif\ifproofdots\proofdotsfalse
\newif\ifdoubleproof\doubleprooffalse
\let\wereinproofbit\relax
\newdimen\shortenproofleft
\newdimen\shortenproofright
\newdimen\proofbelowshift
\newbox\proofabove
\newbox\proofbelow
\newbox\proofrulename
\def\shiftproofbelow{\let\next\relax\afterassignment\setshiftproofbelow\dimen0 }
\def\shiftproofbelowneg{\def\next{\multiply\dimen0 by-1 }%
\afterassignment\setshiftproofbelow\dimen0 }
\def\setshiftproofbelow{\next\proofbelowshift=\dimen0 }
\def\setproofrulebreadth{\proofrulebreadth}

\def\prooftree{
\ifnum  \lastpenalty=1
\then   \unpenalty
\else   \onleftofproofrulefalse
\fi
\ifonleftofproofrule
\else   \ifinsideprooftree
        \then   \hskip.5em plus1fil
        \fi
\fi
\bgroup
\setbox\proofbelow=\hbox{}\setbox\proofrulename=\hbox{}%
\let\justifies\proofover\let\leadsto\proofoverdots\let\Justifies\proofoverdbl
\let\using\proofusing\let\[\prooftree
\ifinsideprooftree\let\]\endprooftree\fi
\proofdotsfalse\doubleprooffalse
\let\thickness\setproofrulebreadth
\let\shiftright\shiftproofbelow \let\shift\shiftproofbelow
\let\shiftleft\shiftproofbelowneg
\let\ifwasinsideprooftree\ifinsideprooftree
\insideprooftreetrue
\setbox\proofabove=\hbox\bgroup$\displaystyle 
\let\wereinproofbit\prooftree
\shortenproofleft=0pt \shortenproofright=0pt \proofbelowshift=0pt
\onleftofproofruletrue\penalty1
}

\def\eproofbit{
\ifx    \wereinproofbit\prooftree
\then   \ifcase \lastpenalty
        \then   \shortenproofright=0pt  
        \or     \unpenalty\hfil         
        \or     \unpenalty\unskip       
        \else   \shortenproofright=0pt  
        \fi
\fi
\global\dimen0=\shortenproofleft
\global\dimen1=\shortenproofright
\global\dimen2=\proofrulebreadth
\global\dimen3=\proofbelowshift
\global\dimen4=\proofdotseparation
\global\mscount=\proofdotnumber
$\egroup  
\shortenproofleft=\dimen0
\shortenproofright=\dimen1
\proofrulebreadth=\dimen2
\proofbelowshift=\dimen3
\proofdotseparation=\dimen4
\proofdotnumber=\mscount
}

\def\proofover{
\eproofbit 
\setbox\proofbelow=\hbox\bgroup 
\let\wereinproofbit\proofover
$\displaystyle
}%
\def\proofoverdbl{
\eproofbit 
\doubleprooftrue
\setbox\proofbelow=\hbox\bgroup 
\let\wereinproofbit\proofoverdbl
$\displaystyle
}%
\def\proofoverdots{
\eproofbit 
\proofdotstrue
\setbox\proofbelow=\hbox\bgroup 
\let\wereinproofbit\proofoverdots
$\displaystyle
}%
\def\proofusing{
\eproofbit 
\setbox\proofrulename=\hbox\bgroup 
\let\wereinproofbit\proofusing
\kern0.3em$
}

\def\endprooftree{
\eproofbit 
  \dimen5 =0pt
\dimen0=\wd\proofabove \advance\dimen0-\shortenproofleft
\advance\dimen0-\shortenproofright
\dimen1=.5\dimen0 \advance\dimen1-.5\wd\proofbelow
\dimen4=\dimen1
\advance\dimen1\proofbelowshift \advance\dimen4-\proofbelowshift
\ifdim  \dimen1<0pt
\then   \advance\shortenproofleft\dimen1
        \advance\dimen0-\dimen1
        \dimen1=0pt
        \ifdim  \shortenproofleft<0pt
        \then   \setbox\proofabove=\hbox{%
                        \kern-\shortenproofleft\unhbox\proofabove}%
                \shortenproofleft=0pt
        \fi
\fi
\ifdim  \dimen4<0pt
\then   \advance\shortenproofright\dimen4
        \advance\dimen0-\dimen4
        \dimen4=0pt
\fi
\ifdim  \shortenproofright<\wd\proofrulename
\then   \shortenproofright=\wd\proofrulename
\fi
\dimen2=\shortenproofleft \advance\dimen2 by\dimen1
\dimen3=\shortenproofright\advance\dimen3 by\dimen4
\ifproofdots
\then
        \dimen6=\shortenproofleft \advance\dimen6 .5\dimen0
        \setbox1=\vbox to\proofdotseparation{\vss\hbox{$\cdot$}\vss}
        \setbox0=\hbox{%
                \kern\dimen6
                $\vcenter to\proofdotnumber\proofdotseparation
                        {\leaders\box1\vfill}$%
                \unhbox\proofrulename}%
\else   \dimen6=\fontdimen22\the\textfont2 
        \dimen7=\dimen6
        \advance\dimen6by.5\proofrulebreadth
        \advance\dimen7by-.5\proofrulebreadth
        \setbox0=\hbox{%
                \kern\shortenproofleft
                \ifdoubleproof
                \then   \hbox to\dimen0{%
                        $\mathsurround0pt\mathord=\mkern-6mu%
                        \cleaders\hbox{$\mkern-2mu=\mkern-2mu$}\hfill
                        \mkern-6mu\mathord=$}%
                \else   \vrule height\dimen6 depth-\dimen7 width\dimen0
                \fi
                \unhbox\proofrulename}%
        \ht0=\dimen6 \dp0=-\dimen7
\fi
\let\doll\relax
\ifwasinsideprooftree
\then   \let\VBOX\vbox
\else   \ifmmode\else$\let\doll=$\fi
        \let\VBOX\vcenter
\fi
\VBOX   {\baselineskip\proofrulebaseline \lineskip.2ex
        \expandafter\lineskiplimit\ifproofdots0ex\else-0.6ex\fi
        \hbox   spread\dimen5   {\hfi\unhbox\proofabove\hfi}%
        \hbox{\box0}%
        \hbox   {\kern\dimen2 \box\proofbelow}}\doll%
\global\dimen2=\dimen2
\global\dimen3=\dimen3
\egroup 
\ifonleftofproofrule
\then   \shortenproofleft=\dimen2
\fi
\shortenproofright=\dimen3
\onleftofproofrulefalse
\ifinsideprooftree
\then   \hskip.5em plus 1fil \penalty2
\fi
}

\xyoption{v2}
\xyoption{2cell}
\xyoption{dvips}
\CompileMatrices

\LaTeXdiagrams

\UseAllTwocells

\ignore{
\def\Multicat{\mbox{${\cal M}\!\mathit{ulticat}$}}
\def\Repmulticat{\mbox{${\cal R}\!\mathit{ep}\Multicat$}}
\def\MonCat{\mbox{${\cal M}\!\mathit{on}\Cat$}}

\def\Pbk{\bold{Pbk}}

\def\wgph{\mbox{$\omega-\mathsf{gph}$}}
\def\wcat{\mbox{$\omega-\mathsf{cat}$}}

\newcommand{\ncat}[1]{\mbox{$#1-\mathsf{cat}$}}
\newcommand{\ngph}[1]{\mbox{$#1-\mathsf{gph}$}}
\newcommand{\intcat}[1]{\mbox{$\mathbf{#1}$}}
\newcommand{\pseudo}[1]{\mbox{${\mathsf{Ps}}\!-\!{#1}$}}
\newcommand{\lax}[1]{\mbox{${\mathsf{Lax}}\!-\!{#1}$}}
\newcommand{\laxrep}[1]{\mbox{${\mathsf{Lax}_{\mathit{rep}}}{#1}$}}
\newcommand{\Alg}[1]{\mbox{$#1\!-\!\mathsf{alg}$}}

\renewcommand{\Mnd}[1]{\mbox{${\mathsf{Mnd}}({#1})$}}
\newcommand{\Mod}[1]{#1\!-\!\mbox{\sf mod}}

\newcommand{\comment}[1]{}
\newcommand{\ignore}[1]{}

\newcommand{\Spn}[1]{\bold{Spn}(#1)}
\newcommand{\SpnT}[2]{\bold{Spn}_{#1}(#2)}

\newcommand{\bimod}[1]{\mbox{$#1\!-\!{\mathsf{bimod}}$}}
\newcommand{\Bimod}[1]{\mbox{${\mathsf{Bimod}}({#1})$}}
\newcommand{\BimodT}[2]{\mbox{${\mathsf{Bimod}}_{{#1}}({#2})$}}
\newcommand{\tensor}{\mathrel{\otimes}}
\newcommand{\lineardual}{o}
\newcommand{\Hom}[1]{\mbox{${\mathsf{Hom}}_{{#1}}$}}
\newcommand{\im}[1]{#1_{\#}}
\newcommand{\subadj}[1]{#1_{\!\scriptscriptstyle\dashv}}

\renewcommand{\twocat}[1]{\mbox{$\mathbf{\mathcal{#1}}$}}

}

\begin{document}

\begin{center}
{\Large\bf  
From coherent structures to universal 
properties}                                  
\vskip5mm
\it Claudio Hermida{$^{*}$}
          \def\thefootnote{}
          \footnotetext{
            \hbox{}\hskip-1.8em{*} 
CMA, Mathematics Department, IST\\
Lisbon, Portugal.\\
e-mail:{\/\tt 
           chermida@math.ist.utl.pt}.}
\end{center}

\begin{center}                           
\number\day\space
\ifcase\month\or January\or February\or March\or April\or May\or
June\or July\or August\or September\or October\or November\or
December\fi \space\number\year
\end{center}

\begin{abstract}
  Given a 2-category $\twocat{K}$ admitting a calculus of bimodules, and a 
2-monad $\mathsf{T}$ on it compatible with such calculus, we construct a 
2-category $\twocat{L}$ with a 2-monad $\mathsf{S}$ on it such that:
\begin{itemize}
\item $\mathsf{S}$ has the adjoint-pseudo-algebra property.
\item The 2-categories of pseudo-algebras of  $\mathsf{S}$ and 
$\mathsf{T}$ are equivalent.
\end{itemize}
Thus, coherent structures (pseudo-$\mathsf{T}$-algebras) are transformed 
into universally characterised ones 
(adjoint-pseudo-$\mathsf{S}$-algebras). The 2-category $\twocat{L}$ 
consists of lax algebras for the pseudo-monad induced by $\mathsf{T}$ 
on the bicategory of bimodules of $\twocat{K}$. We give an intrinsic 
characterisation of pseudo-$\mathsf{S}$-algebras in terms of {\em 
representability\/}. Two major consequences of the above transformation 
are the classifications of lax and strong morphisms, with the attendant 
coherence result for pseudo-algebras. We apply the theory in the context 
of internal categories and examine monoidal and monoidal globular 
categories (including their {\em monoid classifiers\/}) as well as 
pseudo-functors into $\Cat$.
\end{abstract}

\tableofcontents

\section{Introduction}
\label{sec:introduction}

In the categorical approach to 
 algebraic structures we deal with them in terms of algebras for a monad.
 In the study of structure 
 borne by a category we are led to consider 2-monads on $\Cat$, the 
 2-category of categories, functors and natural transformations 
 \cf.\cite{BlackwellKellyPower89}. However, the strict associativity 
 axiom for algebras is too restrictive to deal with the structures of 
 interest, \eg. completeness or cocompleteness, in which the operations
 are associative only up to isomorphism. This prompts the consideration 
 of {\em pseudo-algebras\/} for a 2-monad $\mathsf{T}$, where the usual 
 $\mathsf{T}$-algebra axioms 
 for a monad are weakened by the introduction of {\em structural 
 isomorphisms\/}, subject to 
 {\em coherence conditions\/}. 
 Likewise, the morphisms of interest are the {\em strong\/} ones (or 
 pseudo-morphisms) which preserve the operations only up to coherent 
 isomorphism.
 When the structure of interest is actually 
 determined by {\em universal properties\/}, \eg. a category $\cat{C}$ 
 possesing certain limits or colimits, the structural isomorphisms in the 
 pseudo-algebra are uniquely determined by such properties (they are {\em 
 canonical\/}) and their coherence conditions are automatically satisfied. 
 In other situations, such as $\cat{C}$ having a monoidal structure
 $I,\tensor,\alpha,\rho,\lambda$, there is nothing universal (in general) 
 about the operations ($I,\tensor$) {\em with respect to \/} $\cat{C}$ to 
 imply the existence of the structural isomorphisms (associativity 
 $\alpha$, left unit $\lambda$ and right unit $\rho$), and therefore they 
 must be provided as part of the data, and their coherence conditions 
 verified. However, if we regard such monoidal structure as algebraic 
 structure on a multicategory $\cat{M}$ (with underlying category 
 $\cat{C}$), it becomes a universal property, namely that $\cat{M}$ be 
 {\em representable\/} in the sense of \cite{Hermida99a}.
 
 The developments in \cite{Hermida99a} lead us to seek a systematic way of achieving 
 such a transformation. That is, given the 2-monad $\mathsf{T}$ on $\Cat$ 
 whose pseudo-algebras are monoidal categories (thus $T\cat{C} = 
 {\mathrm free monoid on }\cat{C}$), we would like to derive 
 the 2-category $\Multicat$ 
 of multicategories and the 2-monad $\mathsf{T}'$ on it whose 
 pseudo-algebras are again monoidal categories, where furthermore 
 $\mathsf{T}'$ has the {\em adjoint-pseudo-algebra property\/}. This 
 latter means that $\morph{x}{T'X}{X}$ bears a pseudo-algebra structure 
 if and only if $x\dashv \eta_X$, \ie. the structure is left adjoint to 
 the unit. This property has been analysed in \cite{Kock95}, and amounts 
 to requiring $\mu\dashv \eta{T}'$, \ie. the adjoint condition holds for 
 the free algebra. Thus, for a 
 given $X$ there is up to isomorphism only one possible pseudo-algebra 
 structure on it, if any such exists, which is completely characterised 
 by the adjunction condition.
 
 To formulate the question of `transforming  coherent structures into 
 universal properties' formally, we make the following  identifications:
 \vskip5mm
 \begin{center}
\fbox{coherent structure $\equiv$ pseudo-algebra for a 2-monad}
\end{center}
\vskip5mm
\begin{center}
\fbox{\begin{minipage}{0.9\textwidth}
structure with universal property $\equiv$ (adjoint-)pseudo-algebra for a 
2-monad with the adjoint-pseudo-algebra property.
\end{minipage}}
\end{center}
\vskip5mm

The problem can be now precisely stated as follows: given a 2-monad 
$\mathsf{T}$ on a 2-category $\twocat{K}$, can we find a 2-category 
$\twocat{K}'$ and a 2-monad $\mathsf{T}'$ on it such that 

\begin{enumerate}
        \item  \label{it:adj-ps-alg} $\mathsf{T}'$ has the adjoint-pseudo-algebra 
        property.
        
        \item  \label{it:ps-alg} The 2-categories of pseudo-algebras, strong 
        morphisms and transformations of $\mathsf{T}$ and $\mathsf{T}'$ are 
        equivalent ?
\end{enumerate}

The main goal of this paper is to give a positive answer to this 
question. We do so under the following hypotheses:

\begin{itemize}
        \item  $\twocat{K}$ admits a calculus of bimodules (\ie. pullbacks of 
        fibrations and cofibrations, pullback stable coidentifiers and Kleisli 
        objects for monads on bimodules, \cf. \sectref{sec:preliminaries}).
        
        \item  $\mathsf{T}$ is cartesian, preserves bimodules, their composites 
        and identities, and their Kleisli objects, \cf. \sectref{sec:bimod-T}.
\end{itemize}

Under these assumptioms the 2-monad $\mathsf{T}$ induces a pseudo-monad 
$\Bimod{\mathsf{T}}$ on $\Bimod{\twocat{K}}$. We show that taking 
$\twocat{K}' = \lax{\Alg{\Bimod{\mathsf{T}}}}$, the 2-category of normal lax 
algebras for this pseudo-monad (with representable bimodules as 
morphisms) provides another basis for the axiomatisation of 
$\Alg{\mathsf{T}}$, \ie. the 2-category of $\mathsf{T}$-algebras is 
monadic over $\twocat{K}'$. The 2-monad $\mathsf{T}'$ induced on 
$\twocat{K}'$ by this adjunction satisfies (\ref{it:adj-ps-alg}) and 
(\ref{it:ps-alg}) above. 

The question we posed above is of a foundational nature, as vehemently 
argued in \cite{Benabou85}. For a structure characterised by a universal 
property, it is its mere {\em existence\/} which matters, regardless of any 
actual {\em choice of representative\/} for the operations. Thus our result provides an 
alternative approach to `avoiding the axiom of choice in category theory' 
to that given in \cite{Makkai96}, which considers all possible choices of 
representatives simultaneously.

We will further establish two important consequences from our 
construction of $\twocat{K}'$ and $\mathsf{T}'$:

\begin{itemize}
        \item  \underline{\textbf{Classification of strong morphisms}:} under a 
        mild additional hypothesis on $\twocat{K}$ and $\mathsf{T}$ (see \sectref
{sec:strict-morph})
        we can reflect pseudo-algebras and strong morphisms into the strict 
        ones. The associated `strong morphism classifier' is a strict 
        $\mathsf{T}$-algebra {\em equivalent\/} to the given pseudo-algebra. 
        This is therefore a strong coherence result.
        
        \item  \underline{\textbf{Classification of lax morphisms}:} the monadic 
        adjunction $\Alg{\mathsf{T}}\rightarrow\twocat{K}'$ induces a 2-comonad
        on $\Alg{\mathsf{T}}$ whose Kleisli 2-category corresponds to that of 
        $\mathsf{T}$-algebras and {\em lax morphisms\/} between them, \cf. \sectref{sec:lax-morph}. 
We thus 
        get an effective and simple description of the classifier of lax 
        morphisms out of a (pseudo-)$\mathsf{T}$-algebra, and we recover some 
        important monoid classifiers, \cf. \sectref{sec:monoid-classifier} and \sectref{sec:globular-monoids}.
\end{itemize}

Our consideration of bimodules was motivated 
by the pursuit of the theory of representable multicategories in 
\cite{Hermida99a}. In trying to understand how to deal with the `hom' of a 
multicategory $\cat{M}$, we observed that the hom-sets of (multi)arrows 
$\cat{M}(\vec{x},y)$  could be organised into a bimodule 
$\rel{M}{T\overline{\cat{M}}}{\overline{\cat{M}}}$, where 
$\overline{\cat{M}}$ is the category of linear morphisms (those with a 
singleton domain) of $\cat{M}$ and $T\overline{\cat{M}}$ is the free 
strict monoidal category on it. Furthermore, the identities and 
composition of $\cat{M}$ endow the bimodule $M$ with the structure of a 
normal $\lax{\Bimod{\mathsf{T}}}$-algebra (or equivalently, with a monad 
structure as we will see in \sectref{sec:lax-algebras-monads}), and $\cat{M}$ can be recovered 
from such data. Thus we arrive at  three alternative views of  the 
structure `monoidal 
category':

\begin{enumerate}
        \item  \label{it:mon1} as a pseudo-algebra for the free-monoid monad on 
        $\Cat$.
        
        \item  \label{it:mon2} as an adjoint pseudo-algebra on a bimodule
         $\rel{M}{T\overline{\cat{M}}}{\overline{\cat{M}}}$ with a monad 
         structure.                                                                                                             
                
        \item    \label{it:mon3} as a representable multicategory, where a 
        multicategory is a monoid structure on a multigraph.
                        
\end{enumerate}

All three views are valuable. The first is the 
traditional one, which by virtue of the classical `all diagrams commute' 
result of Mac Lane \cite{MacLane71} admits a {\em finite 
presentation\/}. The other two approaches have the advantage of enabling 
us to reason by universal arguments and dispense with the coherence 
axioms. The third approach is technically the simplest and we will 
exploit it in \ref{sec:monoid-classifier} to give an explicit description of 
the monoid classifier. The second approach is the one that leads to the 
general theory we develop in the present paper and exhibits the 
correspondence with pseudo-algebras in a more penetrating way than the 
account based on multigraphs. The fact that the bimodule so constructed 
from a multicategory $\cat{M}$ corresponds to its `hom' is made precise 
in \cite{Hermida99b}, which exhibits it as part of a Yoneda structure on 
multicategories, and develops a theory for fibrations for them so that 
the Yoneda object on $\cat{M}$ does correspond to the object of discrete 
fibrations on it. Without indulging in details here, let us
point out that this provides a natural example of a 2-category where the 
appropriate notion of fibration is {\em not\/} the representable one 
advocated in \cite{Street73}.

In fact the second and third approaches above can be formally related 
(and shown equivalent) in the context of internal category theory (see 
\sectref{sec:internal-bimodules}). The second part of the paper is 
devoted to three basic examples in this context: the one already 
mentioned of monoidal categories, monoidal globular categories and 
pseudo-functors into $\Cat$. We show their coherence 
results  (via a technique introduced in \cite{Hermida99a}specifically for 
adjoint-pseudo-algebras and developed here in 
\sectref{sec:strict-morph})
and in the case of monoidal and monoidal globular categories
give explicit descriptions of their monoid classifiers.

\vspace{1ex}
\noindent \underline{\textbf{Overview of the paper}:}
Part \ref{part:gral-theory} deals with the general construction producing a 2-monad with the adjoint pseudo-algebra property from a given one, while Part \ref{part:applications} instantiates this general framework in the context of internal category theory and studies three basic examples.

\noindent \underline{Part \ref{part:gral-theory}:}
In \sectref
{sec:preliminaries} we recall the basic results about bimodules (2-sided discrete fibrations) in a 2-category which we need. In \sectref{sec:monads-Kleisli} we review the explicit description of Kleisli objects for monads on bimodules in $\Cat$, while in \sectref{sub:Kleisli-coinserter} we give a construction of the kind of Kleisli objects we need (which we christened {\em representable\/} since we require universality only in the context of representable bimodules) in terms of coinserters and coequifiers.

In \sectref{sec:bimod-T} we introduce the kind of 2-monad we deal with, namely {\em cartesian\/} ones in which the 2-functor preserves comma-objects and coidentifiers. We derive some intrinsic properties concerning the behaviour of such monad in the context of bimodules (Lemma \ref{lemma:cart-fib}) which allows us to deduce an induced pseudo-monad on the bicategory of bimodules (Corollary \ref{cor:pseudo-monad}). 

In \sectref{sec:lax-algebras-monads} we introduce a 2-category of lax algebras relative to the pseudo-monad previously constructed (Definition \ref{def:LaxBimodT}), the particular point of interest being the definition of 2-cells. We provide an alternative characterisation of this 2-category (Definition \ref{def:mnd-bimod} and Proposition \ref{prop:laxalg-mnd}) as the 2-category of monads in a Kleisli bicategory of bimodules and give an algebraic description of its 2-cells in terms of equivariant morphisms.

\sectref{sec:pseudo-alg-prop} is the technical core of the paper. In \sectref{sub:adjunction} we set up the fundamental 2-adjunction between $\Alg{\mathsf{T}}$ and $\lax{\Alg{\Bimod{T}}}$ (Proposition \ref{prop:fund-adj}), the main point of note being the use of representable Kleisli objects to obtain a free $\mathsf{T}$-algebra
on a lax $\Bimod{T}$-algebra. In \sectref{sec:adjoint-representable} we show that the 2-monad induced by this adjunction has the adjoint-pseudo-algebra property (Proposition \ref{prop:Kock-prop}), the first important intrinsic result of the present work, whose proof involves some technical delicacy. We then establish one of our main results (Theorem \ref{thm:representable-characterisation}), viz. the characterisation of lax $\Bimod{T}$-algebras which bear an adjoint-pseudo-algebra structure in terms of {\em representability} of their underlying bimodules and the invertibility of their structural associator 2-cell. In \sectref{sec:monadicity} we complete our basic general theory establishing the monadicity of the fundamental 2-adjunction (Theorem \ref{thm:monadicity}) and the correspondence of its psuedo-algebras with pseudo-$\mathsf{T}$-algebras. The proof of the theorem relies in the representable characterisation of adjoint-pseudo-algebras mentioned above.

In \sectref{sec:lax-morph} we exhibit an interesting byproduct of the setup in \sectref{sec:pseudo-alg-prop}, namely an explicit construction of the classifier of lax morphisms between (pseudo-)$\mathsf{T}$-algebras (Theorem \ref{thm:classification-lax-morphisms}).

\sectref{sec:strict-morph} deals with the classification of strong morphisms in terms of strict ones. In order to do so, we establish a couple of technical results regarding Kleisli objects (Propositions \ref{prop:Kleisli-lax} and \ref{prop:Kleisli-K}) as well as recalling a key technical lemma (Lemma \ref{key-lemma}) which shows that we can construct the relevant coinverter (for the classification of strong morphisms) using Kleisli objects. Theorem \ref{thm:classification-strong-morphisms} is our final main result of the general theory, showing that the classification sought yields the coherence result, namely every pseudo-algebra is equivalent to a strict one (Corollary
\ref{cor:coherence}). Notice that we provide a full-fledged coherence statement, encompassing morphisms and 2-cells.

\vspace{1ex}
\noindent \underline{Part \ref{part:applications}:}
In \sectref{sec:internal-bimodules} we show that given a category with pullbacks $\cat{B}$ and a cartesian
monad $\mathsf{T}$ on it, the 2-category $\Cat({\cat{B}})$ of internal categories and the induced 2-monad $\Cat(\mathsf{T})$ satisfy the hypothesis of our general theory. This requires a review of internal bimodules, in particular the construction of Kleisli objects (in \sectref{sub:Kleisli-int-bim}). The main original result 
here is Theorem \ref{thm:bimod-spn} which shows that in this context we can dispense with bimodules in favour of spans and their simple pullback composition. This in turn enables an easy explicit description of monoid classifiers for monoidal and monoidal globular categories in the subsequent sections.

\sectref{sec:multicat} reviews the basic results of \cite{Hermida99a} as consequences of the theory in Part
\ref{part:gral-theory}. The main novelty is the recovery of the well-known monoid classifier for monoidal categories (\ie. the category of finite ordinals and monotone maps) from our general classification of lax morphisms (Theorem \ref{thm:delta}). \sectref{sec:monoidal-globular-categories} shows that Batanin's monoidal globular categories fit into our framework, since they are pseudo-algebras for a 2-monad induced by a cartesian monad on the category of $\omega$-graphs (\sectref{sec:pseudo-algebras-globular-categories}). Hence we can give a universal characterisation of monoidal-globular categories in terms of  {\em representable multicategories\/} on $\omega$-graphs (Corollary \ref{cor:rep-w-multicat}),
recover their basic coherence result (Corollary \ref{coherence-monoidal-glob}), and their globular monoid classifier (Corollary \ref{cor:globular-monoid-classifier}, the basic result of \cite{BataninStreet98}). 

\sectref{sec:pseudo-functors} shows how the other basic example of a coherent structure, namely pseudo-functors into $\Cat$, are dealt with in our present theory. An interesting aspect of this example is that it exhibits clearly the representable characterisation of adjoint-pseudo-algebras (Remark \ref{rmk:representable-pseudo-functors}). It also shows how lax functors into $\Bimod{\Cat}$ arise naturally in this context, which explains their relevance for fibred category theory (Remark \ref{rmk:fibrations}).

\part{General Theory}
\label{part:gral-theory}
\section{Preliminaries on bimodules}
\label{sec:preliminaries}

Since there is no comprehensive account of the elementary properties of 
bicategories of bimodules internal to a 2-category, we collect some basic 
facts in this section.
Background material on bimodules can be found in 
\cite{Street73,Street80,Wood82,Wood85,CarboniJohnsonStreetVerity94} and 
\cite[\S 1]{BournCordier80}.
 In order to build a bicategory of bimodules internal 
to a given 2-category $\twocat{K}$, 
we assume the following:

\begin{enumerate}
\item $\twocat{K}$ admits pullbacks of fibrations and cofibrations

\item $\twocat{K}$ admits comma-objects

\item $\twocat{K}$ admits coidentifiers, stable under pullback along fibrations.
\end{enumerate}

\begin{remark}
  If $\twocat{K}$ admits pullbacks, it also admits comma-objects if and 
only if it admits cotensors with the arrow category
$^{\rightarrow}$.
\end{remark}

It is convenient to view a bicategory of bimodules as a 2-dimensional 
version of a bicategory of relations. From this perspective, our 
assumptions on \twocat{K} amount to {\em regularity\/} for an ordinary 
category. With \twocat{K} as above, recall that a \textbf{bimodule} or 
\textbf{discrete fibration} from $X$ to $Y$ (objects of $\twocat{K}$) is 
a span 

               \[\xymatrixrowsep{1pc}
        \xymatrixcolsep{1pc}
        \begin{diagram}
        &\dlto_{d_R}{R}\drto^{c_R}& \\
 {X}& &{Y}
               \end{diagram}
        \]
\noindent which is discrete as an object in $\Spn{\twocat{K}}(X,Y)$, \ie. 
has discrete fibres, and bears an algebra structure for the monad

\[ \morph{\Hom{X}\circ (\_ )\circ 
\Hom{Y}}{\Spn{\twocat{K}}(X,Y)}{\Spn{\twocat{K}}(X,Y)}\]      
        
 \noindent where $\_\circ\_$ is composition of spans (via pullbacks) and 
 the following diagram is a comma-object:
 
              \[\xymatrixrowsep{1pc}
        \xymatrixcolsep{1pc}
        \begin{diagram}
        &\dlto_{d}{\Hom{X}}\drto^{c}& \\
 {X}\drto_{\mathit{id}}&\Rightarrow &\dlto^{\mathit{id}}{Y}\\
& {Y}&
               \end{diagram}
        \]
\noindent Therefore $\Hom{X}$ amounts to the cotensor $X^{\rightarrow}$. 
Such an algebra structure on the span $R$ means that $\morph{d_R}{R}{X}$ 
is a split fibration, 
$\morph{c_R}{R}{Y}$ is  a split cofibration and their actions are 
compatible. Hence, such an algebra structure is unique up to isomorphism 
if it exists.
    
\begin{remark}[fibrations in $\Cat$] The internal definition of fibrations in a 2-category, as introduced in \cite{Street73}, reflects the situation in $\Cat$: given a functor
$\morph{p}{\cat{E}}{\cat{B}}$, the free fibration over it (in $\Cat/\cat{B}$) is given by the projection $\morph{\bar{p}}{\CommaObj{\mathit{id}}{p}}{\cat{B}}$ out of the comma-object. This defines a 2-monad $\morph{\CommaObj{\mathit{id}}{\_}}{\Cat/\cat{B}}{\Cat/\cat{B}}$ with the adjoint-pseudo-algebra property (\cf.\sectref{sec:introduction}). Thus, $p$ is a fibration iff the unit $\morph{\eta_p}{p}{\CommaObj{\mathit{id}}{p}}$ has a right-adjoint. This latter amounts to a choice of \textit{cleavage\/} for $p$. This characterisation of fibration (and its simple extension to deal with 2-sided discrete fibrations or bimodules) can be internalised in any 2-category with comma-objects, and it is this internal version we work with throughout.
  
\end{remark}

Now
we can define the bicategory of 
bimodules in such a 2-category \twocat{K}:

\begin{definition}
\label{def:bimod}
The 
{\bfseries 
        bicategory of bimodules} 
        $\Bimod{\twocat{K}}$ consists of 
        \begin{description}
                \item[objects]  those of $\twocat{K}$
                
                \item[morphisms]  a morphism from $X$ to $Y$ is a 
                bimodule from $X$ to $Y$, which we write $\rel{R}{X}{Y}$.
              
                  \item[2-cells]  a 2-cell between morphisms is a morphism 
between the top objects of the spans, commuting with the domain and 
codomain morphisms and the actions of the fibrations and cofibrations:
           \[\xymatrixrowsep{1pc}
        \xymatrixcolsep{1pc}
        \begin{diagram}
        &\dlto_{d_R}{R}\ddto^{\alpha}
        \drto^{c_R}& \\
 {X}& &{Y} \\
 &\ulto^{d_S}{S}\urto_{c_S}
               \end{diagram}
        \]

\noindent In short, $\alpha$ is a $\Hom{X}\circ (\_ )\circ 
\Hom{Y}$-algebra morphism.
        \end{description}
        
        The identity span on $X$ is 
$\rel{\CommaObj{\mathit{id}}{\mathit{id}} = \Hom{X}}{X}{X}$. Composition 
is 
given by
        
$$\prooftree
{\xymatrixrowsep{1pc}
 \xymatrixcolsep{1pc}
        \begin{diagram}
        &\dlto_{d_R}{R}\drto^{c_R}& &   &\dlto_{d_S}{S}\drto^{c_S}& \\
 {X}& &{Y}                         &  {Y}& &{Z} 
               \end{diagram}
}
\justifies
{\xymatrixrowsep{1pc}
 \xymatrixcolsep{1pc}
\begin{diagram}
&\dlto {R{\bullet}S}\drto &\\
{X}& &{Z} \\
\end{diagram}
}
\thickness=0.08em
\endprooftree $$        
        
\noindent where $R{\bullet}S$ is given by the the following 
        coequalizer
        
        \[\xymatrixrowsep{2pc}
 \xymatrixcolsep{2pc}
\begin{diagram}
 {R\circ \Hom{Y}\circ S}\ar@<1ex>[r]^-{l\circ S}
                                 \ar@<-1ex>[r]_-{R\circ{r}}&{R\circ 
                                 S}\ar@{>>}[r]& {R{\bullet} S}
\end{diagram}
\]

\noindent with $l$ being the action of the cofibration $\morph{c_R}{R}{Y}$
and $r$ being the action of the fibration $\morph{d_S}{S}{Y}$. 
 Horizontal composition of 
2-cells is canonically induced by that of morphisms, while their 
vertical 
composition is inherited from that of 1-cells in $\twocat{K}$. When drawing diagrams, we display bimodules by bent arrows, reserving straight arrows for 1-cells in $\twocat{K}$.
        
\end{definition}

We recall the following properties of the bicategory $\Bimod{\twocat{K}}$:

\begin{enumerate}

\item A morphism $\morph{f}{X}{Y}$ in $\twocat{K}$ gives rise to 
the \textbf{representable bimodule} $\rel{f_{\#}}{X}{Y}$
defined by 
the comma-object 
\[\xymatrixrowsep{1pc}
        \xymatrixcolsep{1pc}
        \begin{diagram}
        &\dlto_{d}{\im{f}}\drto^{c}& \\
 {X}\drto_{f}&{\Rightarrow} &\dlto^{\mathit{id}}{Y}\\
& {Y}&
               \end{diagram}
        \]
and to $\rel{f^{*}}{Y}{X}$ given by the comma-object
\[\xymatrixrowsep{1pc}
        \xymatrixcolsep{1pc}
        \begin{diagram}
        &\dlto_{d}{f^{*}}\drto^{c}& \\
 {Y}\drto_{\mathit{id}}&{\Rightarrow} &\dlto^{f}{X}\\
& {Y}&
               \end{diagram}
        \]
\noindent Furthermore, $\im{f}\dashv f^{*}$.
\item \underline{\textbf{Embeddings}:} There is a homomorphism 
$\morph{(\_)_{\#}}{\twocat{K}^{\mathit{co}}}{\Bimod{\twocat{K}}}$:
which sends a morphism $\morph{f}{X}{Y}$ to the representable bimodule 
$\rel{f_{\#}}{X}{Y}$, and a homomorphism 
$\morph{(\_)^{*}}{\twocat{K}^{\mathit{op}}}{\Bimod{\twocat{K}}}$ with 
action
$(\morph{f}{X}{Y}) \mapsto (\rel{f^{*}}{Y}{X})$.
Both these homomorphisms are locally fully faithful (that 
is, full and faithful on 2-cells).

\item \label{comma-object}
There is a particular instance of composition which deserves special 
interest: given morphisms $\morph{f}{X}{Z}$ and 
$\morph{g}{Y}{Z}$, the composite bimodule 
$\rel{\im{f}{\bullet}g^{*}}{X}{Y}$ is given by the top span of the 
comma-object
\[\xymatrixrowsep{1pc}
        \xymatrixcolsep{1pc}
        \begin{diagram}
        &\dlto_{p}{\CommaObj{f}{g}}\drto^{q}& 
\\ {X}\drto_{f}&{\Rightarrow} &\dlto^{g}{Y}\\
& {Z}&
               \end{diagram}
        \]
\noindent meaning that the canonical comparison 
$\CommaObj{f}{g}\rightarrow \im{f}{\bullet}g^{*}$ is an isomorphism. We 
have therefore the following exactness property:

\begin{displaymath}
  p^{*}\bullet\im{q}\iso \im{f}{\bullet}g^{*}\quad \mathrm{canonically}
\end{displaymath}

\item \label{it:fullfaith} It follows from the previous item that, given a 
morphism $\morph{f}{X}{Y}$, the unit 
$\cell{\tilde{f}}{\Hom{X}}{\im{f}{\bullet}f^{*}}$
of the adjunction $\im{f}\dashv f^{*}$ is an isomorphism iff $f$ is 
(representably) fully faithful.

\item Every bimodule 

              \[\xymatrixrowsep{1pc}
        \xymatrixcolsep{1pc}
        \begin{diagram}
        &\dlto_{d_R}{R}\drto^{c_R}& \\
 {X}& &{Y}
\end{diagram}
        \]
\noindent has a canonical factorisation $R \iso 
d_R^{*}{\bullet}\im{(c_R)}$. This factorisation is compatible with 
2-cells, in the sense that given 
        \[\xymatrixrowsep{1pc}
        \xymatrixcolsep{1pc}
        \begin{diagram}
        &\dlto_{d_R}{R}\ddto^{\alpha}
        \drto^{c_R}& \\
 {X}& &{Y} \\
 &\ulto^{d_S}{S}\urto_{c_S}&
               \end{diagram}
        \]
\noindent it induces $\morph{\alpha^{*}}{d_R^{*}}{d_S^{*}}$
and $\morph{\im{\alpha}}{\im{(c_R)}}{\im{(c_S)}}$ by the universal 
property of comma-objects, and then 
     \[\xymatrixrowsep{1.5pc}
        \xymatrixcolsep{1.5pc}
        \begin{diagram}
       {R}\dto_{\alpha}\rto^-{\sim}&{d_R^{*}\bullet\im{(c_R)}}
       \dto^{\alpha^{*}\bullet\im{\alpha}}\\
       {S}\rto^-{\sim}&{d_S^{*}\bullet\im{(c_S)}}
               \end{diagram}
        \]

\item\label{coidentifier-composition}
 Let us analyse what the composition of a representable bimodule and a 
dual of such amount to:

\[\xymatrixrowsep{2pc}
        \xymatrixcolsep{1.5pc}
        \begin{diagram}
&&{f^{*}\circ\im{g}}\dlto\drto\ar@{>>}[rrrr]^-{q}&&&&
{\Pi^{0}_{Y,Z}(f^{*}\circ\im{g})}\\
        &\dlto {f^{*}}\drto & &\dlto {\im{g}}\drto & &\\
 {Y}\drto_{\mathit{id}}&{\Rightarrow} &\dlto^{f}{X}
\drto_{g}&{\Rightarrow} &\dlto^{\mathit{id}}{Z}&\\
& {Y}&&{Z} &&
               \end{diagram}
        \]

\noindent thus $f^{*}\bullet\im{g} = \Pi^{0}_{Y,Z}(f^{*}\circ\im{g})$ 
where 
the latter is the `connected components' (\ie. the reflection into a 
discrete object) 
of $f^{*}\circ\im{g}$ in $\Spn{\twocat{K}}(Y,Z)$. Thus, we only require 
such coidentifiers (rather than general coequalizers) for composition of 
bimodules 
$$\prooftree
{\xymatrixrowsep{1pc}
 \xymatrixcolsep{1pc}
        \begin{diagram}
        &\dlto_{d_R}{R}\drto^{c_R}& &   &\dlto_{d_S}{S}\drto^{c_S}& \\
 {X}& &{Y}                         &  {Y}& &{Z} 
               \end{diagram}
}
\justifies
{\xymatrixrowsep{1pc}
 \xymatrixcolsep{2pc}
\begin{diagram}
&{X}\ar@/^1ex/[r]^-{(d_R)^{*}}&
{R}\ar@/^1ex/[r]^-{\im{(c_R)}}&
{Y}\ar@/^1ex/[r]^-{(d_R)^{*}}&
{S}\ar@/^1ex/[r]^-{\im{(c_S)}}&{Z}\\
{=} 
&{X}\ar@/^1ex/[r]^-{(d_R)^{*}}&
{R}\ar@/^1ex/[r]^-{p^{*}}&
{\scriptstyle\CommaObj{c_R}{d_S}}\ar@/^1ex/[r]^-{\im{q}}&
{S}\ar@/^1ex/[r]^-{\im{(c_S)}}&
{Z} \\
{=}
&{X}\ar@/^1ex/[rr]^-{(d_R\circ{p})^{*}}&&
{\scriptstyle\CommaObj{c_R}{d_S}}\ar@/^1ex/[rr]^-{\im{(c_S\circ{q})}}&&
{Z}
\end{diagram}
}
\thickness=0.08em
\endprooftree $$

\noindent where $p$ and $q$ are the projections out of the comma-object 
$\CommaObj{c_R}{d_S}$
as in (\ref{comma-object}).

\item \label{bimod-refl-span}
There is a reflection 

  \begin{displaymath}
    \xymatrixcolsep{3pc}
\begin{diagram}
{\Bimod{\twocat{K}}(X,Y)}\ar@/_2ex/[r]^-{\perp}&
 \ar@/_2ex/[l]_{^{\dag}}{\Spn{\twocat{K}}(X,Y)}
\end{diagram}
  \end{displaymath}

\noindent given as follows:
\[\xymatrixrowsep{1pc}
 \xymatrixcolsep{1.5pc}
\left( 
  \begin{diagram}
    &\dlto_{f}S\drto^{g}&\\
X& &{Y}
  \end{diagram}\right)^{\dag} = 
\begin{diagram}
  X\ar@/^1ex/[r]^-{f^{*}}&{S}\ar@/^1ex/[r]^-{\im{g}}&{Y}
\end{diagram}
\]

\noindent Therefore this reflection sends  the 2-cells  of 
$\Spn{\twocat{K}}(X,Y)$, 

\begin{displaymath}
\xymatrixrowsep{1pc}
 \xymatrixcolsep{1pc}
  \begin{diagram}
       &\dlto_{f}{S}\ddtwocell^{b}_{a}{\alpha}\drto^{g}&\\
X& &{Y}\\
&\ulto^{h}{T}\urto_{k}&
  \end{diagram}
\end{displaymath}
\noindent to identities, $\alpha^{\dag} = \mathit{id}: a^{\dag} = 
b^{\dag}$.

\item \underline{\textbf{Duality}:}
Every bimodule $\rel{R}{X}{Y}$ has associated its 
\textbf{dual} $\rel{R^\lineardual}{Y}{X}$ given by $R^\lineardual =  
c_R^{*}{\bullet}\im{(d_R)}$. 

The assignment $R\mapsto R^{\lineardual}$ extends to an 
identity-on-objects biequivalence 
$\morph{(\_)^\lineardual}{\Bimod{\twocat{K}}}{\Bimod{\twocat{K}}^{\mathit{op}}}$.

To describe its action on 2-cells, given $\cell{\alpha}{R}{S}$, 

\[ \cell{\alpha^{\lineardual} = \alpha^{*}{\bullet}\im{\alpha}}
{c_R^{*}{\bullet}\im{(d_R)}}{c_S^{*}{\bullet}\im{(d_S)}} \]

\noindent where $\cell{\alpha^{*}}{c_R^{*}}{c_S^{*}}$ and 
$\cell{\im{\alpha}}{\im{(d_R)}}{\im{(d_S)}}$ are induced by universality 
of comma-objects, as above.

\item \underline{\textbf{Change of base}:} Given a bimodule 
$\rel{R}{X}{Y}$ and 
morphisms $\morph{f}{X'}{X}$
and $\morph{g}{Y'}{Y}$, we obtain a bimodule 
$\rel{\im{f}{\bullet}R{\bullet}g^{*}}{X'}{Y'}$. We thus have a 
\textbf{change 
of base} functor 
$\morph{(f,g)^{*}}{\Bimod{\twocat{K}}(X,Y)}{\Bimod{\twocat{K}}(X',Y')}$,
whose action can be  described more simply via pullbacks; the bimodule 
$(f,g)^{*}(R)$
is obtained as indicated in the following limit diagram:

\[
  {\xymatrixrowsep{1.5pc} \xymatrixcolsep{1.5pc}
    \let\objectstyle=\scriptstyle
        \begin{diagram}          
{X'}\dto_-{f}&\ar@{->}[l]_-{d'}{(f,g)^{*}(R)}\ar@{->}[d]\ar@{->}[r]^-{c'}&{Y'}\dto^-{g}\\
{X}&\lto^-{d}{R}\rto_-{c}&{Y}
               \end{diagram}
               }
\]

\noindent In $\twocat{K} = \Cat$, we have 
$(f,g)^{*}(R)(x',y') = 
R(fx',gy')$.
\end{enumerate}

\subsection{Monads and Kleisli objects}
\label{sec:monads-Kleisli}
In our intended application of bimodules, monads and their associated 
Kleisli objects in $\Bimod{\twocat{K}}$ play a central role. As a helpful 
analogy, recall that for a given commutative ring $R$, an $R$-algebra 
$A$ can be presented either as a monoid in $\Mod{R}$ or as a 
ring together with a ring homomorphism $R\rightarrow A$ (whose image lies in the center of $A$). 
This latter view 
amounts to taking the Kleisli object of the monad (in the bicategory of 
rings and bimodules) corresponding to the former presentation.

We recall the explicit description of Kleisli objects for monads on 
bimodules in $\Cat$.
Given a monad $(\rel{M}{X}{X},\eta,\mu)$ in  $\Bimod{\twocat{\Cat}}$, its 
Kleisli object is given by the category $\underline{M}$ with 

\begin{description}
\item[objects] those of $X$.

\item[morphisms] $\underline{M}(x,y) = M(x,y)$, the fibre of the bimodule 
$M$ over the objects $x, y$.

\item[identities] $\underline{\mathit{id}}_x = \eta_{x,x}(\mathit{id}_x) 
\in M(x,x)$

\item[composition] 
$$\prooftree
f \in M(x,y) \qquad g \in M(y,z)
\justifies
g{\circ}f = \mu_{x,z}([f,g]) \in M(x,z)
\endprooftree $$
\noindent where $[f,g]$ represents the equivalence class of the pair 
$(f,g)$ in $M{\bullet}M$.
\end{description}

\noindent The unit $\cell{\eta}{\Hom{X}}{M}$ induces a functor 
$\morph{J}{X}{\underline{M}}$, which is the identity on objects. There is 
a 2-cell
$\cell{\rho}{M{\bullet}\im{J}}{\im{J}}$ which sets up a universal lax 
cocone. Notice that we recover the monad $M$ as the composite 
$\im{J}{\bullet}J^{*} \iso \CommaObj{J}{J}$, which is the resolution of 
the monad obtained via its Kleisli object. Just as in the case of algebras 
mentioned above, to give a monad on $X$ is equivalent then to give a 
category $M$ and an identity-on-objects functor $\morph{J}{X}{M}$.

\begin{remark}
\label{rmk:Kleisli-representable}
        It follows from the above explicit description that given a lax 
cocone 
        $\rel{\cell{\gamma}{M{\bullet}L}{L}}{X}{Y}$ in $\Bimod{\Cat}$,
        if $\rel{L}{X}{Y}$ is representable  ($L = \im{l}$) so is the induced 
        mediating bimodule $\rel{\hat{L}}{\underline{M}}{Y}$ ($\hat{L} = \im{\hat{l}}$),
and this preservation of representability is 2-functorial in the sense that given a morphism
$\morph{k}{Y}{Z}$ in $\twocat{K}$, we have $\widehat{k{\circ}l} = \widehat{k}{\circ}\widehat{l}$ and similarly for 2-cells between such morphisms 
 \cf. 
        \cite[Axiom 5]{Wood85}.  We say that Kleisli objects in $\Bimod{\Cat}$ {\em preserve representability\/} 
\end{remark}

\subsection{Kleisli objects for bimodules via coinserters and coequifiers}
\label{sub:Kleisli-coinserter}
In order to set up the monadic adjunction in \sectref{sub:adjunction} below we must assume that Kleisli objects for monads in $\Bimod{\twocat{K}}$ behave like those in $\Bimod{\Cat}$. In fact, we only need 
the universal property of
Kleisli objects (\ie. initial lax cocone) with respect to lax cocones 
$\rel{\cell{\gamma}{M{\bullet}L}{L}}{X}{Y}$ with $L$ {\em representable}, with the induced mediating bimodules representable as well (as in Remark \ref{rmk:Kleisli-representable}). 
We then say that $\Bimod{\twocat{K}}$ \textbf{admits representable Kleisli objects}.

\begin{assumption}
\label{assumptionKleisliBimod}
  $\Bimod{\twocat{K}}$ admits representable Kleisli objects and
\linebreak
$\morph{\Bimod{T}}{\Bimod{\twocat{K}}}{\Bimod{\twocat{K}}}$ preserves them.
\end{assumption}

We give an explicit construction of such representable Kleisli objects in $\twocat{K}$ using colimits, so as to have sufficient conditions on $\twocat{K}$ and $T$ for the above assumption to hold. Before embarking on the abstract construction, it might be helpful for the reader to think of it as a 2-dimensional analogue of `taking a quotient by an equivalence relation', which in a regular category is achieved by taking the coequalizer of the pair of morphisms (\textit{domain},\textit{codomain}) of the relation. Here, we first add a 2-cell between such pair and then impose the `lax cocone' equations (thus coequalizing at the level of 2-cells, rather than 1-cells). Thus our present work may be seen as some foundational 2-dimensional algebra intrinsic to a `regular 2-category'. 
We rely on the notions of coinserter and coequifier in a 2-category, which can be found in 
the survey article \cite{Kelly89}.

\begin{proposition}\label{prop:Kleisli-coinserter}
  If $\twocat{K}$ admits coinserters and coequifiers, then   $\Bimod{\twocat{K}}$ admits representable Kleisli objects.
\end{proposition}
\begin{proof}
  Consider a monad $(\rel{M}{X}{X},\eta,\mu)$ in $\Bimod{\twocat{K}}$. Let $M[\lambda]$ be the (object part of) the coinserter of the pair $\morph{d,c}{M}{X}$:
  \begin{displaymath}
\xymatrixrowsep{1.5pc}
 \xymatrixcolsep{2pc}
    \begin{diagram}
      M\rrcompositemap<2>_{d}^{J}{\lambda}
\rrcompositemap<-2>_c^J
&&{M[\lambda]}
    \end{diagram}
  \end{displaymath}

\noindent Notice that the 2-cell $\lambda$ sets up a lax cocone over $M \simeq d^{*}{\bullet}\im{c}$ as follows:
\begin{displaymath}
  \prooftree{
\cell{\lambda}{Jd}{Jc}}
\Justifies
{
\[
{\cell{\im{\lambda}}{\im{c}{\bullet}\im{J}}{\im{d}{\bullet}\im{J}}}
\Justifies
{\cell{\rho = (\im{\check{\lambda}})}{\im{J}}{d^{*}{\bullet}\im{c}{\bullet}\im{J}}}
\]}
\thickness=0.08em
\endprooftree
\end{displaymath}

But this lax cocone need not satisfy the equations involving $\eta$ and $\mu$, so we must enforce them.
For the unit, pasting $\rho$ and $\eta$, we obtain the 2-cell 
$\cell{\rho\cdot\im{J}\eta}{\im{J}}{\im{J}}$, 
which by (contravariant) local full and faithfullness of $\im{\_}$, 
corresponds to a (unique) 2-cell in $\twocat{K}$, written $\cell{\kappa}{J}{J}$. Consider the coequifier of this 2-cell and the identity:
 \begin{displaymath}
\xymatrixrowsep{1.5pc}
 \xymatrixcolsep{1.5pc}
    \begin{diagram}
      X\ar@/^2ex/[rrr]^-{J}
\ar@/_2ex/[rrr]_-{J}\xtwocell[0,2]{}\omit{\kappa}
&&&
{M[\lambda]}\xtwocell[0,-2]{}\omit{^\mathit{id}}\rto^-{\pi}&{M[\lambda]_{/(\kappa\equiv\mathit{id})}}
    \end{diagram}
  \end{displaymath}

\noindent We obtain thus a new cocone $\morph{\pi{J}}{X}{M[\lambda]_{/(\kappa\equiv\mathit{id})}}$ with $\cell{\pi\rho}{M{\bullet}\im{\pi{J}}}{\im{\pi{J}}}$. Again\footnote{If $\twocat{K}$ admits sums, we could perform both `coequifications' in one step.}, we must impose on this data the condition for $\mu$. We have two 
2-cells $$\cell{\pi\rho\cdot\pi{\rho}M}{M{\bullet}M}{M} \qquad
\cell{\pi\rho\cdot\mu}{M{\bullet}M}{M}$$ Using the following comma square

\[\xymatrixrowsep{1pc}
        \xymatrixcolsep{1pc}
\let\labelstyle=\scriptstyle
        \begin{diagram}
        &\dlto_{p}{\CommaObj{c}{d}}\drto^{q}& \\
 {M}\drto_{c}&\Rightarrow &\dlto^{d}{M}\\
& {X}&
               \end{diagram}
        \]
\noindent the above 2-cells determine, by adjoint transposition, 2-cells
$$\cell{\check{(\pi\rho\cdot\pi{\rho}M)}}{\im{q}{\bullet}\im{c}{\bullet}\im{J}}{\im{p}{\bullet}\im{d}{\bullet}\im{J}} \qquad
\cell{\check{(\pi\rho\cdot\mu)}}{\im{q}{\bullet}\im{c}{\bullet}\im{J}}{\im{p}{\bullet}\im{d}{\bullet}\im{J}}$$

\noindent which by (contravariant) local full and faithfullness of $\im{\_}$, correspond to two parallel 2-cells $\cell{\mu_l,\mu_r}{Jdp}{Jcq}$, which we coequify:

\begin{displaymath}
\xymatrixrowsep{1.5pc}
 \xymatrixcolsep{1.5pc}
    \begin{diagram}
     {\CommaObj{c}{d}}\ar@/^2ex/[rrr]^-{Jdp}
\ar@/_2ex/[rrr]_-{Jcq}\xtwocell[0,2]{}\omit{\mu_l}
&&&
{M[\lambda]_{/(\kappa\equiv\mathit{id})}}\xtwocell[0,-2]{}\omit{^\mu_r}\rto^-{\varpi}&{M[\lambda]_{/(\kappa\equiv\mathit{id},\mu_l\equiv\mu_r)}}
    \end{diagram}
  \end{displaymath}
\noindent so that $\underline{M} = M[\lambda]_{/(\kappa\equiv\mathit{id},\mu_l\equiv\mu_r)}$ is the Kleisli object, with
universal lax cocone $\morph{\varpi\pi{J}}{X}{\underline{M}}$ and $\cell{\varpi\pi\rho}{M{\bullet}\im{\varpi\pi{J}}}{\varpi\pi{J}}$. Universality and preservation of representability follow at once. 

\qed
\end{proof}

\begin{corollary}\label{cor:Kleisli-coinserters}
  If $\morph{T}{\twocat{K}}{\twocat{K}}$ preserves coinserters and coequifiers,
then \linebreak
$\morph{\Bimod{T}}{\Bimod{\twocat{K}}}{\Bimod{\twocat{K}}}$ preserves representable Kleisli objects.
\qed
\end{corollary}

\begin{remarks}
\hfill{ }
\begin{itemize}
  \item
 Notice that when $\twocat{K}$ admits a Yoneda structure (in the sense of \cite{StreetWalters78}), 
so that a bimodule 
$\rel{M}{X}{Y}$ is classified by a 1-cell $\morph{\hat{M}}{Y}{PX}$ (\eg. for $\twocat{K} = \Cat$, the Yoneda object is $P\cat{X} = [\cat{X}^{\mathit{op}},\Set]$), this restricted universal property with respect to representables entails the general one, so that $\Bimod{\twocat{K}}$ admits Kleisli objects in full generality (and moreover, representability is preserved).

\item In \sectref{sub:Kleisli-int-bim} we show that for $\twocat{K} = \Cat(\cat{B})$, the 2-category of internal categories in $\cat{B}$, $\Bimod{\twocat{K}}$ admits Kleisli objects preserving representability, although it does not admit a Yoneda structure in general, \ie. internal bimodules are not classifiable by internal functors.

  \end{itemize}
\end{remarks}

\section{Bicategory of bimodules 
associated to a cartesian 2-monad}
\label{sec:bimod-T}

Given a 2-category $\twocat{K}$ admitting a calculus of bimodules as in 
\sectref{sec:preliminaries}, we consider a 2-monad $(T,\eta,\mu)$ on it 
which induces a pseudo-monad on $\Bimod{\twocat{K}}$. More precisely, we 
assume 

\begin{enumerate}
\item \label{it:one} $\morph{T}{\twocat{K}}{\twocat{K}}$ preserves 
pullbacks, comma-objects and the coidentifiers 
$f^{*}{\circ}\im{g}\rightarrow\Pi^{0}(f^{*}{\circ}\im{g})$ of 
\sectref{sec:preliminaries}.(\ref{coidentifier-composition}).
\item $\cell{\eta}{1}{T}$ and $\cell{\mu}{T^{2}}{T}$ are cartesian 
2-natural transformations, \ie. the naturality squares are pullbacks.
\end{enumerate}

Notice that by (\ref{it:one}) $T$ preserves bimodules, their identities 
and composites. It does therefore induce a homomorphism
\[\xymatrixrowsep{1pc}
        \xymatrixcolsep{1pc}
        \begin{diagram}        
{\Bimod{T}:}&{\Bimod{\twocat{K}}}\ar@{>}[rrr]&&&{\Bimod{\twocat{K}}}&\\
        &\dlto_{d_R}{R}\ddto^{\alpha}
        \drto^{c_R}& {\mapsto}& & \dlto_{Td_R}{TR}\ddto^{T\alpha}
        \drto^{Tc_R}&\\
 {X}& &{Y} &{TX}& &{TY} \\
 &\ulto^{d_S}{S}\urto_{c_S}& &&\ulto^{Td_S}{TS}\urto_{Tc_S}&
               \end{diagram}
        \]

We want such a 2-monad to induce both a pseudo-monad 
$(\Bimod{T},\im{\eta},\im{\mu})$ and a pseudo-comonad 
$(\Bimod{T},{\eta}^{*},{\mu}^{*})$ 
on $\Bimod{\twocat{K}}$ (see \eg.\cite{DayStreet97,Lack99,Marmolejo99} for 
details on 
pseudo-monads). Hence we must show that the cartesian 2-natural 
transformations $\eta$ and $\mu$ induce pseudo-natural transformations. 
This follows from the technical lemma below.

\begin{lemma}
\label{lemma:cart-fib}
        Let $\morph{F,G}{\twocat{L}}{\twocat{K}}$ be 2-functors,  with $G$ 
preserving cotensors with $^{\rightarrow}$,
 and
$\cell{\alpha}{F}{G}$ be a cartesian 2-natural transformation. The 
following properties  hold: 
        
        \begin{enumerate}
                \item  \label{it:fib-cofib}Every component 
$\morph{\alpha_x}{Fx}{Gx}$ is both 
a split fibration and 
        cofibration.
                
                \item \label{canonical-iso} For $\morph{f}{x}{y}$, the 
naturality square
\[\xymatrixrowsep{1pc}
        \xymatrixcolsep{1pc}
        \begin{diagram}
   {Fx}\dto_{\alpha_x}\rto^-{Ff}&{Fy}\dto^{\alpha_y}\\
{Gx}\rto_{Gf}&{Gy}
               \end{diagram}
        \]

induces an isomorphism 
$\isomorph{\theta_f}
{\im{(\alpha_x)}\bullet\im{Gf}}
{\im{Ff}\bullet\im{(\alpha_y)}}$. 

Its adjoint
transpose 
$\morph{\check{\theta}_f}{(Ff)^{*}\bullet\im{(\alpha_x)}}{\im{(\alpha_y)}\bullet(Gf)^{*}}$ 
is also an isomorphism.

        \end{enumerate}
        
\end{lemma}

\begin{proof}

  \begin{enumerate}
  \item Recall that for a given object $x$
        
         \[\xymatrixrowsep{1.5pc} \xymatrixcolsep{1.5pc}
         \let\objectstyle=\scriptstyle
        \begin{diagram}          
&\ar@/_2ex/[ddl]_-{\mathit{id}}{x}\dto|-{\tuple{\mathit{id}}}\ar@/^2ex/[ddr]^-{\mathit{id}}&\\
&\dlto_{d}{\Hom{x}}\drto^{c}&\\          
{x}\drto_{\mathit{id}}&{\stackrel{\lambda}{\Rightarrow}}&\dlto^{\mathit{id}}{x}\\          
&{x}&
               \end{diagram}
               \]
               
               we have the string of adjunctions

         \[\xymatrixrowsep{1.5pc} \xymatrixcolsep{2.5pc}
         \let\objectstyle=\scriptstyle
        \begin{diagram}
         { \Hom{x}}\ar@/^4ex/ [d]^-{c}
          \ar@/_4ex/ [d]_-{d}\\
          {x}\uto|{\tuple{\mathit{id}}}^{\dashv\;\;}_{\;\;\dashv}
               \end{diagram}
               \]

\noindent 
\comment{For $\morph{\alpha_x}{Fx}{Gx}$, cartesianness implies the 
following equality
         \[\xymatrixrowsep{1.5pc} \xymatrixcolsep{1.5pc}
         \let\objectstyle=\scriptstyle \let\labelstyle=\scriptstyle
        \begin{diagram}
          {Fx}\ddrto_-{\alpha_x}\ar@{->}[r]^-{\eta}&{\CommaObj{Gx}{\alpha_x}}\ddto\rto\xtwocell[2,1]{}\omit{^}
&
{Fx}\ddto^-{\alpha_x}\ar@{}[ddr]|{=}&
          {Fx}\dto_{\alpha_x}\ar@{->}[r]^{F\tuple{\mathit{id}}}&{F\Hom{x}}\dto|{\alpha_{\Hom{X}}}\rto^{c}&
          {Fx}\dto_{\alpha_x}\\
          &&&{Gx}\rto^-{G\tuple{\mathit{id}}}\drto_-{\mathit{id}}
          &{G\Hom{x}}\dto|-{Gd}\rto^-{Gc}
\xtwocell[1,1]{}\omit{^{G\lambda}}
&{Gx}\dto^-{\mathit{id}}\\
          &{Gx}\rto_-{\mathit{id}}&{Gx}&{\mbox{ }} &{Gx}\rto_-{\mathit{id}}&{Gx}
               \end{diagram}
               \]
               
               Hence $F\Hom{x}\simeq\CommaObj{Gx}{\alpha_x}$ and $F\tuple{\mathit{id}}\dashv Fc$ provides the right
               adjoint to $\eta = F\tuple{\mathit{id}}$ which gives
               $\alpha_x$ the structure of a split fibration. The same
               argument in $\twocat{K}^{\mathit{co}}$ yields the split
               cofibration structure.
}

For $\morph{\alpha_x}{Fx}{Gx}$, cartesianness implies that
\[
\xymatrixrowsep{1.5pc} \xymatrixcolsep{1.5pc}
         \let\objectstyle=\scriptstyle \let\labelstyle=\scriptstyle
        \begin{diagram}
{F\Hom{x}} \rto^-{Fc} \dto_-{\alpha_{\Hom{X}}} & 
  {Fx} \dto^-{\alpha_x} \\
{G\Hom{x}} \rto^-{Gc} \dto_-{Gd} \drtwocell<\omit>{^G\lambda} & {Gx}\dto^-{\mathit{id}} \\
{Gx} \rto_-{\mathit{id}} & {Gx}
\end{diagram}
\]
is the comma object $Gx\downarrow\alpha_x$. Thus since 
$\morph{F\tuple{\mathit{id}}}{Fx}{F\Hom{x}}$ is left adjoint to 
$\morph{Fc}{F\Hom{x}}{Fx}$, 
also $\morph{\eta}{Fx}{Gx\downarrow\alpha_x}$ is left adjoint to the projection $Gx\downarrow\alpha_x\rightarrow Fx$; and so $\alpha_x$ is a split fibration. 
The same argument in $\twocat{K}^{\mathit{co}}$ yields the split cofibration structure.

\item Recall that $\rel{\im{(\alpha_y)}\bullet(Gf)^{*}}{Fy}{Gx}$ is the 
comma-object

\[\xymatrixrowsep{1pc}
        \xymatrixcolsep{1pc}
        \begin{diagram}
        &\dlto_{p}{\scriptstyle\CommaObj{\alpha_y}{Gf}}\drto^{q}& \\
 {Fy}\drto_{\alpha_y}&\Rightarrow &\dlto^{Gf}{Gx}\\
& {Gy}&
               \end{diagram}
        \]

\noindent and therefore the naturality square induces a canonical morphism 
\linebreak
$\morph{\zeta}{Fx}{\CommaObj{\alpha_y}{Gf}}$
in ${\Spn{\twocat{K}}}(Fy,Gx)$. Since $\alpha_y$ is a split cofibration, the 
canonical morphism
$\morph{\eta}{Fy}{\CommaObj{\alpha_y}{Gy}}$ has a left adjoint $l\dashv 
\eta$ over $Gy$. Considering the following diagram
\[\xymatrixrowsep{1pc}
        \xymatrixcolsep{1pc}
       \let\objectstyle=\scriptstyle \let\labelstyle=\scriptstyle
        \begin{diagram}
{Fx}\rto^-{Ff}\dto_{\zeta}\ar@/_4ex/[dd]_-{\alpha_x}
&{Fy}\dto^{\eta}\ar@/^4ex/[dd]^-{\alpha_y}\\
{{\alpha_y}\downarrow{Gf}}\rto\dto&{{\alpha_y}\downarrow{Gy}}\dto \\
{Gx}\rto_{Gf}&{Gy}
               \end{diagram}
        \]
\noindent where both squares and the outer rectangle are pullbacks, we see 
that $\zeta$ is the pullback of $\eta$ along $Gf$, and therefore\footnote{Pullback along $Gf$ yields a 2-functor $\morph{Gf\ups}{\twocat{K}/Gy}{\twocat{K}/Gx}$ which takes the left adjoint $l$ to a left adjoint $Gf\ups{l}$ of $Gf\ups(\eta) = \zeta$} has a 
left adjoint $Gf\ups (l) \dashv \zeta$. This adjunction is taken by the 
reflection 
\sectref{sec:preliminaries}.(\ref{bimod-refl-span}) to an isomorphism 
$\isomorph{\zeta^{\dag}}{(Ff)^{*}\bullet\im{(\alpha_x)}}{\im{(\alpha_y)}\bullet(Gf)^{*}}$ 
which corresponds to 
${\check{\theta}_f}$.

\end{enumerate}

\qed
\end{proof}

\begin{corollary}
\label{cor:pseudo-monad}
  A cartesian 2-natural transformation $\cell{\alpha}{F}{G}$ as in Lemma 
\ref{lemma:cart-fib},
with both $\morph{F,G}{\twocat{L}}{\twocat{K}}$ preserving pullbacks, 
comma-objects and the relevant coidentifiers,
induces pseudo-natural transformations 
$$\cell{\Bimod{\im{\alpha}}}{\Bimod{F}}{\Bimod{G}}$$
\noindent and $$\cell{\Bimod{\alpha^{*}}}{\Bimod{G}}{\Bimod{F}}.$$
\end{corollary}
\begin{proof}
\mbox{ }\hfill{ }
  \begin{itemize}
  \item \underline{$\Bimod{\im{\alpha}}$:} Given a bimodule $\rel{M\iso
      d^{*}\bullet{\im{c}}}{X}{Y}$ define the invertible 2-cell
    $\cell{\alpha_M}{FM\bullet\im{(\alpha_y)}}{\im{(\alpha_x)}\bullet{GM}}$
    as the pasting composite

\[\xymatrixrowsep{1.5pc}
        \xymatrixcolsep{1.5pc}
       \let\objectstyle=\scriptstyle \let\labelstyle=\scriptstyle
        \begin{diagram}
{Fx}\ruppertwocell<1>^{(Fd)^{*}}{<3>{\check{\theta_f}}}\ar@/_1ex/[d]_-{\im{(\alpha_x)}}
&{FM}\ar@/^1ex/[d]|-{\im{(\alpha_M)}^{}}\ruppertwocell<1>^{\im{(Fc)}}{<3>{\alpha_c}}&{Fy}\ar@/^1ex/[d]^-{\im{(\alpha_y)}}\\
{Gx}\ar@/_1ex/[r]_-{(Gd)^{*}}&{GM}\ar@/_1ex/[r]_-{\im{(Gc)}}&{Gy}
             \end{diagram}
        \]

\noindent where the left isomorphism is that of Lemma 
\ref{lemma:cart-fib}.(\ref{canonical-iso}) and the right one is given by the 
homomorphism $\im{(\_)}$ applied to the naturality square for $\alpha$.

\item \underline{$\Bimod{\alpha^{*}}$:} Use the same argument 
in $\twocat{K}^{\mathit{co}}$.
\end{itemize}
\qed
\end{proof}

\begin{corollary}
        Given $(T,\eta,\mu)$ a cartesian 2-monad on $\twocat{K}$, with $T$ 
        preserving comma-objects,  it induces both a pseudo-monad
        $(\Bimod{T},\im{\eta},\im{\mu})$ and a pseudo-comonad 
$(\Bimod{T},{\eta}^{*},{\mu}^{*})$ 
on $\Bimod{\twocat{K}}$.
\end{corollary}

\begin{remark}[Coherence assumptions]
  \label{rmk:coherence-assumptions}
We make the following simplifying assumptions allowed by coherence for 
bicategories:
\begin{itemize}
\item We regard $\Bimod{\twocat{K}}$ as a 2-category, thereby ignoring the 
associativity constraints.

\item We regard 
$\morph{\Bimod{T}}{\Bimod{\twocat{K}}}{\Bimod{\twocat{K}}}$, 
$\morph{(\_)_{\#}}{\twocat{K}^{\mathit{co}}}{\Bimod{\twocat{K}}}$ and
$\morph{(\_)^{*}}{\twocat{K}^{\mathit{op}}}{\Bimod{\twocat{K}}}$ as 
2-functors, thereby ignoring the coherent structural isomorphisms for 
composites and identities.
\end{itemize}
\end{remark}

Now we can consider the Kleisli bicategory 
$\BimodT{\bold{T}}{\twocat{K}}$ associated to this pseudo-comonad, which 
we describe explicitly as follows:

  \begin{description}
                \item[objects]  those of $\twocat{K}$
                
                \item[morphisms]  a morphism from $X$ to $Y$ is a span
                
                \[\xymatrixrowsep{2pc}
        \xymatrixcolsep{1.5pc}
        \begin{diagram}
        &\dlto_{d_R}{R}\drto^{c_R}& \\
 {TX}& &{Y}
               \end{diagram}
        \]
                
\noindent which is a discrete fibration from $TX$ to $Y$, or more simply a 
1-cell in $\Bimod{\twocat{K}}$ between these objects.                
                
                \item[2-cells]  a 2-cell is a morphism of bimodules 
so that 
$$\BimodT{\bold{T}}{\twocat{K}} (X,Y) = 
               \Bimod{\twocat{K}} (TX,Y).$$
        \end{description}
        
\noindent The identity span on $X$ is $\rel{\eta_X^{*}}{TX}{X}$
 and composition is given by
        
$$\prooftree
{\xymatrixrowsep{2pc}
 \xymatrixcolsep{2pc}
        \begin{diagram}
{TX}\ar@/^1ex/ [r]^{R}&{Y}& &{TY}\ar@/^1ex/ [r]^{S}&{Z}
               \end{diagram}
}
\justifies
{\xymatrixrowsep{2pc}
 \xymatrixcolsep{2pc}
\begin{diagram}
{TX}\ar@/^1ex/ [r]^{\mu_X^{*}}&{T^{2}X} \ar@/^1ex/ [r]^{TR}&{TY}\ar@/^1ex/ 
[r]^{S}&{Z}
\end{diagram}
}
\thickness=0.08em
\endprooftree $$        
        
\noindent Horizontal composition of 
2-cells is clearly induced by that of morphisms, while the vertical 
composition is inherited from $\Bimod{\twocat{K}}$. 

\begin{remark}
  The bicategory $\BimodT{\bold{T}}{\twocat{K}}$ is analogous to our 
bicategory $\SpnT{\bold{T}}{\cat{B}}$ associated to a cartesian monad on  
a category with pullbacks $\cat{B}$ in
\cite[Appendix A]{Hermida99a}.
\end{remark}

We record the following  properties of the unit  of a cartesian 2-monad:

\begin{proposition}
\label{pr:unit}
  Consider a 2-monad $(T,\eta,\mu)$ on a 2-category $\twocat{K}$ admitting 
a calculus of bimodules 
  \begin{enumerate}
  \item \label{it:mono} If $\eta$ is cartesian, then it is monic.

  \item \label{it:ff} If $T$ preserves cotensors with $^{\rightarrow}$ and 
$\eta$ is cartesian, then $\eta$ is fully faithful.
  \end{enumerate}
\end{proposition}

\begin{proof}
  \begin{enumerate}
  \item The square
{\xymatrixrowsep{1.5pc}
 \xymatrixcolsep{1.5pc}
        \begin{diagram}
{X}\dto_{\eta_X}\rto^-{\eta_X}&{TX}\dto^{T\eta_X}\\
{TX}\rto_-{\eta_{TX}}&{T^{2}X}
               \end{diagram}
}
is a pullback and the lower arrow is a split monic ($\mu_X\eta_{TX} = 
\mathit{id}$), hence the top arrow is monic.

\item  Given `global generic' elements $\morph{a,b}{Z}{X}$, we must 
stablish a 1-1 correspondence 
  \begin{displaymath}
    (\cell{\alpha}{a}{b}) \leftrightarrow 
(\cell{\eta_X\alpha}{\eta_Xa}{\eta_Xb})
  \end{displaymath}
The 2-cells into $X$ are classified by 1-cells into $\Hom{X}$ and those 
into $TX$ are classified by 1-cells into $\Hom{TX} \iso T(\Hom{X})$.
 In the following diagram

   \[\xymatrixrowsep{1.5pc} \xymatrixcolsep{1.5pc}
     \let\objectstyle=\scriptstyle
        \begin{diagram}          
{X}\dto_-{\eta_X}&\lto_-{d}{\Hom{X}}\dto|-{\eta_{\Hom{X}}}\rto^-{c}&{X}\dto^-{T\eta_X}\\
{TX}&\lto^-{Td}{T\Hom{X}}\rto_-{Tc}&{TX}
               \end{diagram}
               \]

\noindent both squares are pullbacks and $\eta_{\Hom{X}}$ is monic, which
yields the desired correspondence.

  \end{enumerate}
\qed
\end{proof}

\section{Lax algebras vs. monads}
\label{sec:lax-algebras-monads} 
Given the pseudo-monad $(\Bimod{T},\im{\eta},\im{\mu})$ on 
$\Bimod{\twocat{K}}$, we consider lax algebras for it, as in \eg. 
\cite{Street73}.

\begin{definition}
  \label{def:lax-alg}
A \textbf{lax} $\Bimod{T}$-\textbf{algebra} consists of an object $X$ of 
$\twocat{K}$, a bimodule $\rel{M}{TX}{X}$ and structural 2-cells
$$\cell{\iota}{\Hom{X}}{\im{(\eta_X)}{\bullet}M}$$ \noindent and 
$$\cell{m}{TM{\bullet}M}{\im{(\mu_X)}{\bullet}M}$$ \noindent satisfying the following 
axioms:


\begin{itemize}

\item unit 

\[\xymatrixrowsep{2pc}
 \xymatrixcolsep{2pc}
\let\objectstyle=\scriptscriptstyle
\let\labelstyle=\scriptscriptstyle
        \begin{diagram}
&{TX}\ar@/_1ex/[dl]_-{\eta_{TX}^{}}
\ar@/_1ex/[d]|{{T\eta_X}^{}}\druppertwocell^{\mathit{id}}{<1>\iota}
& & & {TX}\ddtwocell^{M}_{M}{=}& &{TX}\ar@{}[dr]|{\cong^{}}
\ar@/_1ex/[d]^-{{\eta_{TX}}^{}}
\ar@/^1ex/[r]^-{M^{}} &{X}\ar@/_1ex/[d]_{\eta_X}\dduppertwocell<10>^{\mathit{id}}{<-3>\iota}\\
{T^2X}\ar@/_1ex/[dr]_-{\mu_X^{}}&{T^2X}\ar@/_1ex/[d]_{\mu_X}\ruppertwocell<1>_{TM}{<4>m}&{TX}\ar@/^1ex/[d]^{M}
&=& & =&
{T^2X}\ar@/_1ex/[d]_{\mu_X}\ruppertwocell<1>^{TM}{<4>m}&{TX}\ar@/_1ex/[d]^{M}\\
&{TX}\ar@/_1ex/[r]_{M}&{X}& & {X} & &{TX}\ar@/_1ex/[r]_{M}&{X}
               \end{diagram}
\]

\newpage
\item associativity 

\[\xymatrixrowsep{1.5pc}
 \xymatrixcolsep{1.5pc}
\let\objectstyle=\scriptscriptstyle
\let\labelstyle=\scriptscriptstyle
               \begin{diagram}
&                 {T^3X}\ar@/_1ex/[dl]_-{{\mu_{TX}}^{}}
\ar@/_1ex/[d]_{T\mu_X}\ruppertwocell<1>^{T^2M}{<3>^{Tm}}
&{T^2X}\ar@/^1ex/[d]^{TM}& &
{T^3X}\ar@{}[dr]|{\cong}\ar@/^1ex/[rr]^{T^2M}\ar@/_1ex/[d]_{\mu_{TX}}
&&{T^2X}\ar@/_1ex/[dl]_{\mu_X}
\ar@/^1ex/[d]^{TM}\xtwocell[2,-1]{}\omit{m}\\
{T^2X}\ar@/_1ex/[dr]_-{{\mu_X}^{}}
&{T^2X}\ar@/_1ex/[d]_{\mu_X}\ruppertwocell<1>_{TM}{<4>m}&{TX}\ar@/^1ex/[d]^{M}&=&
{T^2X}\ar@/_1ex/[d]_{\mu_X}\ruppertwocell<1>_{TM}{<4>m}&{TX}\ar@/^1ex/[d]^{M}&{TX}\ar@/^1ex/[dl]^{M}\\
&{TX}\ar@/_1ex/[r]_{M}&X&&{TX}\ar@/_1ex/[r]_{M}&X&
               \end{diagram}
\]

\end{itemize}

When no confusion is likely we simply say that $\rel{M}{TX}{X}$ is a lax 
algebra, understanding the rest of the data. The lax algebra is 
\textbf{normal} if $\iota = \mathit{id}_{\Hom{X}}$ (which of course 
requires $\rel{\Hom{X} = \im{(\eta_X)}{\bullet}M}{X}{X}$).
\end{definition}


Given lax algebras $\rel{M}{TX}{X}$ and $\rel{N}{TY}{Y}$, a morphism 
between them is given by a morphism $\morph{f}{X}{Y}$
 in $\twocat{K}$ and a 2-cell 
$\cell{\theta_f}{(Tf)^{*}{\bullet}M}{N{\bullet}f^{*}}$ in 
$\Bimod{\twocat{K}}(TY,X)$, compatible with the structural 2-cells.

\begin{remark}
  The notion of morphism between lax algebras above is a special instance 
of the notion of lax morphism of lax algebras in \cite{Street73}, where we 
require $f$ to be a representable bimodule (one coming from a morphism in 
$\twocat{K}$). Indeed, taking the adjoint mate of $\theta_f$ across 
$\im{Tf}\dashv(Tf)^{*}$ and 
 $\im{f}\dashv f^{*}$ we get 

\[\xymatrixrowsep{1.5pc}
 \xymatrixcolsep{1.5pc}
               \begin{diagram}
   {TX}
\ar@/_1ex/[d]_{\im{Tf}}
\ruppertwocell<1>^{M}{<3>{{\widehat{\theta_f}}}^{}}&{X}
\ar@/^1ex/[d]^{\im{f}}\\
{TY}\ar@/_1ex/[r]_-{N}&{Y}
\end{diagram}
\]
\noindent but the above presentation makes it easier to fit 2-cells. In 
our constructions however, we will use whichever direction is more 
convenient, without making a fuss over it.
\end{remark}

\begin{definition}
  \label{def:LaxBimodT} The 2-category $\lax{\Alg{\Bimod{T}}}$ of lax 
algebras 
for the pseudo-monad $(\Bimod{T},\im{\eta},\im{\mu})$ on 
$\Bimod{\twocat{K}}$
consists of:

\begin{description}
\item[objects] normal lax algebras $\rel{M}{TX}{X}$.

\item[morphisms] Morphisms of lax algebras $\morph{(f,\theta_f)}{M}{N}$.

\item[2-cells] Given morphisms $\morph{(f,\theta_f),(g,\theta_g)}{M}{N}$, 
a 2-cell 
$\cell{\alpha}{(f,\theta_f)}{(g,\theta_g)}$ consists of a 2-cell 
$\cell{\alpha}{f}{g}$ in $\twocat{K}$ such that

\newpage

\begin{displaymath}
\xymatrixrowsep{3pc}
 \xymatrixcolsep{4pc}
\let\labelstyle=\scriptstyle
  \begin{diagram}
 {TX}\ruppertwocell<1>^{M}{<4>\theta_g}&{X} \ar@{}[drr]|{=}
& &{TX}\ruppertwocell<1>^{M}{<4>\theta_f}&{X}\\
{TY}\utwocell<5>^{(Tf)^{*}}_{(Tg)^{*}}{{(T\alpha)^{*}}}\ar@/_1ex/[r]_-{N}&{Y}\ar@/_1ex/[u]_{g^{*}}
& &{TY}\ar@/^1ex/[u]^-{(Tg)^{*}}\ar@/_1ex/[r]_-{N}&{Y}
\utwocell<5>^{f^{*}}_{g^{*}}{{\alpha^{*}}}
  \end{diagram}
\end{displaymath}

\item[identities] $\morph{(\mathit{id}_X,\mathit{id}_M)}{M}{M}$


\item[composition]  Composition of 1-cells is given by composition of 1-cells in $\twocat{K}$ and the pasting

\begin{displaymath}
\xymatrixrowsep{1.5pc}
 \xymatrixcolsep{1.5pc}
\let\objectstyle=\scriptscriptstyle
\let\labelstyle=\scriptscriptstyle
  \begin{diagram}
 {TX}\ruppertwocell<1>^{{M}^{}}{<3>\theta_f}&{X} \\
{TY}\ar@/^1ex/[u]^-{(Tf)^{*}}
\ruppertwocell<1>_{{N}^{}}{<3>\theta_h}&{Y}\ar@/_1ex/[u]_{f^{*}}\\
{TZ}\ar@/^1ex/[u]^-{(Th)^{*}}\ar@/_1ex/[r]_-{P}&{Z}\ar@/_1ex/[u]_{h^{*}}
  \end{diagram}
\end{displaymath}

\noindent while compositions of 2-cells are inherited from $\twocat{K}$.

\end{description}
\end{definition}
 
The 2-category $\lax{\Alg{\Bimod{T}}}$ would be the basis on which 
pseudo-$T$-algebras become properties. We will give an alternative 
description of this 2-category in terms of 
$\BimodT{\mathbf{T}}{\twocat{K}}$, which, besides being of interest in 
itself, will be helpful to recapture our guiding example of 
multicategories, \cf. Theorem \ref{thm:bimod-spn} in Part \ref{part:applications}. 

\begin{definition}
\label{def:mnd-bimod}
 The 2-category $\Mnd{\BimodT{\mathbf{T}}{\twocat{K}}}$
of (normal) \textbf{monads} in  
  $\BimodT{\mathbf{T}}{\twocat{K}}$ consists of:

  \begin{description}
  \item[objects] monads in $\BimodT{\mathbf{T}}{\twocat{K}}$, that is 
    \begin{itemize}
    \item An object $X$ in $\twocat{K}$ and a bimodule $\rel{M}{TX}{X}$
    
    \item Unit $\cell{\iota}{\eta_X^{*}}{M}$ and multiplication 
$\cell{m}{\mu_X^{*}{\bullet}TM{\bullet}M}{M}$ satisfying the unit and 
associativity axioms.
    \end{itemize}
We furthermore require the following \textbf{normality} condition: the 
adjoint transpose of $\iota$ across $\im{(\eta_X)}\dashv \eta_X^{*}$ is 
$\morph{\mathit{id}}{\Hom{X}}{\Hom{X}}$. Henceforth we call monads 
satisfying this condition \textbf{normal}.

\item[morphisms] A morphism from $\rel{M}{TX}{X}$ to $\rel{N}{TY}{Y}$ is 
given by a pair of morphisms $\morph{f}{X}{Y}$ and $\morph{f_h}{M}{N}$ in 
$\twocat{K}$ such that


\begin{displaymath}
\xymatrixrowsep{1.5pc}
 \xymatrixcolsep{1.5pc}
\let\objectstyle=\scriptscriptstyle
\let\labelstyle=\scriptscriptstyle
  \begin{diagram}
 {TX}\dto_-{Tf}&\lto_-{d}{M}\dto|-{f_h}\rto^-{c}&{X}\dto^-{f} \\
{TY}&\lto^-{d'}{N}\rto_-{c'}&{Y}
  \end{diagram}
\end{displaymath}
\noindent such $f_h$ preserves the bimodule structure and is compatible 
with the units and multiplications of $M$ and $N$.

\item[2-cells] Given morphisms $\morph{(f,f_h),(g,g_h)}{M}{N}$, a 2-cell 
$\cell{\phi}{(f,f_h)}{(g,g_h)}$ consists of a 2-cell 
$\cell{\phi}{M}{(Tf,g)^{*}N}$ in 
$\Bimod{\twocat{K}}(TX,X)$ which is a morphism of $\bimod{(M,M)}$ (left 
$M$-right $M$ modules). This last requirement makes sense:
$M$, being a monoid in $\BimodT{\mathbf{T}}{\twocat{K}}(X,X)$ acts on 
itself by composition, while $N$ has a similar $\bimod{(N,N)}$ structure, 
which is transferred by change-of-base along the morphisms 
$\morph{(f,f_h),(g,g_h)}{M}{N}$.
  \end{description}

\end{definition}

\begin{remarks}
\hfill{ } 
\begin{itemize}
\item  Our definition of morphism of monads above is a special case of that of 
\cite{Street72}, where we have restricted the 1-cells to be representable 
bimodules. This is in accordance with the similar restriction we imposed 
on morphisms of lax algebras.

\item  Our definition of 2-cells mimics that  of 
transformation between morphisms of multicategories 
\cite[Def.~6.6]{Hermida99a}.

\end{itemize}

\end{remarks}

  Consistently with our treatment of morphisms of multicategories in \cite{Hermida99a}, we say that a 1-cell
in $\Mnd{\BimodT{\mathbf{T}}{\twocat{K}}}$ as above is \textbf{full and faithful} if the corresponding morphism of bimodules
$\morph{f_h}{M}{\im{Tf}{\bullet}N{\bullet}g^{*}}$ is an isomorphism, so that we have a change of base situation.

We now set about to show that the 2-categories $\lax{\Alg{\Bimod{T}}}$ and
$\Mnd{\BimodT{\mathbf{T}}{\twocat{K}}}$ are essentially the same. The only 
subtlety in this correspondence is the identification of 2-cells. Given a 
(normal) lax algebra $\rel{M}{TX}{X}$ with 
structural 2-cells
$\cell{\iota}{\Hom{X}}{\im{(\eta_X)}{\bullet}M}$ and 
$\cell{m}{TM{\bullet}M}{\im{(\mu_X)}{\bullet}M}$, we obtain a (normal) 
monad
$\rel{M}{TX}{X}$ with unit and multiplication obtained from $\iota$ and 
$m$ by transposing across the adjunctions $\im{(\eta_X)}\dashv\eta_X^{*}$ 
and $\im{(\mu_X)}\dashv\mu_X^{*}$. Let us denote the resulting data 
$(\check{M},\check{\iota},\check{\mu})$.

\begin{proposition}
  \label{prop:laxalg-mnd}
There is an isomorphism of 2-categories 
$$\isomorph{\check{(\_)}}{\lax{\Alg{\Bimod{T}}}}{\Mnd{\BimodT{\mathbf{T}}{\twocat{K}}}}$$
\end{proposition}

\begin{proof}
  We must first verify that $(\check{M},\check{\iota},\check{\mu})$ is 
indeed a monad. This follows by a routine calculation using the unit and 
associativity axioms for a lax algebra 
(Def.~\ref{def:lax-alg}). The normality condition is straightforward. It 
is furthermore clear that starting with a normal monad $(M,\iota,m)$ we 
obtain the data for a normal lax algebra by taking adjoint transposes of 
$\iota$ and $m$, and this correspondence is inverse to $\check{(\_)}$.

As for morphisms, to give 
\begin{displaymath}
\xymatrixrowsep{1.5pc}
 \xymatrixcolsep{1.5pc}
\let\objectstyle=\scriptscriptstyle
\let\labelstyle=\scriptscriptstyle
  \begin{diagram}
 {TX}\dto_-{Tf}&\lto_-{d}{M}\dto|-{f_h}\rto^-{c}&X\dto^-{f} \\
{TY}&\lto^-{d'}{N}\rto_-{c'}&{Y}
  \end{diagram}
\end{displaymath}

\noindent between monads is the same as to give a morphism of bimodules 
$\cell{\hat{f_h}}{M}{(Tf,f)^{*}(N)}$ (by change-of-base). But 
$(Tf,f)^{*}(N) = \im{(Tf)}{\bullet}N{\bullet}f^{*}$, and so we have the 
following correspondence

\[
\prooftree{
\cell{\hat{f_h}}{M}{\im{(Tf)}{\bullet}N{\bullet}f^{*}}}
\Justifies{
\cell{\theta_f}{(Tf)^{*}{\bullet}M}{N{\bullet}f^{*}}}
\using
{\im{Tf}\dashv{f^{*}}}
\thickness=0.08em
\endprooftree
\]

\noindent which sets up the bijective correspondence between morphisms of 
lax algebras $\morph{(f,\theta_f)}{M}{N}$ and morphisms of monads 
$\morph{(f,f_h)}{\check{M}}{\check{N}}$.

As for 2-cells, given one in $\lax{\Alg{\Bimod{T}}}$, $\cell{\alpha}{f}{g}$ for parallel 
morphisms $\morph{(f,\theta_f),(g,\theta_g)}{M}{N}$, we get one in 
$\Mnd{\BimodT{\mathbf{T}}{\twocat{K}}}$ by precomposing with the adjoint transpose of $\theta_f$: $$\underline{\alpha} \equiv 
\begin{diagram}
\let\labelstyle=\scriptscriptstyle
\xymatrixcolsep{2pc}
  {M}\ar@{=>}[r]^-{\check{\theta_f}}&{\im{Tf}{\bullet}N{\bullet}f^{*}}
\ar@{=>}[rr]^-{{\im{Tf}{\bullet}N{\bullet}\alpha^{*}}}&&
{\im{Tf}{\bullet}N{\bullet}g^{*}}
\end{diagram}
$$
\noindent and in the opposite direction, given $\cell{\phi}{M}{\im{Tf}{\bullet}N{\bullet}g^{*}}$ the composite 

$$
\xymatrixrowsep{2.5pc}
 \xymatrixcolsep{2.5pc}
\let\objectstyle=\scriptscriptstyle
\let\labelstyle=\scriptscriptstyle
  \begin{diagram}
 {X}\ar@/^3ex/[rrr]^(.3){{\Hom{X}}^{}}
\ar@/_1ex/[dr]_-{{\im{f}}^{}}\ar@{}[drr]|{\tiny \cong}
\ar@/^1ex/[r]_-{{\eta^{*}_X}^{}}&{TX}\rrtwocell^{M}_{{\im{Tf}{\bullet}N{\bullet}g^{*}}^{}}{{\phi}^{}}
& &{X}\\
&{Y}\ar@/_2ex/[r]_-{{\Hom{Y}}^{}}
\rcompositemap<1>_{{\eta^{*}_Y}^{}}^{N}{\omit}&{Y}\ar@/_1ex/[ur]_-{{g^{*}}^{}}&
  \end{diagram}
 $$

\noindent where the isomorphism is that of pseudo-naturality of $\eta^{*}$, yields
a 2-cell \linebreak
$\cell{\overline{\phi}}{\Hom{X}}{\im{f}{\bullet}g^{*} = \CommaObj{f}{g}}$
and therefore a 2-cell $\cell{\overline{\phi}}{f}{g}$ in $\twocat{K}$ which is furthermore well-defined as a 2-cell in $\lax{\Alg{\Bimod{T}}}$. These correspondences of 2-cells are readily verified to be mutually inverse.
\qed
\end{proof}

\begin{remarks}
\hfill{ }
  \begin{itemize}
  \item The best way to understand the correspondence at the level of 2-cells in the above proof is to look at the analysis of ordinary natural transformations in $\Cat$ in \cite[\S 6.1]{Hermida99a}.

  \item The correspondence at the object level, namely lax algebras versus monads, relies purely on the fact that the unit and multiplication of $\Bimod{T}$ have right adjoints. However, the rest of the setup relies heavily on the relationship between $\twocat{K}$ and $\Bimod{\twocat{K}}$, construing morphisms of the former as maps in the latter. Notice in particular that we define 
2-cells for monads by means of `tents' 
(commuting diagrams of 1-morphisms) using the fact that bimodules are spans. 
  \end{itemize}
 
\end{remarks}

\section{Pseudo-algebras and properties}
\label{sec:pseudo-alg-prop}

Having introduced our basic new gadget, namely the 2-category $\Mnd{\BimodT{T}{\twocat{K}}}$ (and its equivalent $\lax{\Alg{\Bimod{T}}}$), we now proceed to our main point: 
the 2-category of $\mathsf{T}$-algebras is {\em monadic\/} over  $\Mnd{\BimodT{T}{\twocat{K}}}$ and furthermore this monad has the {\em adjoint-pseudo-algebra\/} property. We then show that a pseudo-$\mathsf{T}$-algebra 
amounts to a {\em universal property\/} of a lax-$\Bimod{T}$-algebra, namely that the unit of the monadic adjunction have a left adjoint. We write $\Alg{\mathsf{T}}$ for the 2-category of $\mathsf{T}$-algebras, strict morphisms and 2-cells compatible with such (\ie. modifications).

\subsection{The adjunction between $\lax{\Alg{\Bimod{T}}}$ and $\Alg{\mathsf{T}}$ }
\label{sub:adjunction}
\noindent
\underline{\textbf{From $\Alg{\mathsf{T}}$ to 
$\lax{\Alg{\Bimod{T}}}$
}:} Given a $\mathsf{T}$-algebra 
$\morph{x}{TX}{X}$, the representable bimodule $\rel{\im{x}}{TX}{X}$ has a 
(normal) lax $\Bimod{T}$-algebra structure:
\begin{itemize}
\item The (identity) unit is $\Hom{X} = \im{(x\circ\eta_X)} = 
\im{(\eta_X)}\bullet\im{x}$, by the unit equation for $x$.
\item The multiplication is similarly obtained from the associativity for 
$x$:
  \begin{displaymath}
    \im{Tx}\bullet\im{x} = \im{(x{\circ}Tx)} = \im{(x{\circ}\mu_X)} = 
\im{(\mu_X)}\bullet\im{x}
  \end{displaymath}
\end{itemize}

A morphism of $\mathsf{T}$-algebras $\morph{f}{x}{y}$ (with 
$\morph{y}{TY}{Y}$) induces a morphism of lax algebras
\begin{displaymath}
\xymatrixrowsep{1.5pc}
 \xymatrixcolsep{1.5pc}
\let\objectstyle=\scriptscriptstyle
\let\labelstyle=\scriptscriptstyle
  \begin{diagram}
 {TX}\dto_-{x}\rto^-{Tf}&{TY}\dto^-{y}\ar@{}[drr]|\mapsto& 
&{TX}\dlowertwocell<-1>_{\im{x}}{^<-4>\theta_f^{}}&\ar@/_1ex/[l]_{Tf^{*}}{TY}
\ar@/^1ex/[d]^-{\im{y}}\\
{X}\rto_-{f}&{Y}& &{X}&\ar@/^1ex/[l]^{f^{*}}{Y}
  \end{diagram}
\end{displaymath}
\noindent by adjoint transposition and similarly a 2-cell 
$\cell{\alpha}{f}{g}$ 
between two such morphisms induces one 
$\cell{\alpha}{f}{g}$ in  $\lax{\Alg{\Bimod{T}}}$. We have thus defined a 
2-functor

\begin{displaymath}
  \morph{R}{\Alg{\mathsf{T}}}{\lax{\Alg{\Bimod{T}}}}
\end{displaymath}

\vspace{1em}
\noindent
\underline{\textbf{The free $\mathsf{T}$-algebra on a lax 
$\Bimod{T}$-algebra}:}
Given a (normal) lax algebra $\rel{M}{TX}{X}$, we consider it as a 
(normal) monad, according to Proposition \ref{prop:laxalg-mnd}.

\begin{lemma}
\label{lemma:Fmonad}
Given a (normal) monad $(\rel{M}{TX}{X},\iota_M,m_M)$ (in 
$\BimodT{T}{\twocat{K}}$), 
 the composite bimodule $\rel{\mu_X^{*}{\bullet}TM}{TX}{TX}$ has a monad 
structure in $\Bimod{\twocat{K}}$ 
\end{lemma}
\begin{proof}
\hfill{ }
  \begin{description}
  \item[unit] 

\begin{displaymath}
\xymatrixrowsep{5pc}
 \xymatrixcolsep{3pc}
\let\objectstyle=\scriptscriptstyle
\let\labelstyle=\scriptscriptstyle
  \begin{diagram}
{TX}\ar@/_1ex/[r]_-{^{\mu_X^{*}}}\ar@/^6ex/[rr]^-{\Hom{X}}&{T^2X}\rtwocell^{T\eta_X^{*}}_{TM}{{}^{\tiny T\iota_M}}&{TX}
\end{diagram}
\end{displaymath}

\newpage

\item[multiplication] 
\begin{displaymath}
\xymatrixrowsep{1.5pc}
 \xymatrixcolsep{2pc}
\let\objectstyle=\scriptscriptstyle
\let\labelstyle=\scriptscriptstyle
  \begin{diagram}
&&{TX}\ar@/^1ex/[dr]^-{^{\mu_X^{*}}}& &\\
& 
{T^2X}\ar@/^1ex/[ur]^-{TM}\ar@/_1ex/[dr]_-{{\mu_{TX}^{*}}^{}}\xtwocell[0,2]{}\omit{{\tiny    \mu_M}^{}}&&
{T^2X}\ar@/^1ex/[dr]^-{{TM}^{}}&\\
{TX}\ar@/^1ex/[ur]^-{{\mu_X^{*}}^{}}
\ar@/_1ex/[dr]_-{{\mu_X^{*}}^{}}&&{T^3X}\uruppertwocell<1>_{{T^2M}^{}}{<4>{}{\tiny    Tm_M}}&&{TX}\\
&{T^2X}\ar@/^1ex/[ur]_-{{{T\mu_X}^{*}}^{}}\ar@/_1ex/[urrr]_-{{TM}^{}}&&&
 \end{diagram}
\end{displaymath}

  \end{description}

\noindent The monad axioms follow from those of $(TM,T\iota_M,Tm_M)$ and 
pseudo-naturality of $\mu^{*}$.

\qed
\end{proof}

\newpage

Given a (normal) monad $(\rel{M}{TX}{X},\iota_M,m_M)$, let 

\begin{displaymath}
\xymatrixrowsep{1.5pc}
 \xymatrixcolsep{2pc}
\let\objectstyle=\scriptscriptstyle
\let\labelstyle=\scriptscriptstyle
 \begin{diagram}
{TX}\ar@/_1ex/[dr]_-{\im{J}}
\ar@/^1ex/[r]^-{\mu_X^{*}}&{T^2X}\xtwocell[1,0]{}\omit{\rho}\ar@/^1ex/[r]^-{TM}&{TX}\ar@/^1ex/[dl]^-{\im{J}}
\\
& {FM}&
 \end{diagram}
\end{displaymath}
 
\noindent be the (representable) Kleisli object of the monad of Lemma \ref{lemma:Fmonad}. 
By Assumption \ref{assumptionKleisliBimod}, $\Bimod{T}$ preserves such Kleisli 
object. Hence, by universality, we have a unique mediating morphism $\morph{\sigma}{T(FM)}{FM}$ in 

\newpage

\begin{displaymath}
\xymatrixrowsep{1pc}
 \xymatrixcolsep{1.5pc}
\let\objectstyle=\scriptscriptstyle
\let\labelstyle=\scriptscriptstyle
  \begin{diagram}
{T^2X}\ar@/_1ex/[dd]_-{^{\im{(\mu_X)}}}
\ruppertwocell<1>^{T\mu_X^{*}}{<6>\mathit{\scriptscriptstyle 
can}^{}}&{T^3X}\ddto|(.3){^{\im{(\mu_{TX})}}}
\ruppertwocell<1>^{{T^2M}^{}}{<6>{\tiny   \mu_M}^{}}&{T^2X}
\ar@/^1ex/[dd]^-{^{\im{(\mu_X)}}}\ar@{}[dddrrr]|{=}
& &{T^2X}\ar@/_1ex/[dr]_-{\im{TJ}} \ar@/^1ex/[r]^{T\mu_X^{*}}
&{T^3X}
\xtwocell[1,0]{}\omit{T\rho}\ar@/^1ex/[r]^{T^2M}&{T^2X}
\ar@/^1ex/[dl]^-{\im{TJ}}
\\
& 
& & & & {T(FM)}\ar@{-->}[dd]^{\im{\sigma}}& \\
{TX}\ar@/_1ex/[dr]_-{\im{J}}
\ar@/_1ex/[r]_-{\mu_X^{*}}&{T^2X}\xtwocell[1,0]{}\omit{\rho^{}}\ar@/_1ex/[r]_-{TM}&{TX}\ar@/^1ex/[dl]^-{\im{J}}\\
& {FM}& & & & {FM}&
\ignore{{TX}\ar@/^1ex/[ur]_-{\im{J}}
\ar@/_1ex/[r]_-{\mu_X^{*}}&{T^2X}\xtwocell[-1,0]{}\omit{^\rho^{}}\ar@/_1ex/[r]_-{TM}&{TX}\ar@/_1ex/[ul]^-{\im{J}}}
\end{diagram}
\end{displaymath}

\noindent where $\mathit{can}$ is the adjoint transpose of the 
associativity square for $\im{\mu}$ and $\mu_M$ is the invertible 2-cell 
corresponding to the pseudo-naturality of $\im{\mu}$. It follows by 
uniqueness of mediating morphisms between universal lax cocones (Kleisli 
objects) that $\morph{\sigma}{T(FM)}{FM}$ is a $\mathsf{T}$-algebra, using the 2-functorial preservation of representability 
of Assumption \ref{assumptionKleisliBimod}. Using the same assumption, we can extend the action of $F$ to morphisms and 2-cells, 
thereby defining a 2-functor
\newpage

\begin{displaymath}
  \morph{F}{\lax{\Alg{\Bimod{T}}}}{\Alg{\mathsf{T}}}
\end{displaymath}

\begin{proposition}
\label{prop:fund-adj}
  The above 2-functors set up an adjunction of 2-categories 
$\morph{F\dashv R}{\Alg{\mathsf{T}}}{\lax{\Alg{\Bimod{T}}}}$ whose unit is 
full and faithful.
\end{proposition}

\begin{proof}
 \mbox{ }\hfill{ }
  \begin{description}
  \item[unit] 
    Given $(\rel{M}{TX}{X},\iota_M,m_M)$, define 
    $\cell{\zeta_M}{M}{RF(M)}$ as follows:

\begin{displaymath}
\xymatrixrowsep{1.5pc}
 \xymatrixcolsep{1.5pc}
\let\objectstyle=\scriptscriptstyle
\let\labelstyle=\scriptscriptstyle
  \begin{diagram}
    &
{TX}\ar@/_1ex/[dl]_-{^{\im{(T\eta_X)}}}\ar@/^1ex/[dr]_-{^{\im{(\eta_{TX})}}}
\ar@/^1ex/[r]^{M}\xtwocell[1,2]{}\omit{{\eta_M}^{}}&{X}\ar@/^1ex/[dr]^-{^{\im{(\eta_{X})}}} 
& \\
    {T^2X}\ar@/_1ex/[d]_-{\im{TJ}} \ar@/_1ex/[dr]^-{\im{(\mu_X)}}&&
{T^2X}\ar@/^1ex/[dl]|-{\im{(\mu_X)}}\ruppertwocell<1>^{TM}{<3>\hat{\rho}^{}} 
&
    {TX}\ar@/^1ex/[d]^-{\im{J}}\\
    {T(FM)}\ar@/_2ex/[rrr]_-{\im{\sigma}}&{TX}\ar@/^1ex/[rr]^-{\im{J}}&
    &{FM} 
  \end{diagram}
\end{displaymath}

\newpage
\item[counit] Given a $\mathsf{T}$-algebra $\morph{x}{TX}{X}$, we have a 
lax cocone
\begin{displaymath}
\xymatrixrowsep{1.5pc}
 \xymatrixcolsep{2pc}
  \begin{diagram}
{TX}\ar@/_1ex/[dr]_-{\im{x}}
\ar@/^1ex/[r]^-{\mu_X^{*}}&{T^2X}\xtwocell[1,0]{}\omit{\lambda}\ar@/^1ex/[r]^-{TM}&{TX}\ar@/^1ex/[dl]^-{\im{x}}
\\
& {X}&
 \end{diagram}
\end{displaymath}
\noindent where $\lambda$ is the adjoint transpose of 
$\im{Tx}\bullet\im{x} = \im{(\mu_X)}\bullet\im{x}$ (associativity for 
$x$). 
Thus we have an induced morphism $\morph{\hat{x}}{FM}{X}$ from the Kleisli 
object $FM$. Once again, associativity 
for $x$
and universality  imply that this induced morphism is a 
$\mathsf{T}$-algebra morphism from $\morph{\sigma}{T(FM)}{FM}$ to 
$\morph{x}{TX}{X}$, which is the desired counit $\epsilon_x = \hat{x}$.

\end{description}

\noindent 

To show $\zeta_{M}$ fully faithful, first let us notice that the 
underlying morphism $J\eta_{X}$ in $\twocat{K}$ is fully faithful as shown 
by the following diagram:

\begin{displaymath}
\xymatrixrowsep{2pc}
 \xymatrixcolsep{2pc}
\let\objectstyle=\scriptscriptstyle
\let\labelstyle=\scriptscriptstyle
  \begin{diagram}
{X}\ar@/_1ex/[d]_-{\im{(\eta_X)}}\ar@/^4ex/[rr]^-{\Hom{X}}
\ar@/^1ex/[r]_(.3){\mu_X^{*}}&
{TX}\ar@/^1ex/[r]_-{M}\xtwocell[1,1]{}\omit{{   \eta^{*}_M}^{}}&{X}\\
{TX}\ar@/^1ex/[ur]|{\Hom{X}}
\ar@{}[drr]|{\cong}\ar@/_1ex/[dr]_-{\im{J}}
\ar@/^1ex/[r]_-{{\mu_X^{*}}^{}}&
{T^2X}
\ar@/_1ex/[u]|{\eta^{*}_{TX}}\ar@/^1ex/[r]_-{{TM}^{}}&{TX}\ar@/_1ex/[u]_{\eta^{*}_X}
\\
& {FM}\ar@/_1ex/[ur]_-{J^{*}}&{  }
 \end{diagram}
\end{displaymath}

\noindent Finally,  full and faithfulness of $\zeta_M$ follows similarly, 
using that of $\eta_X$ (\cf. Proposition \ref{pr:unit}.(\ref{it:ff})).

\qed
\end{proof}

\subsection{Adjoint psuedo-algebras and representable bimodules}
\label{sec:adjoint-representable}

Let 
$\subadj{{\mathsf T}} = (\morph{\subadj{T}}{\lax{\Alg{\Bimod{T}}}}{\lax{\Alg{\Bimod{T}}}},\zeta,\nu)$ 
be the 2-monad induced by the adjuntion of 
Proposition \ref{prop:fund-adj}. In this subsection we will show that its 
pseudo-algebras are characterised as left adjoints to units and give an 
intrinsic characterisation of them in terms of representability.

\begin{proposition}
  \label{prop:Kock-prop}
The 2-monad $(\subadj{T},\zeta,\nu)$ has the adjoint-pseudo-algebra 
property, \ie. $\nu_{M}\dashv \zeta_{\subadj{T}M}$ for all normal lax 
algebras $M$.
\end{proposition}
\begin{proof}
  Since the counit of the required adjunction is the identity, we must 
simply define the unit.
\newpage

Notice that $\rel{\subadj{T}M = \im{(\sigma)}}{T(FM)}{FM}$
and $\morph{\nu_M}{\subadj{T^2}M}{\subadj{T}M}$ is the morphism of 
algebras 
$\morph{\widehat{\sigma}}{F(\subadj{T}M)}{FM}$ uniquely determined in
\begin{displaymath}
\xymatrixrowsep{1.5pc}
 \xymatrixcolsep{2pc}
\let\objectstyle=\scriptscriptstyle
\let\labelstyle=\scriptscriptstyle
  \begin{diagram}
{T(FM)}\ar@/^3ex/[drr]^-{\im{\sigma}^{}}\ar@/^1ex/[dr]|{{\im{J'}}^{}}& & \\
{T^2(FM)}\xtwocell[0,1]{}\omit{{\rho'}^{}}
\ar@/^1ex/[u]^{\im{(T\sigma)}}&{F(\subadj{T}M)}\ar@/^1ex/[r]|-{\im{\widehat{\sigma}}^{}}&{FM} 
\\
{T(FM)}\ar@/^1ex/[u]^{\mu^{*}_{FM}}\ar@/_1ex/[ur]|{{\im{J'}}^{}}\ar@/_3ex/[urr]_{\im{\sigma}}
\end{diagram}
\end{displaymath}

\noindent by the lax cocone 
$\cell{\alpha}{\mu^{*}_{FM}\bullet{\im{(T\sigma)}}\bullet\im{(\sigma)}}{\im{(\sigma)}}$
given by the adjoint transpose of 
$\cell{\mathit{id}}{{\im{(T\sigma)}}\bullet\im{(\sigma)}}{\im{(\mu_{FM})}\bullet\im{(\sigma)}}$, 
while the underlying morphism of 
$\cell{\zeta_{\subadj{T}M}}{\subadj{T}M}{\subadj{T^2}M}$ is 
$\morph{J'\eta_{FM}}{FM}{F(\subadj{T}M)}$. To give a 2-cell 
$1\Rightarrow J'\eta_{FM}\widehat{\sigma}$ amounts to give, by the 
2-dimensional universal property of the Kleisli object,  a modification 
$J'\Rightarrow J'\eta_{FM}\sigma$ between the respective lax cocones. 
To define such modification, notice that the adjoint 
 transpose 
$\cell{\im{\underline{\rho'}}}{\im{(T\sigma)}\bullet\im{J'}}{\im{(\mu_{FM})}\bullet\im{J'}}$ 
in $\Bimod{\twocat{K}}$ corresponds to a 2-cell 
$\cell{\underline{\rho'}}{J'\mu_{FM}}{J'T\sigma}$. Define 
$\cell{\kappa}{J'}{J'\eta_{FM}\sigma}$ as the 
composite

\newpage

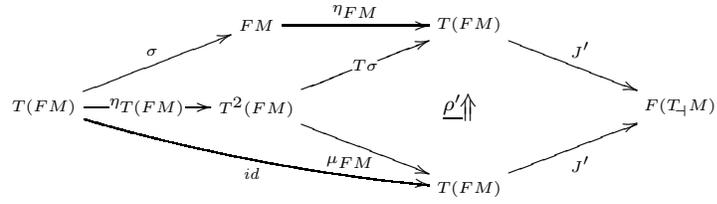
\begin{displaymath}
\xymatrixrowsep{1.5pc}
 \xymatrixcolsep{4pc}
\let\objectstyle=\scriptscriptstyle
\let\labelstyle=\scriptscriptstyle
  \begin{diagram}
&{FM}\rto^-{\eta_{FM}}&{T(FM)}\drto^-{J'}& \\
{T(FM)}\ar@/_1ex/[drr]_-{\mathit{id}}\urto^-{\sigma}\rto|-{\eta_{T(FM)}}&{T^2(FM)}\drto_{\mu_{FM}}\urto|-{T\sigma}\xtwocell[0,2]{}\omit{^{\underline{\rho'}}}&&F(\subadj{T}M)\\
&&{T(FM)}\urto_-{J'}& 
 \end{diagram}
\end{displaymath}

The verification that the 2-cell so defined is indeed a modification is 
quite delicate, so we outline the details. The equation
$$\rho' \circ (\mu^{*}_{FM}\bullet\im{(T\sigma)}\bullet\im{\kappa}) = 
\im{\kappa} \circ (\alpha\bullet\im{(\eta_{FM})}\bullet\im{J'})
$$
amounts to the equality (omitting objects to simplify the diagram)

\newpage

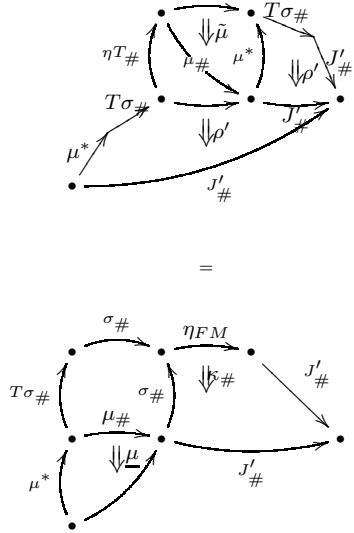
\begin{displaymath}
\xymatrixrowsep{2pc}
 \xymatrixcolsep{2pc}
\let\objectstyle=\scriptscriptstyle
\let\labelstyle=\scriptscriptstyle
  \begin{diagram}
&{\bullet}\ruppertwocell<1>{<1.5>{\tiny \tilde{\mu}}^{}}\ar@/_1ex/[dr]|-{\im{\mu}}&
{\bullet}\drcompositemap<2>_{\im{T\sigma}}^{\im{J'}}{\omit}&  \\
&{\bullet}\ar@/^1ex/[u]^-{\im{\eta{T}}}\ruppertwocell<-1>{<2.5>{\tiny \rho'}^{}}&{\bullet}
\ar@/_1ex/[u]^-{\mu^{*}}
\rlowertwocell<-1>_{\im{J'}}{<-2>{\tiny \rho'}^{}}&{\bullet} \\
{\bullet}\urcompositemap<0>_{\mu^{*}}^{\im{T\sigma}}{\omit}
\ar@/_2ex/[urrr]_-{\im{J'}}\ar@{}[ddrrr]|{=}
& & & \\
& & &\\
{\bullet}\ar@/^1ex/[r]^-{\im{\sigma}}& 
{\bullet}\ruppertwocell<1>^{\eta_{FM}}{<2>{\tiny    \im{\kappa}}^{}}&
{\bullet}\drto^-{\im{J'}} 
& { } \\
{\bullet}\ar@/^1ex/[u]^-{\im{T\sigma}}\ruppertwocell<1>^{\im{\mu}}{<1.5>
{\tiny \underline{\mu}}^{}}
&{\bullet}\ar@/_1ex/[u]^-{\im{\sigma}}\ar@/_1ex/[rr]_-{\im{J'}}& &{\bullet}\\
{\bullet}\ar@/^1ex/[u]^-{\mu^{*}}\ar@/_1ex/[ur] & & &
 \end{diagram}
\end{displaymath}

\noindent the upper pasting can be simplified using the fact that $\rho'$ 
is a lax cocone for the monad $\mu^*_X\bullet\im{(T\sigma)}$:
\begin{displaymath}
  \rho' \cdot \rho'\im{T\sigma}\mu^{*}_X = \rho'\cdot \im{J'}\overline{m}
\end{displaymath}
\noindent where $\overline{m} = Tm_{\sigma}\mu^{*}_X \cdot \im{T\sigma}\mu_{\sigma}\mu^{*}_X$ is the multiplication for the monad $\im{T\sigma}\mu^{*}_X$, \cf.~Lemma \ref{lemma:Fmonad}.
 So we will 
be done if we can show the following equality\footnote{See Diagram 1, attached.}:

 $$ \im{T\sigma}(T\underline{\mu}\mu^{*}{\circ}\mu_{\sigma}\mu^{*}{\circ}
\tilde{\mu}(\mu^*{\bullet}\im{T\sigma}{\bullet}\im{\eta{T}})) =
\im{T\sigma}(\tilde{\mu}\im{\eta{T}}{\circ}\im{\eta{T}}\underline{\mu})
$$

\noindent as pasting both sides with $\rho'$ gives the desired equality. 
\ignore{
\begin{displaymath}
\xymatrixrowsep{2pc}
 \xymatrixcolsep{2pc}
\let\objectstyle=\scriptscriptstyle
\let\labelstyle=\scriptscriptstyle
  \begin{diagram}
{\bullet}\ruppertwocell<1>{<1>\tilde{\mu}}\ar@/_1ex/[dr]|-{\im{\mu}}&
{\bullet}\ar@/^1ex/[r]^-{\im{(T\sigma)}}&{\bullet}  \\
{\bullet}\ar@/^1ex/[u]^-{\im{\eta{T}}}\ar@{}[dr]|{\cong}
\ar@/_1ex/[r]&{\bullet}\ar@/^1ex/[u]|-{{\mu^{*}}^{}}&  \\
{\bullet}\ar@/^1ex/[u]^{\im{(T\sigma)}}\ar@/_1ex/[r]^-{\mu{T}^{*}}&{\bullet}\ar@/_3ex/[uu]_{\im{{(T^2\sigma)}}}\ar@/^1ex/[r]^(.7){\im{T\mu}}&{\bullet}\ar@/_3ex/[uu]_{\im{(T\sigma)}}\\
{\bullet}\ar@{}[ddrr]|-{=}
\ar@/^1ex/[u]^-{\mu^{*}}\ar@/_1ex/[r]^-{\mu^{*}}&{\bullet}\ar@/^1ex/[u]^-{T\mu^{*}}
\urlowertwocell<-3>{<-1.5>{\tiny  T\underline{\mu}}^{}} & \\
& & \\
{\bullet}\ar@/^1ex/[r]_-{\im{(\sigma)}}& {\bullet}
\ar@/^1ex/[r]_-{\im{\eta}}&{\bullet} \\
{\bullet}\ar@/^1ex/[u]^-{\im{(T\sigma)}}\ar@/^1ex/[r]^{\im{\mu}}
&{\bullet}\ar@/_1ex/[d]|(.7){\im{\eta{T}}}&{\bullet}\ar@/_1ex/[u]_-{\im{(T\sigma)}} 
\\
{\bullet}\ar@/^1ex/[u]^-{\mu^{*}}\urlowertwocell<-2>{<-1.5>{\tiny   \underline{\mu}}^{}}
&{\bullet}\ar@/_1ex/[r]_-{\im{\mu}}\uruppertwocell<1>{<1.5>{\tiny  \tilde{\mu}}^{}}&
{\bullet}\ar@/_1ex/[u]_-{\mu^{*}}
 \end{diagram}
\end{displaymath}

\newpage
}
The trick to prove this latter equality is to realise the morphisms of spans which induce 
these 2-cells. The left one is induced by 


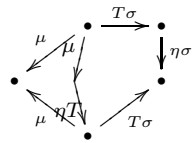
\begin{displaymath}
\xymatrixrowsep{1pc}
 \xymatrixcolsep{1.5pc}
\let\objectstyle=\scriptscriptstyle
\let\labelstyle=\scriptscriptstyle
  \begin{diagram}
& 
\dlto_-{\mu}{\bullet}\ddcompositemap<-1>_{\mu}^{\eta{T}}{\omit}\rto^-{T\sigma}&{\bullet}
\dto^-{\eta\sigma} \\
{\bullet}& & {\bullet} \\
&\ulto^-{\mu}{\bullet}\urto_-{T\sigma}
 \end{diagram}
\end{displaymath}
while the right one is induced by 


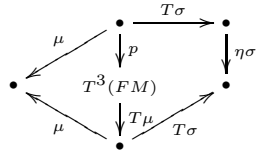
\begin{displaymath}
\xymatrixrowsep{1pc}
 \xymatrixcolsep{1.5pc}
\let\objectstyle=\scriptscriptstyle
\let\labelstyle=\scriptscriptstyle
  \begin{diagram}
& 
\dlto_-{\mu}{\bullet}\dto^-{p}\rto^-{T\sigma}&{\bullet}\dto^-{\eta\sigma}
\\
{\bullet}&{T^3(FM)}\dto^-{T\mu} & {\bullet} \\
&\ulto^-{\mu}{\bullet}\urto_-{T\sigma}
 \end{diagram}
\end{displaymath}
\noindent where in turn

\begin{displaymath}
\xymatrixrowsep{1.5pc}
 \xymatrixcolsep{1.5pc}
\let\objectstyle=\scriptscriptstyle
\let\labelstyle=\scriptscriptstyle
  
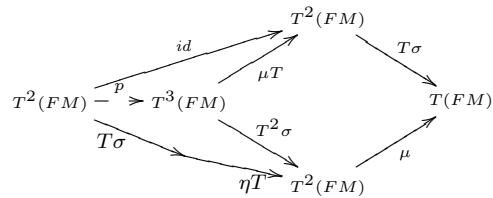
\begin{diagram}
&&{T^2(FM)}\drto^-{T\sigma}&\\
{T^2(FM)}\urrto^{\mathit{id}}\drrcompositemap<-1>^{{\eta{T}}^{}}_{{T\sigma}^{}}{\omit}
\ar@{-->}[r]^-{p}&{T^3(FM)}\urto_-{{\mu{T}}^{}}\drto^-{{T^2\sigma}^{}}& 
&{T(FM)} \\
& & {T^2(FM)}\urto_-{\mu} & 
 \end{diagram}
\end{displaymath}
\noindent using cartesiannes of $\mu$. Finally, naturality shows 
$\eta_{T(FM)}\circ\mu_{FM} = T\mu_{FM}\circ p$ and thus the morphisms of 
bimodules induced from those of  spans are equal, as desired. 

To conclude the proof, we must verify the adjunction equations for \linebreak
the 2-cell $\cell{\underline{\kappa}}{1}{J'\eta_{FM}\widehat{\sigma}}$ induced by the modification $\kappa$, namely

\begin{itemize}
\item $\widehat{\sigma}\kappa = \mathit{id}$, which is immediate by the 
definition of the lax cocone with 2-cell $\alpha$
which induces $\widehat{\sigma}$.

\item $\kappa J'\eta_{FM} = \mathit{id}$, which is established using $\rho' 
\circ (\iota\bullet\im{J'}) = \mathit{id} $ (lax cocone condition for the 
unit of the monad 
$\cell{\iota}{\Hom{T(FM)}}{\mu^*_{FM}\bullet\im{(T\sigma)}}$).
\end{itemize}

Finally, universality of (representable) Kleisli objects guarantees that $\underline{\kappa}$ is a 
well-defined 2-cell in $\lax{\Alg{\Bimod{T}}}$.

\qed
\end{proof}

Now we establish an important intrinsic characterisation of the 
adjoint-pseudo-algebras for $(\subadj{T},\zeta,\nu)$.
In fact, it was such a characterisation which motivated the present theory.

\begin{theorem}
\label{thm:representable-characterisation}
  A normal lax algebra $\rel{M}{TX}{X}$ with structure 2-cell 
\linebreak
$\cell{m}{TM{\bullet}M}{\im{(\mu_X)}{\bullet}M}$
is an adjoint-pseudo-algebra for the 2-monad $(\subadj{T},\zeta,\nu)$ (\ie.  
$\morph{\zeta_M}{M}{\subadj{T} M}$ admits a left adjoint) iff the 
following two conditions hold:
\begin{enumerate}
\item \label{it:representability} $\rel{M}{TX}{X}$ is a representable 
bimodule.

\item \label{it:structural-iso} 
$\cell{m}{TM{\bullet}M}{\im{(\mu_X)}{\bullet}M}$ is an isomorphism.
\end{enumerate}

\end{theorem}

\begin{proof}
\mbox{ }\hfill{ }
  \begin{description}
  \item[$\Longrightarrow$] Let $\rel{M}{TX}{X}$ have an 
adjoint-pseudo-algebra structure $\morph{(s,\bar{s})}{\subadj{T}M}{M}$, 
with 
$\cell{\tau}{\mathit{id}}{\zeta_{M}{\circ}s}$ and 
$\cell{\varepsilon}{s{\circ}\zeta_{M}}{\mathit{id}}$ the unit and 
(invertible) counit of the adjunction $s\dashv\zeta_M$.

Recall that the underlying morphism of $\zeta_M$ is $
\begin{diagram}
  {X}\rto|-{\eta_X}&{TX}\rto|-{J}&FM
\end{diagram}
$. Define $\morph{x = s{\circ}J}{TX}{X}$. We now intend to show that $x$ 
represents $M$. Define $\cell{\theta}{\im{x}}{M}$ as the pasting 


\begin{displaymath}
\xymatrixrowsep{2pc}
 \xymatrixcolsep{2pc}
\let\objectstyle=\scriptscriptstyle
\let\labelstyle=\scriptscriptstyle
  \begin{diagram}
{TX}\ar@{-->}[dr]\ar@/^5ex/[rrr]
\ar@/^1ex/[r]_-{{\im{(T\eta_X)}}^{}}&{T^2X}\ar@/_1ex/[d]^-{\mu_X^{}}
\rlowertwocell<-1>_{{\im{TJ}}^{}}{^<-2>{\tiny \im{T\varepsilon^{-1}}}^{}}
&{T(FM)}\xtwocell[1,1]{}\omit{^\bar{s}^{}}\ar@/_1ex/[d]_-{\im{\sigma}}
\ar@/^1ex/[r]_-{{\im{Ts}}^{}}&{TX}\ar@/^1ex/[d]^-{M}\\
&{TX}\ar@/_5ex/[rr]_-{{\im{x}}^{}}
\ar@/_1ex/[r]^-{{\im{J}}^{}}&{FM}\ar@/_1ex/[r]^-{{\im{s}}^{}}&{X}
 \end{diagram}
\end{displaymath}
 
We claim $\theta$ is an isomorphism: its inverse $\theta^{-1}$ is given by the pasting
\newpage

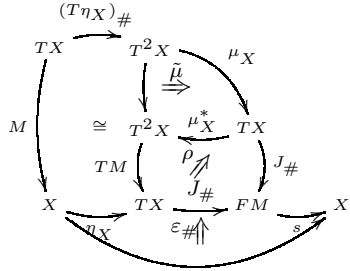
\begin{displaymath}
\xymatrixrowsep{1.5pc}
 \xymatrixcolsep{1.5pc}
\let\objectstyle=\scriptscriptstyle
\let\labelstyle=\scriptscriptstyle
  \begin{diagram}
{TX}\ar@{}[ddr]|-{\cong}\ar@/_1ex/[dd]_-{M^{}}
\ar@/^1ex/[r]^-{\im{(T\eta_X)}^{}}&{T^2X}\dlowertwocell<-1>{^<-2>{\tilde{\mu}}}
\ar@/^2ex/[dr]^-{\mu_X^{}}& & \\
&{T^2X}\xtwocell[1,1]{}\omit{^{\rho}^{}}\ar@/_1ex/[d]_{{TM}^{}}
&\ar@/^1ex/[l]_-{{\mu^{*}_X}^{}} {TX}\ar@/^1ex/[d]^-{\im{J}}& \\
{X}\ar@/_5ex/[rrr]
\ar@/_1ex/[r]_{\eta_X^{}}&{TX}\ruppertwocell<-1>^{\im{J}}{^<2>{\tiny 
\im{\varepsilon}}^{}}&{FM}\ar@/_1ex/[r]_{s}&{X}
\end{diagram}
\end{displaymath}

\noindent $\theta{\circ}\theta^{-1} = \mathit{id}$ because the counit is a 2-cell in 
$\lax{\Alg{\Bimod{T}}}$ and $\theta^{-1}{\circ}\theta = \mathit{id}$ by 
the adjunction equations for $\im{(J\eta_X)XS}\dashv s$ (since $\im{(J\eta_X)}\varepsilon^{-1} = \tau\im{(J\eta_X)XS}$).
We have thus shown that $M$ is representable. To see that its 
structure 2-cell $m$ is an isomorphism, we use the fact that 
$\morph{(s,\bar{s})}{\im{\sigma}}{M}$ commutes which the structural 
2-cells of these lax algebras. We have therefore

\newpage

\begin{displaymath}
\xymatrixrowsep{1.5pc}
 \xymatrixcolsep{1.5pc}
\let\objectstyle=\scriptscriptstyle
\let\labelstyle=\scriptscriptstyle
  \begin{diagram}
&{\bullet}\ar@/^1ex/[dr]^-{\im{\sigma}^{}}&  & &    
&{\bullet}\ar@/^1ex/[dd]|-{\im{Ts}^{}}
\xtwocell[3,1]{}\omit{{\tiny  \bar{s}}^{}}
\ar@/^1ex/[dr]^-{\im{\sigma}^{}}& \\
{\bullet}\ar@/_1ex/[dd]_{{\im{T^2s}}^{}}\ar@/^1ex/[dr]|-{\im{(\mu_X)}^{}}
\ar@/^1ex/[ur]^-{\im{T\sigma}^{}}& &{\bullet}
\ar@/^1ex/[dd]^{\im{s}^{}}
\ar@{}[drr]|{=}& & 
{\bullet}
\xtwocell[1,1]{}\omit{{\scriptscriptstyle  T\bar{s}}^{}}\ar@/_1ex/[dd]|-{{\im{T^2s}}^{}}
\ar@/^1ex/[ur]^-{\im{T\sigma}^{}}& &{\bullet}\ar@/^1ex/[dd]^{\im{s}^{}} \\
&{\bullet}\ar@/^1ex/[ur]_-{{\im{\sigma}}^{}}
\ar@/_1ex/[dd]|-{\im{Ts}^{}}\xtwocell[1,1]{}\omit{{\tiny  \bar{s}}^{}}&  & 
& &{\bullet}\ar@/^1ex/[dr]_-{M^{}} &  \\
{\bullet}\ar@{}[ur]|{\cong}\ar@/_1ex/[dr]_{\im{(\mu_X)}^{}}
& & {\bullet}& 
&{\bullet}\ar@/_1ex/[dr]_{\im{(\mu_X)}^{}}\ar@/^1ex/[ur]_-{{TM}^{}}
\xtwocell[0,2]{}\omit{{\tiny  m}^{}}&&{\bullet}\\
&{\bullet}\ar@/_1ex/[ur]_-{M^{}}& & & &{\bullet}\ar@/_1ex/[ur]_-{M^{}} &
 \end{diagram}
\end{displaymath}

\noindent Since $\theta$ is an isomorphism, so is $\theta\im{(Ts)}$,  as well as $\im{(Ts)}\im{T\tau}$ by adjointness, because $\varepsilon$ is an isomorphism and therefore $\varepsilon^{-1}s = s\tau$. From this adjointness 
we also conclude that the pasting of $\theta\im{(Ts)}$ and 
$\im{(Ts)}\im{T\tau}$
is $\bar{s}$. Hence,

\begin{equation}
  \cell{\bar{s}}{\im{\sigma}{\bullet}{\im{s}}}{\im{Ts}{\bullet}M}\mbox{ is an isomorphism} \label{eqn:bar-s}
\end{equation}

\noindent Thus, the equality of the pasting diagrams above imply that $m \im{T^2s}$ is an isomorphism and so is $m \im{T^2s}\im{(T^2J  T^2\eta_X)}$. Therefore,
\begin{displaymath}
  (m \im{T^2s}\im{(T^2J  T^2\eta_X)})\cdot \im{(\mu_X)}\im{T^2\varepsilon}
= m \cdot (M  TM  \im{T^2\varepsilon})
\end{displaymath}
\noindent is an isomorphism. Finally, $M  TM  \im{T^2\varepsilon}$ being an isomorphism allows us to conclude that $m$ is one as well.


\item[$\Longleftarrow$] Given $\morph{x}{TX}{X}$, the isomorphisms $\cell{\theta}{\im{x}}{M}$ and
$\cell{m}{TM{\bullet}M}{\im{(\mu_X)}{\bullet}M}$ endow $x$ with a pseudo-$\mathsf{T}$-algebra structure. We should therefore show that any such pseudo-algebra does endow $\rel{\im{x}}{TX}{X}$ with a pseudo-$\subadj{\mathsf{T}}$-algebra structure (which a fortriori would be an adjoint one, by Proposition \ref{prop:Kock-prop}). Let 
$\cell{\iota}{\mathit{id}}{x{\circ}\eta_X}$ and $\cell{\alpha}{x{\circ}Tx}{x{\circ}\mu_X}$
be the structural isomorphisms. The adjoint transpose of $\im{\alpha}$ produces a lax cocone

\newpage

\begin{displaymath}
\xymatrixrowsep{1.5pc}
 \xymatrixcolsep{2pc}
\let\objectstyle=\scriptscriptstyle
\let\labelstyle=\scriptscriptstyle
 \begin{diagram}
{TX}\ar@/_1ex/[dr]_-{\im{x}}
\ar@/^1ex/[r]^-{\mu_X^{*}}&{T^2X}\xtwocell[1,0]{}\omit{\check{\alpha}}\ar@/^1ex/[r]^-{TM}&{TX}\ar@/^1ex/[dl]^-{\im{x}}
\\
& {X}&
 \end{diagram}
\end{displaymath}
\noindent and we get an induced morphism $\morph{\hat{x}}{F(\im{x})}{X}$, which is equipped with structural isomorphisms

\begin{description}
\item[unit] $\cell{\im{\iota^{-1}}}{\Hom{X}}{\im{(x{\circ}\eta_X)} = \im{(J\eta_X)}\bullet\im{\hat{x}}}$

\item[associativity] $\cell{\hat{\alpha}}{\im{\subadj{T}\hat{x}}\bullet\im{\hat{x}}}{\nu_{\im{x}}}\bullet\im{\hat{x}}$ uniquely determined by the given $\cell{\alpha}{x{\circ}Tx}{x{\circ}\mu_X}$ and the 2-dimensional universal property of the Kleisli objects involved. Furthermore, universality of Kleisli objects ensure that the induced canonical isomorphisms satisfy the pseudo-$\mathsf{T}$-algebra axioms.

\end{description}

\end{description}

\qed
\end{proof}

A clear example of this theorem is given in the case of pseudo-functors into $\Cat$ in Part \ref{part:applications}, \cf. Remark 
\ref{rmk:representable-pseudo-functors}.

\subsection{Monadicity}
\label{sec:monadicity}

Let $\subadj{\twocat{K}} = \lax{\Alg{\Bimod{\mathsf{T}}}}$.
There is an evident `underlying object' 2-functor $\morph{U}{\subadj{\twocat{K}}}{\twocat{K}}$ with action 
$(\rel{M}{TX}{X}) \mapsto X$. Furthermore, there is a 2-natural transformation

\begin{displaymath}
\xymatrixrowsep{1.5pc}
 \xymatrixcolsep{2pc}
 \begin{diagram}
{\subadj{\twocat{K}}}\xtwocell[1,1]{}\omit{{J}^{}}
\dto_-{{\subadj{T}}^{}}
\rto^-{U}&{\twocat{K}}\dto^-{T}\\
{\subadj{\twocat{K}}}
\rto_-{U}&{\twocat{K}}
 \end{diagram}
\end{displaymath}
\noindent given by $J_M = 
\begin{diagram}
  TX\rto|-{J}&FM
\end{diagram}
$, the Kleisli morphism into the Kleisli object $FM = U(\subadj{T}(M))$. The 2-naturality of $J$ follows from the universality of Kleisli objects (with the preservation of representability of our Assumption \ref{assumptionKleisliBimod}).

\begin{proposition}
  The pair $(U,J)$ is a morphism of 2-monads from $(\subadj{T},\zeta,\nu)$ to $(T,\eta,\mu)$.
\end{proposition}
\begin{proof}
 We must verify the compatibility with units and multiplications:
 \begin{itemize}
 \item 
\begin{displaymath}
\let\objectstyle=\scriptscriptstyle
\let\labelstyle=\scriptscriptstyle
\xymatrixrowsep{2pc}
 \xymatrixcolsep{2pc}
 \begin{diagram}
{\subadj{\twocat{K}}}\xtwocell[1,1]{}\omit{<.8>{J}^{}}
\dto_-{{\subadj{T}}^{}}
\rto^-{U}&{\twocat{K}}\dtwocell_{T}^{\mathit{id}^{}}{{\eta}^{}}
\ar@{}[drr]|{=}&&{\subadj{\twocat{K}}}\dtwocell_{{\subadj{T}}^{}}^{\mathit{id}^{}}{{\zeta}^{}}\rto^-{U}&{\twocat{K}}\dto^{\mathit{id}}\\
{\subadj{\twocat{K}}}
\rto_{U^{}}&{\twocat{K}}&&{\subadj{\twocat{K}}}\rto_{U^{}}&{\twocat{K}}
 \end{diagram}
\end{displaymath}
\noindent follows immediately from the definition of $\zeta$ that $U(\zeta_M) = J{\circ}\eta_X$ (for $\rel{M}{TX}{X}$).


\item 
\begin{displaymath}
\let\objectstyle=\scriptscriptstyle
\let\labelstyle=\scriptscriptstyle
\xymatrixrowsep{1.5pc}
 \xymatrixcolsep{1.5pc}
 \begin{diagram}
{\subadj{\twocat{K}}}\xtwocell[1,1]{}\omit{^{J}^{}}
\dto_-{U}\rruppertwocell<8>^{{\subadj{T}}^{}}{^<-2>{{\tiny \nu}^{}}}
\rto_-{{\subadj{T}}^{}}&{\subadj{\twocat{K}}}\xtwocell[1,1]{}\omit{^{J}^{}}
\dto^-{U}\rto_-{{\subadj{T}}^{}}
&{\subadj{\twocat{K}}}\dto|-{U}\ar@{}[drr]|{=}&&{\subadj{\twocat{K}}}\xtwocell[0,2]{}\omit{^{{J}^{}}}
\dto_-{U}\ar@/^3ex/[rr]^-{{\subadj{T}}^{}}&&
{\subadj{\twocat{K}}}\dto^-{U}\\
{\twocat{K}}
\rto_{T^{}}&
{\twocat{K}}
\rto_{T^{}}&{\twocat{K}}&&
{\twocat{K}}
\rruppertwocell<8>^{T^{}}{^<-2>{\mu}^{}}\rto_-{T^{}}&{\twocat{K}}\rto_-{T^{}}&{\twocat{K}}
 \end{diagram}
\end{displaymath}

\noindent we have by definition of $\nu$ that the right pasting instantiated at $\rel{M}{TX}{X}$ is
the upper morphism in the following commuting diagram,

\begin{displaymath}
\let\objectstyle=\scriptscriptstyle
\let\labelstyle=\scriptscriptstyle
  \xymatrixrowsep{1.5pc}
 \xymatrixcolsep{2pc}
\begin{diagram}
{T^2X}\rto^-{TJ}
\rrcompositemap<-2>_{\mu^{}}^{J^{}}{\omit}
&{T(FM)}\rto^-{\sigma^{}}
\rcompositemap_{{J'}^{}}^{{\hat{\sigma}}^{}}{\omit}&{FM}
\end{diagram}
\end{displaymath}

\noindent while the bottom map is the corresponding instance of the left pasting.

 \end{itemize}
\qed

\end{proof}


\begin{theorem}[Monadicity]
  \label{thm:monadicity}

The morphism of 2-monads $\morph{(U,J)}{\subadj{\mathsf{T}}}{\mathsf{T}}$ induces 2-equivalences
\vskip5mm
\begin{center}
  \framebox[1.1\width]{$\alg{(U,J)}:\alg{\subadj{\mathsf{T}}}\simeq\alg{\mathsf{T}}
$}
\end{center}

\vskip5mm
\begin{center}
  \framebox[1.1\width]{$\pseudo{\alg{(U,J)}}:\pseudo{\alg{\subadj{\mathsf{T}}}}\simeq\pseudo{\alg{\mathsf{T}}}
$}
\end{center}

\end{theorem}
\begin{proof}
  The details are essentially contained in the proof of Theorem \ref{thm:representable-characterisation}. 
Given a (pseudo-)$\mathsf{T}$-algebra $\morph{x}{TX}{X}$ (with structural isomorphisms $\iota$ and $\alpha$) we have 
$$
\begin{diagram}
  {TX}\rto^{J}&{FM}\rto^{\hat{x}}&X & = & {TX}\rto^{x}&{X}
\end{diagram}
$$


\noindent by definition of $\hat{x}$. In the other direction, 
$$
\xymatrixrowsep{1.5pc}
 \xymatrixcolsep{1.5pc}
\begin{diagram}
 {TX}\xtwocell[1,1]{}\omit{\theta^{}}
\ar@/^1ex/[r]^{\im{x}}\ar@/_1ex/[d]
&{X} \ar@/^1ex/[d]\\
{TX}\ar@/_1ex/[r]_-{M}&{X}
\end{diagram}$$
\noindent is an isomorphism of (pseudo-)$\subadj{T}$-algebras.
\qed
\end{proof}

\begin{remark}
  The underlying object 2-functor $\morph{U}{\subadj{\twocat{K}}}{\twocat{K}}$ is locally conservative, a fibration at the 1-cell level (by change-of-base for bimodules) and has a left 2-adjoint with action $X \mapsto (\rel{\eta^{*}_X}{TX}{X})$.
\end{remark}


\section{Classification of lax morphisms}
\label{sec:lax-morph}
We use the adjunction 
$\morph{F\dashv{R}}{\Alg{\mathsf{T}}}{\lax{\Alg{\Bimod{T}}}}$ to obtain 
an explicit classification of lax morphisms between (pseudo-) $T$-algebras.

Given (pseudo-)$T$-algebras $\morph{x}{TX}{X}$ and $\morph{y}{TY}{Y}$, 
consider a morphism $\morph{(f,\theta_f)}{\im{x}}{\im{y}}$ in 
$\lax{\Alg{\Bimod{T}}}$. We have the following correspondences

\begin{displaymath}
\xymatrixrowsep{1.5pc}
 \xymatrixcolsep{1.5pc}
\let\objectstyle=\scriptscriptstyle
\let\labelstyle=\scriptscriptstyle
  \begin{diagram}
 {TX}\ar@/^1ex/[r]^{\im{x}}\xtwocell[1,1]{}\omit{{\tiny \theta_f}}
 &{X}
 \ar@{}[drr]|{\leftrightarrow}&&
 {TX}\ar@/_1ex/[d]_-{\im{Tf}}
 \ar@/^1ex/[r]^{\im{x}}\xtwocell[1,1]{}\omit{{\tiny \check{\theta_f}}}
 &{X}\ar@/^1ex/[d]^-{\im{f}}
 \ar@{}[drr]|{\leftrightarrow}&&
 {TX}\dto_-{Tf}
 \rto^{x}\xtwocell[1,1]{}\omit{^{\tiny \check{\theta_f}}}&{X}\dto^-{f} \\
 {TY}\ar@/^1ex/[u]^{Tf^{*}}
 \ar@/_1ex/[r]_-{\im{y}}&{Y}\ar@/_1ex/[u]_-{f^{*}}& &
 {TY} \ar@/_1ex/[r]_-{\im{y}}&{Y}& &{TY}\rto_-{y}&{Y}
  \end{diagram}
\end{displaymath}


\noindent the first one by taking adjoint mates and the second 
one using that \linebreak
$\morph{\im{(\_)}}{\twocat{K}^{\mathit{co}}}{\Bimod{\twocat{K}}}$ is 
locally fully faithful. The 
resulting data amounts to a \textbf{lax morphism} between the 
(pseudo-)$T$-algebras. There is a similar (and quite clear) 
identification at the level of 2-cells, so that we have the following 
result:

\begin{theorem}[Classification of lax morphisms]
\label{thm:classification-lax-morphisms}
\mbox{ }\hfill{ }
\begin{enumerate}
        \item  
        Given $T$-algebras $\morph{x}{TX}{X}$ and $\morph{y}{TY}{Y}$, 
        the adjunction  
$$\morph{F\dashv{R}}{\Alg{\mathsf{T}}}{\lax{\Alg{\Bimod{T}}}}$$ 
        induces the following isomorphisms
        \[ \Alg{\mathsf{T}}_l(x,y) \iso 
\lax{\Alg{\Bimod{T}}}(\im{x},\im{y}) \iso
        \Alg{\mathsf{T}}(FRx,y)\] 
        \noindent and therefore the following  isomorphism of 2-categories
        \[  \Alg{\mathsf{T}}_l \iso \Alg{{\mathsf T}}_G\]
        \noindent where $G=FR$ is the 2-comonad induced on 
$\Alg{{\mathsf T}}$ 
        and the second 2-category above its associated Kleisli 
construction.
        
        \item  Given pseudo-$T$-algebras $\morph{x}{TX}{X}$ and 
$\morph{y}{TY}{Y}$, 
        the biadjunction  
$\morph{F\dashv{R}}{\pseudo{\Alg{\mathsf{T}}}}{\lax{\Alg{\Bimod{T}}}}$ 
        induces the following equivalences
        \[ \pseudo{\Alg{\mathsf{T}}}_l(x,y) \simeq 
        \lax{\Alg{\Bimod{T}}}(\im{x},\im{y}) \simeq 
\pseudo{\Alg{\mathsf{T}}}(FRx,y)\] 
        \noindent and therefore the following  biequivalence of 
2-categories
        \[  \pseudo{\Alg{\mathsf{T}}}_l \simeq\pseudo{\Alg{{\mathsf T}}}_G\]
        \noindent where $G=FR$ is the pseudo-comonad induced on 
$\pseudo{\Alg{{\mathsf T}}}$ 
        and the second 2-category above its associated Kleisli 
construction.
\end{enumerate}
\end{theorem}

\begin{remark}
        A classification of lax morphisms for $T$-algebras for 2-monads 
{\em with 
        a rank\/} appears in \cite{BlackwellKellyPower89}. Besides relying 
on a 
        quite different hypothesis from ours, the work in \ibid.  does not provide any 
explicit 
        description of the corresponding free objects (admitedly, that paper
has 
        an altogether different flavour and purpose from the present 
one).  
        As we will see in Part 
        \ref{part:applications}, our construction  allows us 
        to recover (or discover) {\em monoid classifiers\/}.
\end{remark}

\section{Classification of strong  morphisms and coherence}
\label{sec:strict-morph}

Given a pseudo-$\subadj{T}$-algebra, we want to freely associate a strict 
algebra to it. In other words, we want to construct a reflection for the 
inclusion $\Alg{\subadj{T}}\hookrightarrow\pseudo{\Alg{\subadj{T}}}$. As 
we will see, we only get a `bireflection' in the sense that we will get 
equivalences rather than isomorphisms at the level of the `Hom' 
categories. Since $\subadj{T}$ has the adjoint-pseudo-algebra property, we 
intend to apply the technique for strictification we introduced in 
\cite[\S 10.2]{Hermida99a}. We reproduce the analysis of \ibid. in this abstract 
setting.

Let $M$ be an adjoint-pseudo-algebra, $s\dashv \zeta_{M}$ with unit 
$\cell{\tau}{1}{\eta_{M}s}$. Suppose $N$ is a strict $\subadj{T}$-algebra 
with structure $\morph{z}{\subadj{T}N}{N}$, and let 
$\morph{(f,\theta_f)}{M}{N}$ a strong morphism. Recall from 
\cite{Street73} that the associativity structure 2-cell 
$\cell{\alpha}{s{\circ}\mu_{M}}{s{\circ}Ts}$ for $M$ is given by the pasting 
(we write $\Upsilon = \subadj{T}$ to simplify notation)

\begin{displaymath}
\xymatrixrowsep{1.5pc}
 \xymatrixcolsep{1.5pc}
\let\objectstyle=\scriptscriptstyle
\let\labelstyle=\scriptscriptstyle
  \begin{diagram}
&{\bullet}\drto|-{\eta}\drrto|-{\eta}\ar@{->}[rrr]  & & &{\bullet} \\
{\bullet}\urto^-{s}\rrlowertwocell<0>{<-1.5>{}^{\tau}}&&{\bullet}\dto|-{^{\eta{\Upsilon}}}&{\bullet}\dlto|-{\Upsilon\eta}\rto^{^{\cong}}_(.3){{\cong}^{}}&{\bullet}\uto_{s}\\
{\bullet}\uto^-{\mu}\rrlowertwocell<0>{<-1.5>{\bar{\tau}}^{}}&&{\bullet}\urrto_-{\Upsilon 
s}&&
 \end{diagram}
\end{displaymath}

\noindent where $\bar{\tau}$ is the unit of 
$\mu_{M}\dashv\eta_{\subadj{T}M}$and the isomorphisms are the counits of the adjunctions.
The fact that $(f,\theta_f)$ is a 
strong morphism implies the equality 

\begin{displaymath}
\xymatrixrowsep{2pc}
 \xymatrixcolsep{2pc}
\let\objectstyle=\scriptscriptstyle
\let\labelstyle=\scriptscriptstyle
  \begin{diagram}
{\bullet}\rto^-{\Upsilon{s}}\dto_-{\mu}
\xtwocell[1,1]{}\omit{^{\alpha}}&{\bullet}\dto^-{s}& 
&{\bullet}\dto_-{\Upsilon{z}}&\lto_{\Upsilon^2f}{\bullet}\dto^{\Upsilon{s}}\\
{\bullet}\dto_-{\Upsilon{f}}\rto|-{s}
\xtwocell[1,1]{}\omit{^{}^{\theta_f}}&{\bullet}\dto^-{f}& {=} 
&{\bullet}\xtwocell[-1,1]{}\omit{{\Upsilon\theta_f}^{}}
\dto_-{z}&\lto|-{\Upsilon{f}}{\bullet}\dto^-{s}\\
{\bullet}\rto_-{z}&{\bullet}& &{\bullet}\xtwocell[-1,1]{}\omit{{\theta_f}^{}}
&\lto^-{f}{\bullet}
 \end{diagram}
\end{displaymath}

The morphism $\morph{f}{M}{N}$ induces
$\morph{\hat{f} = z{\circ}Tf}{TM}{N}$.  A little 
fiddling with the above diagrams shows that 

\[\xymatrixrowsep{2pc}
 \xymatrixcolsep{2.5pc}
\begin{diagram}
 {\subadj{T}{M}}\rruppertwocell^{1}{\omit}\relax
\rrcompositemap<-2>_{s^{}}^{\eta_{\scriptstyle{M}}^{}}{<-1>\tau}&&\relax
{\subadj{T}{M}}\rto^-{\hat{f}}&{N} &= & 
{\subadj{T}{M}}\rrcompositemap<2>_{s^{}}^{f_{}}{\omit}
\rrcompositemap<-2>_{{\subadj{T}{f}}^{}}^{z^{}}{{}^{\tiny     \theta_f^{-1}}}&&{N}
 \end{diagram}
\]

\noindent Hence $\morph{\hat{f}}{\subadj{T}M}{N}$ inverts $\tau$. Notice 
that $\subadj{T}M$ has a strict $\subadj{T}$-algebra structure (the free 
such over $M$). So a good object candidate for the free strict 
$\subadj{T}$-algebra for $M$ is the {\em coinverter\/} of $\tau$. It would 
give us the desired reflection if we can endow it with a 
$\subadj{T}$-algebra structure. We reproduce, without proof,  the 
key technical lemma \cite[Lemma 10.4]{Hermida99a}:

\begin{lemma}
\label{key-lemma}
Consider an adjunction $\adj{\eta}{\varepsilon}{l}{r}{C}{D}$ in a 
2-category, with $r$ full and faithful (which is equivalent to 
$\varepsilon$ being an isomorphism). Consider the coinverter of the unit
\[\xymatrixrowsep{2pc}
 \xymatrixcolsep{2.5pc}
\begin{diagram}
 D\rruppertwocell^{1}{\omit}\relax
 \rrcompositemap<-2>^{r^{}}_{l^{}}{<-1>\eta}&&\relax
 D\rto^-{q}&D[\eta^{-1}]
 \end{diagram}
\]
and the unique morphism $\morph{l'}{D[\eta^{-1}]}{C}$ induced by $l$. 
\begin{enumerate}
        \item  \label{coinv-equiv} 
        The morphisms $\morph{l'}{D[\eta^{-1}]}{C}$ 
        and $\morph{qr}{C}{D[\eta^{-1}]}$ form an adjoint equivalence.
        
        \item  \label{Kleisli-coinv}
The coinverter $D[\eta^{-1}]$ is the Kleisli object $\underline{(rl)}$ of the (idempotent) monad $\morph{rl}{D}{D}$ induced on $D$ by the given 
adjunction, so that 
        there is a canonical isomorphism 
$$\xymatrixrowsep{1.5pc}
 \xymatrixcolsep{1.5pc}
\begin{diagram}
 D\drto_{J}\rto^-{q}&D[\eta^{-1}]\ar@{-->}[d]\\
 &{\underline{(rl)}}
 \end{diagram}
$$ 
        
\end{enumerate}

\end{lemma}

In order to apply this lemma, we must show that we have the relevant Kleisli 
objects. In principle we are working in the ambient 2-category 
$\lax{\Alg{\Bimod{T}}}$, but the object part of the monad whose Kleisli 
object we should provide is a (free) strict $T$-algebra. Hence the 
following result will suffice for our purposes:

\begin{proposition}
\label{prop:Kleisli-lax}
Given a 2-monad $\mathsf{T}$ on $\twocat{K}$, if $\twocat{K}$ has Kleisli 
objects (for monads) and $\morph{T}{\twocat{K}}{\twocat{K}}$ preserves 
them, 
then the forgetful 2-functor 
$\morph{U}{\Alg{\mathsf{T}}_{\mathit{l}}}{\twocat{K}}$ creates Kleisli 
objects in the 2-category of algebras and lax morphisms.
\end{proposition}
\begin{proof}
 Let
\begin{displaymath}
\xymatrixrowsep{1.5pc}
 \xymatrixcolsep{1.5pc}
\let\objectstyle=\scriptscriptstyle
\let\labelstyle=\scriptscriptstyle
  \begin{diagram}
 {TX}\dto_-{x}\ruppertwocell<0>^{Tm}{<3>{}^{\theta_m}}&{TX}\dto^-{x} \\
{X}\rto_-{m}&{X}
  \end{diagram}
\end{displaymath}
\noindent be an endomorphism in $\Alg{\mathsf{T}}_{\mathit{l}}$ with 
$((m,\theta_m),\eta,\mu)$ a monad. Let 
\begin{displaymath}
\xymatrixrowsep{1pc}
 \xymatrixcolsep{2pc}
\let\objectstyle=\scriptscriptstyle
\let\labelstyle=\scriptscriptstyle
  \begin{diagram}
 {X}\drto_-{J}\rruppertwocell<0>^{m}{<2>{\rho}}&&{X}\dlto^-{J} \\
&{\underline{X}}&
  \end{diagram}
\end{displaymath}
\noindent be the Kleisli object of $(m,\eta,\mu)$ in $\twocat{K}$. Since 
$T$ preserves it, there is a uniquely determined morphism 
$\morph{\underline{x}}{T\underline{X}}{\underline{X}}$ mediating between 
lax cocones in the following diagram:
\begin{displaymath}
\xymatrixrowsep{2.5pc}
 \xymatrixcolsep{2pc}
\let\objectstyle=\scriptscriptstyle
\let\labelstyle=\scriptscriptstyle
  \begin{diagram}
  &{T\underline{X}}\ar@{-->}[rr]^-{\underline{x}}&&{\underline{X}}\\
  {TX}\urto^{TJ}\rrto^{x}&&X\urto^{J}&\\
  {TX}\uto^-{Tm}\uurto_(.3){TJ}\xtwocell[-2,1]{}\omit{<-1>{T\rho}}
  \rrlowertwocell<0>_{x}{<-3>{}^{\theta_m}}&&
  {X}\uto^-{m}\uurto_{J}\xtwocell[-2,1]{}\omit{<-1>{\rho}}
  &
  \end{diagram}
\end{displaymath}
\noindent Uniqueness of mediating morphisms between lax cocones means 
that the algebra axioms for 
$\morph{\underline{x}}{T\underline{X}}{\underline{X}}$ follow from those 
of $\morph{x}{TX}{X}$. The above lax cocone in 
$\Alg{\mathsf{T}}_{\mathit{l}}$ is the Kleisli object of 
$((m,\theta_m),\eta,\mu)$. To verify its universality, consider another 
lax cocone

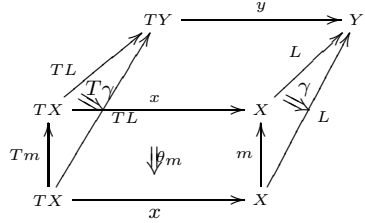
\begin{displaymath}
\xymatrixrowsep{2pc}
 \xymatrixcolsep{2pc}
\let\objectstyle=\scriptscriptstyle
\let\labelstyle=\scriptscriptstyle
  \begin{diagram}
  &{TY}\rrto^-{y}&&{Y}\\
  {TX}\urto^(.3){TL}\rrto^{x}&&X\urto^{L}&\\
  {TX}\uto^-{Tm}
\uurto_{TL}\xtwocell[-2,1]{}\omit{<-1>{T\gamma}}
  \rrlowertwocell<0>_{x}{<-3>{}^{\theta_m}}&&
  X\uto^-{m}
\uurto_{L}\xtwocell[-2,1]{}\omit{<-1>{\gamma}}
  &
  \end{diagram}
\end{displaymath}

\noindent we have a mediating morphism 
$\morph{\hat{L}}{\underline{X}}{Y}$ induced by universality of 
$\underline{X}$ and a 2-cell
\begin{displaymath}
\xymatrixrowsep{1.5pc}
 \xymatrixcolsep{1.5pc}
\let\objectstyle=\scriptstyle
\let\labelstyle=\scriptscriptstyle
  \begin{diagram}
 {T\underline{X}}\dto_-{\underline{x}}\ruppertwocell<0>^{{T\hat{L}}^{}}{<3>{}^{\hat{\gamma}}}
&{TY}\dto^-{y} \\
{\underline{X}}\rto_-{{\hat{L}}^{}}&{Y}
  \end{diagram}
\end{displaymath}
\noindent induced by 2-dimensional universal property of 
$T\underline{X}$, which furthermore guarantees the validity of the axioms making the 
above diagram a morphism in $\Alg{\mathsf{T}}_{\mathit{l}}$. The 
2-dimensional property of the Kleisli object follows similarly.

\qed
\end{proof}

\begin{proposition}
\label{prop:Kleisli-K}
If $\Bimod{\twocat{K}}$ admits representable Kleisli objects, then
        $\twocat{K}$ admits Kleisli objects. Furthermore, if $\Bimod{T}$ preserves (representable) Kleisli objects, so does $\morph{T}{\twocat{K}}{\twocat{K}}$
\end{proposition}

\begin{proof}
        Given a monad $(\morph{t}{X}{X},\eta,\mu)$ in 
        $\twocat{K}$, we transform it into a monad $(\rel{t^{*}}{X}{X},{\eta^{*}},{\mu^{*}})$
        in $\Bimod{\twocat{K}}$. We claim that the resulting representable Kleisli object
        $\underline{t^{*}}$ in $\Bimod{\twocat{K}}$ yields one for $t$ in
        $\twocat{K}$. Indeed, we obtain a lax cocone for $t$ as follows


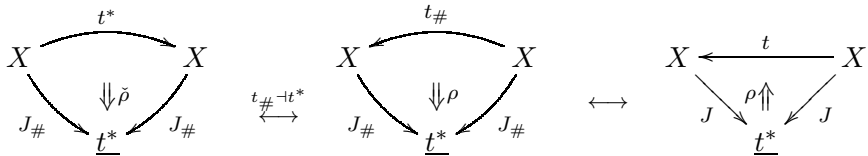
\begin{displaymath}
\xymatrixrowsep{1.5pc}
 \xymatrixcolsep{1.5pc}
  \begin{diagram}
{X}\ar@/_1ex/[dr]_-{\im{J}}
\ar@/^2ex/[rr]^{t^{*}}
\xtwocell[0,2]{}
\omit{<3>\check{\rho}}&&{X}\ar@/^1ex/[dl]^-{\im{J}}
\ar@{}[drr]|-{\stackrel{\im{t}\dashv{t^{*}}}{\longleftrightarrow}}
& &
{X}\ar@/_1ex/[dr]_-{\im{J}}&&\ar@/_2ex/[ll]_{\im{t}}
\xtwocell[0,-2]{}\omit{^<-3>\rho}{X}
\ar@/^1ex/[dl]^-{\im{J}}
\ar@{}[drr]|-{{\longleftrightarrow}}
& &
{X}\drto_-{J}&&\llto_{t}
\xtwocell[0,-2]{}\omit{<-3>\rho}{X}\dlto^-{J}
\\
& {\underline{t^{*}}}& & &&{\underline{t^{*}}}& & & &{\underline{t^{*}}}&
 \end{diagram}
\end{displaymath} 

        \noindent the latter correspondence by 
        $\morph{\im{\_}}{\twocat{K}^{\mathit{co}}}{\Bimod{\twocat{K}}}$ 
        being locally fully faithful. 
The universal property of such a representable Kleisli object 
        in $\Bimod{\twocat{K}}$
         restricts appropriately to $\twocat{K}$ thanks to the preservation of representability, \cf. Remark \ref{rmk:Kleisli-representable}.
         Preservation by $T$ is now evident, since the action of $\Bimod{T}$ on maps (representable bimodules) is simply that of $T$ on the corresponding morphisms of $\twocat{K}$.
\qed
\end{proof}

We are finally in position to state the important classification theorem 
for strong morphisms in terms of strict ones:

\begin{theorem}[Classification of strong morphisms]
\label{thm:classification-strong-morphisms}
\hfill{ }

\noindent  The inclusion $\morph{J}{\Alg{\subadj{T}}}{\pseudo{\Alg{\subadj{T}}}}$ 
has a left biadjoint 
$\morph{(\_)^{\sigma}}{\pseudo{\Alg{\subadj{T}}}}{\Alg{\subadj{T}}}$ whose 
unit is a (pseudo-natural) equivalence.
\end{theorem}
\begin{proof}
Notice first that for a given pseudo-$\subadj{T}$-algebra 
$M$, the (morphism part of the) monad 
$\morph{\zeta_M\circ{s}}{\im{(\nu_M)}}{\im{(\nu_M)}}$ corresponds to a 
morphism in $\Alg{\mathsf{T}}_{\mathit{l}}$. 
Now, our Assumption \ref{assumptionKleisliBimod} and Proposition \ref{prop:Kleisli-K} allows to apply  Proposition \ref{prop:Kleisli-lax}.
The corresponding  Kleisli 
object yields the required left biadjoint by virtue of the analysis preceeding
Lemma \ref{key-lemma}, while this latter guarantees that the unit of the biadjunction is an equivalence, as required.
\qed
\end{proof}

\begin{corollary}[Coherence for pseudo-algebras]
\label{cor:coherence}
\mbox{ }\hfill
  \begin{enumerate}
  \item Every pseudo-$\subadj{T}$-algebra is equivalent to a strict one.

  \item Every pseudo-$T$-algebra is equivalent to a strict one.
  \end{enumerate}
 \qed
\end{corollary}

\begin{remark}
Just about the only other result in the literature about coherence at this 
level of generality is that of \cite{Power89}. The main difference in terms 
of prerequisites is that we assume the good behaviour of $T$ with  
respect to bimodules, which allows its transformation into a doctrine with 
the adjoint-pseudo-algebra property. Thus, we require $T$ to be cartesian and preserve $\Hom{\_}$.
Besides that, we require $\Bimod{T}$ to 
preserve representable Kleisli objects (which implies $T$ does as well), while Power requires $T$ to preserve {\em 
bijective-on-objects\/} functors. Since 
\begin{quote}
 ` Kleisli morphism = bijective-on-objects + right adjoint '
\end{quote}
whenever these terms make sense (\eg. in any 2-category endowed with a 
Yoneda structure \cite{StreetWalters78}) and 2-functors preserve adjoints, 
our requirement is in principle formally weaker than that of Power (modulo 
our general assumption about bimodules). One advantage of our construction 
is that it usually provides an {\em explicit\/} description of the 
associated strict algebra, at least to the extent the relevant Kleisli object 
can be so described, which is certainly the case for our applications in Part \ref{part:applications}.

\end{remark}

\newpage
\part{Applications in Internal Category Theory}
\label{part:applications}

In Part \ref{part:gral-theory} we developed a theory allowing us to 
transform coherent structures into adjoint-pseudo-algebras, in the context 
of a 2-category $\twocat{K}$ admitting a calculus of bimodules and a 
2-monad $\mathsf{T}$ on it preserving bimodules, their composites and 
their Kleisli objects. We will now provide a general construction of an 
important kind of example of such a $(\twocat{K},\mathsf{T})$ pair: given 
a category with pullbacks $\cat{B}$ and a cartesian monad
$\bold{T}$ on it, the 2-category $\Cat(\cat{B})$ of internal categories 
equipped 
with the induced 2-monad $\Cat(\bold{T})$ satisfies our hypothesis. We 
proceed to illustrate the resulting theory at work in this internal 
setting by examining three important examples: monoidal categories and 
pseudo-functors (into $\Cat$) (the quintessential 
`coherent  structures') and the more novel {\em monoidal globular 
categories\/} introduced by Batanin for the development of weak 
higher-dimensional categorical structures.

\section{Bimodules for Internal Categories}
\label{sec:internal-bimodules}

In this section we consider 
$\cat{B}$ a category with pullbacks, and ${\mathbf{T}} = 
(T,\eta,\mu)$ a \textbf{cartesian monad} on it, \ie. 
        \begin{itemize}
                \item  The functor $\morph{T}{\cat{B}}{\cat{B}}$ preserves 
pullbacks
                
                \item  The transformations $\cell{\eta}{\mathit{id}}{T}$ 
and 
                $\cell{\mu}{T^2}{T}$ are cartesian, \ie. the naturality 
squares 
                are pullbacks.
        \end{itemize}

We get a cartesian 2-monad 
$\morph{\Cat({\mathbf{T}})}{\Cat(\cat{B})}{\Cat(\cat{B})}$ on the 2-category 
of internal categories, functors and natural transformations in $\cat{B}$, 
since we have a 2-functor

\[ \morph{\Cat(\_)}{\Pbk}{2-\!\Cat} \]

 \noindent where $\Pbk$ is the 2-category of categories with pullbacks, 
pullback-preserving functors and cartesian transformations.

Provisionally, we also assume $\cat{B}$ has pullback-stable coequalizers 
and that $T$ preserves them,  so 
we can apply our preceeding theory with $\twocat{K} = \Cat(\cat{B})$.  In 
particular we have an adjoint string 

  \[\xymatrixrowsep{2.5pc} \xymatrixcolsep{1.5pc}
         \let\objectstyle=\scriptstyle
        \begin{diagram}
         { \Cat(\cat{B})}\ar@/^5ex/ [d]^-{(\_)_0}
          \ar@/_5ex/ [d]_-{\pi^{0}}\\
          {\cat{B}}\ar@{-->}[u]|{\Delta}^{\dashv\;\;}_{\;\;\dashv}
               \end{diagram}
               \]

\noindent where $\Delta$ takes an object to the discrete category on it, 
and the `connected components' of an internal category $\intcat{C} =
\begin{diagram}
 {C_0}&\lto|{d}{C_1}\rto|{c} &{C_0} 
\end{diagram}
$ 
given by the coequalizer

       \[\xymatrixrowsep{2pc}
 \xymatrixcolsep{2pc}
\begin{diagram}
 {C_1}\ar@<1ex>[r]^-{d}
                                 \ar@<-1ex>[r]_-{c}&{C_0}
                                 \ar@{>>}[r]& {\pi^0(C)}
\end{diagram}
\]

For a primer on internal bimodules, \ie. bimodules in $\Cat(\cat{B})$, 
see \cite[Ch.~2]{Johnstone77}. Given categories $\intcat{C}$ and $\intcat{D}$ in 
$\cat{B}$, a bimodule $\rel{M}{\intcat{C}}{\intcat{D}}$ amounts to a span

      \[\xymatrixrowsep{1.5pc}
 \xymatrixcolsep{1.5pc}
\begin{diagram}
&\dlto_{s}{M}\drto^{t} &\\
{C_0}& & {D_0}
\end{diagram}
\]
\noindent equipped with an action 
$\morph{\alpha}{C_1\circ{M}\circ{D_1}}{M}$ in $\Spn{\cat{B}}(C_0,D_0)$, 
commuting with the monoid structure of $\intcat{C}$
and $\intcat{D}$. Alternatively, $M$ is equipped with a pair of actions
$\morph{\alpha_l}{C_1{\circ}M}{M}$ and $\morph{\alpha_r}{M{\circ}D_1}{M}$
compatible in the sense that the following diagram commutes
      \[\xymatrixrowsep{1.5pc}
 \xymatrixcolsep{1.5pc}
\begin{diagram}
{C_1{\circ}M{\circ}D_1}\dto_-{^{C_1{\circ}\alpha_r}}\rrto^-{^{\alpha_l{\circ}D_1}}&&{M{\circ}D_1}\dto^-{\alpha_r}\\
{C_1{\circ}M}\rrto_-{\alpha_l}&&{M}
\end{diagram}
\]

\begin{proposition}
\label{prop:homs}
If $T$ preserves pullbacks, $\Cat({T})$ preserves comma-objects. 
Furthermore, if $T$ preserves coequalizers, $\Cat({T})$ preserves 
bimodule composition.
\end{proposition}
\begin{proof}
        For an internal category $\intcat{C}$, the corresponding `hom' bimodule is
        $\Hom{\intcat{C}} = \begin{diagram}
 C_0&\lto|{d}C_1\rto|{c} &C_0 
\end{diagram}$ and therefore $T\Hom{\intcat{C}} = \Hom{T\intcat{C}}$. 
\qed
\end{proof}

\vskip5mm
\noindent
\underline{\textbf{Embeddings}}:  Given an internal functor
$\morph{f}{\intcat{C}}{\intcat{D}}$, its associated representable bimodule
$\rel{\im{f}}{\intcat{C}}{\intcat{D}}$ is given by the span

     \[\xymatrixrowsep{1.5pc}
 \xymatrixcolsep{1.5pc}
\begin{diagram}
{\im{f}}\ar@/^5ex/[drr]^-{t}\dto_{s}\rto_-{\bar{f}}&{D_1}\ar@{-->}[d]^{d}\drto^{c}\\
{C_0}\ar@{-->}[r]_-{^{f_0}}&{D_0}&{D_0}
\end{diagram}
\]
\noindent where the square is a pullback, and its dual 
$\rel{f^{*}}{\intcat{D}}{\intcat{C}}$ has underlying span

     \[\xymatrixrowsep{1.5pc}
 \xymatrixcolsep{1.5pc}
\begin{diagram}
&\dlto_{d}\ar@{-->}[d]_{c}{D_1}&\lto^-{{\bar{f}}}{f^{*}}\dto^{t}\ar@/_5ex/[dll]_-{s}\\
{D_0}&{D_0}&\ar@{-->}[l]^-{^{f_0}}{C_0}
\end{diagram}
\]
\noindent where the square is a pullback. In particular, for the 
cartesian transformations $\eta$ and $\mu$ we have 
   \[\xymatrixrowsep{1.5pc}
 \xymatrixcolsep{1.5pc}
\begin{diagram}
{\eta_{\intcat{C}}^{*}}\ar@{}[dr]|{=}&&\dlto_{Td}{TC_1}\ar@{-->}[d]_{Tc}&\lto_{^{\eta_{C_1}}}{C_1}\dto^{c}\\
&{TC_0}&{TC_0}&\ar@{-->}[l]^{\eta_{C_0}}{C_0}
\end{diagram}
\]
 \[\xymatrixrowsep{1.5pc}
 \xymatrixcolsep{1.5pc}
\begin{diagram}
{\mu_{\intcat{C}}^{*}}\ar@{}[dr]|{=}&&\dlto_{Td}{TC_1}\ar@{-->}[d]_{Tc}&
\lto_{^{\mu_{C_1}}}{T^2C_1}\dto^{^{T^2c}} \\
&{TC_0}&{TC_0}&\ar@{-->}[l]^{\mu_{C_0}}{T^2C_0}
\end{diagram}
\]

We want to relate the 2-category $\BimodT{\Cat(T)}{\Cat(\cat{B})}$ with 
the simpler $\SpnT{\mathsf{T}}{\cat{B}}$ of \cite{Hermida99a}, which we have
used  as a framework for  representable multicategories. There is in 
fact a  lax functor 
$\morph{\|\_\|}{\BimodT{\Cat({\mathbf{T}})}{\Cat(\cat{B})}}{\SpnT{\mathsf{T}}{\cat{B}}}$
with action

\[\xymatrixrowsep{1.5pc}
 \xymatrixcolsep{1.5pc}
\begin{diagram}
& & & & &{M}\ddto^-{\theta}\dlto_-{s}\drto^-{t} & \\
{T\intcat{C}}\rrtwocell^{M}_{M'}{\theta}& &{\intcat{D}}&\mapsto&
{TC_0}& &{D_0}\\
& & & &  &{M'}\ulto^-{s}\urto_-{t} &
 \end{diagram}
\]
\noindent with structural 2-cells  (for $\rel{N}{T\intcat{D}}{\intcat{E}}$)
\begin{description}
        \item[$\cell{\gamma_{\intcat{C}}}{\mathit{id}}{\|\eta_{\intcat{C}}^{*}\|}$:]  
          
\[\xymatrixrowsep{1.5pc}
 \xymatrixcolsep{1.5pc}
\begin{diagram}
& & \dllto_-{\eta_{C_0}}{C_0}\dto_-{\iota}\ar@{-->}[dl]|{\mathit{id}}\drto^{\mathit{id}} & \\
{TC_0}&\lto^-{\eta_{C_0}}{C_0}&\lto^-{d}{C_1}\rto_-{c}&{C_0}
\end{diagram}
\]
        
        \item[$\cell{\delta_{M,N}}{\|M\|{\bullet}\|N\|}{\|M{\bullet}N\|}$:]  

        given by the coequalizer which realizes the composition $TM{\bullet}N$ in
        $\Bimod{\Cat(\cat{B})}$

 \[\xymatrixrowsep{2pc}
 \xymatrixcolsep{1.5pc}
\begin{diagram}
 {TM{\circ}TD_1{\circ}N}\ar@<1ex>[rr]^-{^{TM{\circ}\alpha^N_l}}
                                 \ar@<-1ex>[rr]_-{^{T\alpha^M_r{\circ}N}}&&{TM{\circ}N}
                                 \ar@{>>}[rr]^-{^{\delta_{M,N}}}&& {TM\bullet{N}}
\end{diagram}
\]
\end{description}

The lax functor $\|\_\|$ therefore takes monads to monads. In fact, it 
induces a 2-functor 

\[ 
\morph{\Mnd{\|\_\|}}{\Mnd{\BimodT{\Cat(T)}{\Cat(\cat{B})}}}
{\Mnd{\SpnT{{\mathbf{T}}}{\cat{B}}}} \]

\noindent where in the target 2-category we drop the requirement of normality on the 
monads (which is not preserved by lax functors), so that 
$$\Mnd{\SpnT{{\mathbf{T}}}{\cat{B}}} \equiv \Multicat_{T}(\cat{B})$$
\noindent the 2-category of multicategories, morphisms of such and 
transformations as defined in \cite[\S 6]{Hermida99a} (where we had left 
implicit the monad $T$ under the assumption of being the free-monoid 
monad $M$, but the definition is purely formal relative to the 
cartesianness of $\mathsf{T}$). Here is where we exploit 
the description of 2-cells in $\Mnd{\BimodT{\mathsf{T}}{\twocat{K}}}$ given in Definition 
\ref{def:mnd-bimod}.

\begin{theorem}
\label{thm:bimod-spn}

        The 2-functor $$\morph{\Mnd{\|\_\|}}{\Mnd{\BimodT{\Cat(T)}{\Cat(\cat{B})}}}
{\Multicat_{T}(\cat{B})} $$ is an isomorphism
\end{theorem}
\begin{proof}
        We define an explicit inverse 
        $\morph{H}{\Multicat_{T}(\cat{B})}{\Mnd{\BimodT{\Cat(T)}{\Cat(\cat{B})}}}$ to the given 2-functor.
        
Given a multicategory $\intcat{M} = \begin{diagram}
 {TC_0}&\lto|{d}{M}\rto|{c} &{C_0} 
\end{diagram}
$ we obtain a category $\overline{\intcat{M}}$ by pulling back along 
$\eta_{C_0}$:
              \[\xymatrixrowsep{1pc}
 \xymatrixcolsep{1pc}
\begin{diagram}
 &\dlto_-{\overline{d}}{C_1}\drto\ar@/^3ex/[ddrr]^-{\overline{c}}& & \\
 {C_0}\drto_-{^\eta_{C_0}}& &\dlto_-{d}{M}\drto^-{c}& \\
 &{TC_0}& &{C_0}
\end{diagram}
\]
\noindent \cf. \cite[Def.~6.7]{Hermida99a}. We then construe $M$ itself as a bimodule
$\rel{H(M)}{T\overline{\intcat{M}}}{\overline{\intcat{M}}}$,
whose action is obtained by 
`restriction' of the composition operation of $M$. Furthermore, the monad 
structure of $M$ as a multicategory endows $H(M)$ with a monad structure
in $\BimodT{\Cat(T)}{\Cat(\cat{B})}$, which is  normal by the 
very definition of $\overline{\intcat{M}}$.

Given a morphism of multicategories 
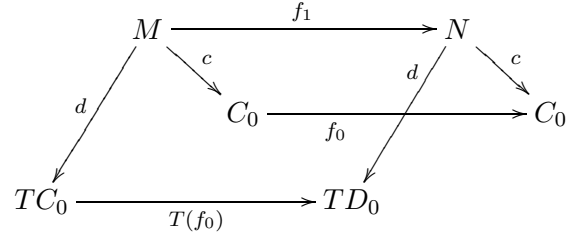
\[\xymatrixrowsep{1.5pc}
 \xymatrixcolsep{1.5pc}
\begin{diagram}
 &\ddlto_{d}M \drto^{c}\ar@{->}[rrr]^{f_1}&&& 
\ddlto_(.3){d}N\drto^{c}&\\
 &&{C_0}\ar@{->}[rrr]_(.3){f_0} 
&&&{C_0}\\
 {{T}C_0}\ar@{->}[rrr]_{{T}(f_0)}&&& {{T}D_0}&&   
\end{diagram}
\]
\noindent we obtain an internal functor 
$\morph{\overline{f}}{\overline{\intcat{M}}}{\overline{\intcat{N}}}$ (since $\overline{(\_)}$ is a 
2-functorial construction) and a morphism of bimodules 
$\morph{\hat{f_1}}{H(M)}{(Tf_0,f_0)^{*}(H(N))}$ commuting with the monads 
structures. Finally, a 2-cell
 \[\xymatrixrowsep{1.5pc}
 \xymatrixcolsep{1.5pc}
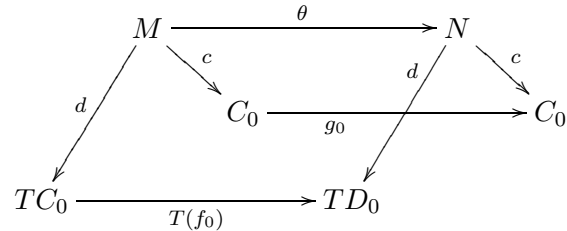
\begin{diagram}
 &\ddlto_{d}M \drto^{c}\ar@{->}[rrr]^{\theta}&&& 
\ddlto_(.3){d}N\drto^{c}&\\
 &&{C_0}\ar@{->}[rrr]_(.3){g_0} 
&&&{C_0}\\
 {{T}C_0}\ar@{->}[rrr]_{{T}(f_0)}&&& {{T}D_0}&&   
\end{diagram}
\]
\noindent gives a 2-cell $\cell{H{\theta}}{\overline{f}}{\overline{g}}$, namely 
$\morph{\hat{\theta}}{M}{(Tf_0,g_0)^{*}(N)}$, which being a morphism of 
$\bimod{(M,M)}$ is also a morphism of $\bimod{(H(M),H(M))}$.    

\qed
\end{proof}

In view of the above 2-isomorphism we can do without 
coequalizers in $\cat{B}$, and work simply with spans and 
their (pullback) composition, rather than the more complex composition of 
bimodules.

\subsection{Kleisli objects in $\Bimod{\Cat(\cat{B})}$}
\label{sub:Kleisli-int-bim}
\begin{proposition}
\label{prop:Kleisli-bimod}
$\Bimod{\Cat(\cat{B})}$ admits Kleisli objects which preserve representability (in the sense of 
Remark \ref{rmk:Kleisli-representable}).
\end{proposition}
\begin{proof}
Given a monad $(\rel{M}{\intcat{C}}{\intcat{C}},\eta,\mu)$, we obtain via 
the lax functor $\|\_\|$ a monad $(\|M\|,\|\eta\|,\|\mu\|)$ in 
$\Spn{\cat{B}}$ which is therefore an internal category. This category 
$\underline{M}$ is  the (vertex of the lax cocone of the) Kleisli 
object of $M$. To define a lax cocone

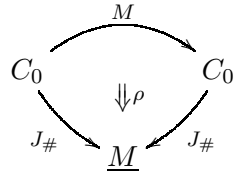
\begin{displaymath}
\xymatrixrowsep{1.5pc}
 \xymatrixcolsep{1.5pc}
  \begin{diagram}
{C_0}\ar@/_1ex/[dr]_-{\im{J}}
\rruppertwocell^{M}{<2>\rho}&&{C_0}\ar@/^1ex/[dl]^-{\im{J}}
\\
& {\underline{M}}&
 \end{diagram}
\end{displaymath}
 \noindent the 2-cell $\cell{\eta}{\Hom{\intcat{C}}}{M}$ yields an 
 identity-on-objects functor  $\morph{J}{\intcat{C}}{\underline{M}}$:
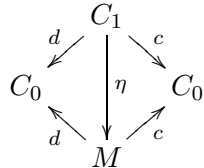
 \begin{displaymath}
\xymatrixrowsep{1pc}
 \xymatrixcolsep{1pc}
  \begin{diagram}
  & \dlto_{d}{C_1}\ddto^{\eta}\drto^{c}& \\
{C_0}&&{C_0}\\
& \ulto^{d}{M}\urto_{c}&
 \end{diagram}
\end{displaymath}
 \noindent \cf. the description of $\Hom{\intcat{C}}$ in the proof of 
 Proposition \ref{prop:homs}. The 2-cell 
 $\cell{\rho}{M{\bullet}\im{J}}{\im{J}}$ is given by $\mu$ since $\im{J}$ 
 corresponds simply to the span $\begin{diagram}
 {C_0}&\lto|{d}{M}\rto|{c} &{C_0}       
 \end{diagram}
 $, because $J$ is identity-on-objects (see the description of $\im{f}$ 
 above).
 
 For universality, given another lax cocone
 \begin{displaymath}
\xymatrixrowsep{1.5pc}
 \xymatrixcolsep{1.5pc}
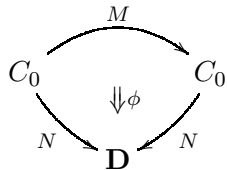
  \begin{diagram}
{C_0}\ar@/_1ex/[dr]_-{N}
\rruppertwocell^{M}{<2>\phi}&&{C_0}\ar@/^1ex/[dl]^-{N}
\\
& {\intcat{D}}&
 \end{diagram}
\end{displaymath}
 \noindent we can take the mediating bimodule 
 $\rel{\hat{N}}{\underline{M}}{\intcat{D}}$ to be $N$ itself, as we have 
 a split coequalizer:
 
     \[\xymatrixrowsep{2pc}
 \xymatrixcolsep{3pc}
  \let\objectstyle=\scriptstyle \let\labelstyle=\scriptstyle
\begin{diagram}
 {M{\circ}M{\circ}N}\ar@<1ex>[r]^-{M{\circ}\|\phi\|}
                                 \ar@<-1ex>[r]_-{\|\mu\|{\circ}N}&{M{\circ}N}
                                 \ar@/^5ex/[l]^-{\|\eta\|{\circ}M{\circ}N}
                                 \ar@{>>}[r]^-{\|\phi\|}& 
                                 \ar@/^2ex/[l]^-{\|\eta\|{\circ}N}{N}
\end{diagram}
\]

 \noindent Clearly, if 
 $N$ is representable, so is the induced $\hat{N}$ and we can force this property to be 2-functorial.
 
 \qed
\end{proof}

Given the explicit description of Kleisli objects in $\Bimod{\Cat(\cat{B})}$ 
above, we obtain the following easy consequence.

\begin{corollary}
        If $\morph{T}{\cat{B}}{\cat{B}}$ preserves pullbacks and coequalizers,
        then the homomorphism $\morph{\Bimod{T}}{\Bimod{\Cat(\cat{B})}}{\Bimod{\Cat(\cat{B})}}$ preserves Kleisli 
        objects.
\end{corollary}

The above properties about Kleisli objects in $\Bimod{\Cat(\cat{B})}$  
show that $\cat{B}$ admits a suitable 
calculus of bimodules and a cartesian monad on it induces a well-behaved 
pseudo-monad on $\Bimod{\Cat(\cat{B})}$. As we already mentioned, by virtue of Theorem
\ref{thm:bimod-spn} we can 
use spans instead of bimodules. Consequently we can dispense with 
Kleisli objects in the adjunction of 
Proposition \ref{prop:fund-adj} and use the 
simpler description of this adjunction in \cite[\S 7]{Hermida99a}. 
However, we do require well-behaved Kleisli objects in $\Cat(\cat{B})$ in 
our classification of strong morphisms and coherence (see Proposition \ref{prop:Kleisli-lax}). We obtain these 
quite simply (and elegantly) from those in $\Bimod{\Cat(\cat{B})}$ as shown in
Proposition \ref{prop:Kleisli-K}.

\begin{corollary}
\label{cor:Kleisli-Cat(B)}
\hfill{ }
\begin{itemize}
\item $\Cat(\cat{B})$ admits Kleisli objects.
\item Given 
         a pullback-preserving functor $\morph{T}{\cat{B}}{\cat{B}}$,
        the induced 2-functor \linebreak 
$\morph{\Cat(T)}{\Cat(\cat{B})}{\Cat(\cat{B})}$ 
        preserves them.
\end{itemize}
\end{corollary}

\begin{remark}
        Notice that by the construction in the proof of Proposition \ref{prop:Kleisli-K}, the internal category 
        $\underline{t^{*}}$ corresponding to a monad $(\morph{t}{\intcat{C}}{\intcat{C}},\eta,\mu)$ in 
        $\Cat(\cat{B})$ has underlying graph
         \[\xymatrixrowsep{1pc}
 \xymatrixcolsep{1pc}
\begin{diagram}
&\dlto_{d}\ar@{-->}[d]_{c}{C_1}&\lto^-{{\overline{t}}^{}}{t^{*}}\dto^{t}\ar@/_5ex/[dll]_-{s}\\
{C_0}&{C_0}&\ar@{-->}[l]^-{^{t}}{C_0}
\end{diagram}
\]
        \noindent so that $\underline{t^{*}}(x,y) = C_1(x,ty)$, just as we 
        expect from the usual construction in $\Cat(\Set)$.
\end{remark}

It is now clear that the construction of Kleisli objects in $\Cat(\cat{B})$, which we obtain from those in 
$\Bimod{\Cat(\cat{B})}$, does only involve pullbacks and no colimits. Hence,  such Kleisli objects are preserved by $\morph{\Cat(T)}{\Cat(\cat{B})}{\Cat(\cat{B})}$, for a pullback-preserving $\morph{T}{\cat{B}}{\cat{B}}$.

\section{Multicategories and monoidal categories}
\label{sec:multicat}

In this section we review our motivating example. Let $\cat{B} = \Set$ 
(we could work in far greater generality, \eg. a topos with a natural numbers object). Let
$\morph{T= (\_^{*})}{\Set}{\Set}$ be the free-monoid monad, so that $X^{*}$ can be explicitly described as the set of finite sequences of elements of $X$. This monad is cartesian \cf. \cite{Benabou90}. An algebra for the  corresponding 2-monad
$\morph{M= \Cat(\_^{*})}{\Cat}{\Cat}$ amounts to giving a strict monoidal category, an algebra morphism is a strict monoidal functor and a 2-cell a monoidal natural transformation. A pseudo-algebra is a monoidal category (with an infinite presentation) and a strong morphism is a strong monoidal functor. We do not elaborate here in the distinction between the `usual' finite presentation  of a monoidal category and its corresponding  infinite presentation (with $n$-fold tensor products for every $n$ and coherent isomorphisms between their multicategory composites) and refer to \cite[\S 9.1]{Hermida99a}. As indicated there, the distinction is inessential for our purposes.

A monad in $\SpnT{\mathsf{M}}{\Set}$ is a multicategory. We furthermore identify the corresponding morphisms and 2-cells so that 
$\Mnd{\SpnT{\mathsf{M}}{\Set}} = \Multicat$ as in \cite[\S 6]{Hermida99a}. Given a multicategory $\cat{C}$,
\[\xymatrixrowsep{1.5pc}
        \xymatrixcolsep{1.5pc}
        \begin{diagram}
        &\dlto_{d_1}{C_1}\drto^{c_1}& \\
 {\bold{M}C_0}& &{C_0}
               \end{diagram}
        \]
with $\overline{\cat{C}}$ the corresponding category of linear morphisms ($\overline{\cat{C}}(x,y) = C_1(\tuple{x},y)$), the corresponding normal monad in $\BimodT{\mathsf{M}}{\Cat}$ is given by the bimodule $\rel{\Hom{\cat{C}}}{M\overline{\cat{C}}}{\overline{\cat{C}}}$ with fibres $\Hom{\cat{C}}(\vec{x},y) = C_1(\vec{x},y)$ and action given by (multicategory) composition in
$\cat{C}$.

The 2-category $\MonCat$ of strict monoidal categories and strict monoidal functors is monadic over $\Multicat$ according to Theorem \ref{thm:monadicity}. The corresponding adjoint pseudo-algebras in $\Multicat$ are the {\em representable multicategories}. These and their correspondence with monoidal categories, as well as the relevant coherence result were analysed in detail in \cite{Hermida99a}. We skip these topics here and turn our attention to the reconstruction of $\Delta$ as the monoid classifier for monoidal categories.

\subsection{Monoid classifier}
\label{sec:monoid-classifier}

Given strict monoidal categories $\cat{C}$ and $\cat{D}$, to give a lax monoidal functor $\morph{f}{\cat{C}}{\cat{D}}$ amounts to give a strict monoidal functor $\morph{f}{FR(\cat{C})}{\cat{D}}$. A \textbf{monoid} in $\cat{C}$ amounts to a lax monoidal functor $\morph{(X,\cdot,e)}{{\mathbf 1}}{\cat{C}}$ while a monoid morphism is a 2-transformation between such. 
In particular, the identity \linebreak
$\morph{\mathit{id}}{FR({\mathbf 1})}{FR({\mathbf 1})}$ yields a monoid
$\morph{G}{{\mathbf 1}}{FR({\mathbf 1})}$
Thus we obtain the following classification:

\begin{corollary}[Monoid classifier for monoidal categories]
  \label{cor:monoid-classifier}
\mbox{ }\hfill{ }
Given a (strict) monoidal category $\cat{C}$, there is an (isomorphism) equivalence of categories 

\begin{displaymath}
  \mathsf{Monoid}(\cat{C}) \equiv (\mathsf{Ps}-)\MonCat(FR({\mathbf 1}),\cat{C})
\end{displaymath}
\noindent realised (from right to left) by precomposition with the \underline{generic monoid} 
$$\morph{G}{{\mathbf 1}}{FR({\mathbf 1})}$$.
\end{corollary}

Let us give an explicit description of $FR({\mathbf 1})$. First $R({\mathbf 1})$ is the terminal multicategory, whose underlying multigraph can be identified with 

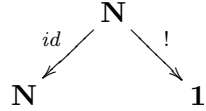
\[\xymatrixrowsep{1.5pc}
        \xymatrixcolsep{1.5pc}
        \begin{diagram}
        &\dlto_{\mathit{id}}{\bold{N}}\drto^{!}& \\
 {\bold{N}}& &{{\mathbf 1}}
               \end{diagram}
        \]

\noindent where $\bold{N}$ is the set of natural numbers. Thus we have a unique arrow $n \succ \bullet$ for every $n$. Now $FR({\mathbf 1})$ has the following underlying graph

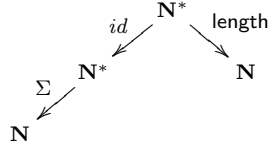
\[\xymatrixrowsep{1pc}
        \xymatrixcolsep{1pc}
\let\objectstyle=\scriptstyle
\let\labelstyle=\scriptstyle
\begin{diagram}
&&\dlto_{\mathit{id}}{\bold{N}^{*}}\drto^{\mathsf{length}}& \\
& {\bold{N}^{*}}\dlto_{\Sigma}& &{\bold{N}}\\
{\bold{N}}& & & \\ 
              \end{diagram}
        \]

\noindent Hence the objects of $FR({\mathbf 1})$ are natural numbers. To give an arrow between $n$ and $m$, we must give a `partition' of $n = n_1 + \ldots + n_m$, which corresponds to the arrow 
$$\morph{\tuple{(n_1\succ \bullet),\ldots,(n_m\succ \bullet)}}{n}{m}$$
\noindent Considering $n$ and $m$ as finite ordinals, such an arrow is then a {\em monotone function\/} from $[n]$ to $[m]$. Conversely, a monotone function $\morph{h}{[n]}{[m]}$ yields a partition of $[n] = h^{-1}(1) + \ldots + h^{-1}(m)$ via its (possibly empty) fibers; the sum is now interpreted as ordinal sum. 

\begin{theorem}
\label{thm:delta}
\mbox{ }\hfill
\begin{center}
\framebox[1.1\width]{
$FR({\mathbf 1}) \equiv \Delta$
}
\end{center}
\noindent where $\Delta$ is the category of finite ordinals and monotone functions, with strict monoidal structure given by ordinal sum. The generic monoid is $1\in \Delta$ with the unique morphisms $\morph{!}{0}{1}$ and $\morph{!}{2}{1}$.
\end{theorem}

We have thus recovered the classical description of the monoid classifier for monoidal categories, usually credited to Lawvere.


\section{Monoidal globular categories}
\label{sec:monoidal-globular-categories}

Monoidal globular categories were introduced in \cite{Batanin98} as a 
framework for weak higher-dimensional categories. Briefly put, they allow the definition of higher-dimensional
{\em operads\/} (see Remark \ref{rmk:operads} below) so that weak $n$-categories are defined as the algebras for a {\em contractible\/} operad in the monoidal globular category of spans. 

In our brief incursion into monoidal globular categories we will exhibit them as pseudo-algebras on {\em globular categories} (\sectref{sec:pseudo-algebras-globular-categories}), give an equivalent adjoint-pseudo-algebra version, namely {\em representable multicategories in the category $\wgph$ of $\omega$-graphs\/} 
and formulate the attendant coherence result (\sectref{sec:rep-mult-wgph})
. In \sectref{sec:globular-monoids} we analyse lax morphisms from the terminal monoidal globular category into a given monoidal globular category $\cat{C}$ (the {\em globular monoids\/} in $\cat{C}$ introduced in  \cite{BataninStreet98}) and give a description of their classifier.

\subsection{Pseudo-algebras on globular categories}
\label{sec:pseudo-algebras-globular-categories}

An $n$-graph $C$ is a commutative diagram of sets 

\[\xymatrixrowsep{1.5pc}
 \xymatrixcolsep{1.5pc}
\begin{diagram}
&\dlto_{d_n}{C_n}\drto^{c_n}&\\
 \dto_{d_{n-1}}{C_{n-1}}\drrto^(.3){c_{n-1}}|!{[d];[rr]}\hole& 
&\dllto_(.3){d_{n-1}}{C_{n-1}}\dto^{c_{n-1}}\\
{\mbox{ }}\ar@{.}[d]& &{\mbox{ }}\ar@{.}[d] \\
 {C_0}& & {C_0}
 \end{diagram}
\]

\noindent which means that the following {\em globularity\/} condition is 
satisfied:
\begin{eqnarray*}
        c_{i-1}c_i & = & c_{i-1}d_i  \\
        d_{i-1}c_i & = & d_{i-1}d_i \qquad 2\leq{i}\leq{n}
\end{eqnarray*}

\noindent This implies that we have well defined domain and codomain 
functions from any dimension $i$ to a lower one $j$, 
$\morph{d_i^{j}}{C_i}{C_j}$ and $\morph{c_i^{j}}{C_i}{C_j}$ by iterated 
composition of $d$'s and $c$'s respectively. We refer to the elements of 
$C_i$ as $i$-{\em cells\/}.

A morphism between $n$-graphs $C$ and $D$ is a collection of functions 
$\morph{f_i}{C_i}{D_i}$ commuting with $d_i$ and $c_i$. We thus have the 
category $\ngph{n}$.
Every $n$-category has an underlying $n$-graph, which yields the forgetful 
functor \linebreak
$\morph{U_n}{\ncat{n}}{\ngph{n}}$ from the category of 
$n$-categories and $n$-functors. This functor is monadic.

We want to consider the monad $\morph{T_\omega}{\wgph}{\wgph}$ on 
$\omega$-graphs induced by the monadic adjunction $\morph{F_\omega\dashv 
U_\omega}{\wcat}{\wgph}$. To analyse its behaviour, we look at its `finite 
approximations' in the filtration (colimit sequence):

\[\xymatrixrowsep{2pc}
 \xymatrixcolsep{2pc}
\begin{diagram}
{\ncat{1}}\ar@/^1ex/[d]^{U_1}\ar@{^{(}->}[r]^-{J_1}&{\ncat{2}}\ar@/^1ex/[d]^{U_2}\ar@{^{(}->}[r]^-{J_2}&{\ldots}&{\wcat}\ar@/^1ex/[d]^{U_\omega}\\
{\ngph{1}}\ar@/^1ex/[u]^{F_1}_{\dashv}\ar@{^{(}->}[r]^-{H_1}&
{\ngph{2}}\ar@/^1ex/[u]^{F_2}_{\dashv}\ar@{^{(}->}[r]^-{H_2}&{\ldots}&
{\wgph}\ar@/^1ex/[u]^{F_\omega}_{\dashv}
 \end{diagram}
\]

\noindent where $\morph{H_n}{\ngph{n}}{\ngph{n+1}}$ considers an $n$-graph 
as an $n+1$-graph with no $n+1$-cells ($C_{n+1} = \emptyset$) and 
$\morph{J_n}{\ncat{n}}{\ncat{n+1}}$ adds only identity $n+1$-cells to an 
$n$-category ($C_{n+1} = C_{n}$). All the categories $\ngph{n}$ and $\wgph$ are presheaf 
categories, and hence have dimension-wise limits and coequalizers. The 
monad $\morph{T_n}{\ngph{n}}{\ngph{n}}$ is cartesian, very much like the 
free-category monad on a graph which is a free-monoid construction (see 
\eg. \cite[\S 11]{Hermida99a}) and so is $T_\omega$ as indicated in 
\cite{Street99b,BataninStreet98}. Consequently we have the 2-monads 
$\morph{\Cat(T_1)}{\Cat(\ngph{1})}{\Cat(\ngph{1})}$ and 
$\morph{\Cat(T_\omega)}{\Cat(\wgph)}{\Cat(\wgph)}$. Notice that $\wgph = 
\mathsf{Glob}$, the category of globular sets in
\cite{Street99b}, and thus $\Cat(\wgph) = \mathsf{GlobCat}$, the category 
of globular categories and globular functors. We want to analyse the 
pseudo-algebras of $\Cat(T_\omega)$. We do so by analysing the algebras 
and pseudo-algebras of $\Cat(T_2)$.

An  object in $\Cat(\ngph{2})$ is a $2$-graph  of categories

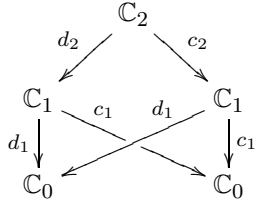
\[\xymatrixrowsep{1.5pc}
 \xymatrixcolsep{1.5pc}
\begin{diagram}
&\dlto_{d_2}{\cat{C}_2}\drto^{c_2}&\\
 \dto_{d_1}{\cat{C}_1}\drrto^(.3){c_1}|!{[d];[rr]}\hole& 
&\dllto_(.3){d_1}{\cat{C}_1}\dto^{c_1}\\
 {\cat{C}_0}& & {\cat{C}_0}
 \end{diagram}
\]

\noindent To simplify the analysis, let us take $\cat{C}_0 = \mathbf{1}$, 
the terminal category, so that we are left with a graph of categories. A 
$\Cat(T_2)$-algebra structure on this graph consists of strict monoidal 
structures $(\cat{C}_1,\tensor_1,I_1)$ and 
$(\cat{C}_2,\tensor_2^{0},I_2^{0})$, with $d_2$ and $c_2$ strict monoidal 
functors, and furthermore the span 
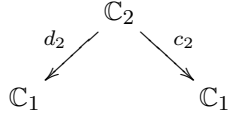
\[\xymatrixrowsep{1.5pc}
 \xymatrixcolsep{1.5pc}
\begin{diagram}
&\dlto_{d_2}{\cat{C}_2}\drto^{c_2}&\\
 {\cat{C}_1}
&  &{\cat{C}_1}
 \end{diagram}
\]
\noindent
has a monoid structure $(\cat{C}_2,\tensor_2^{1},I_2^{1})$ in 
$\Spn{\Cat}(\cat{C}_1,\cat{C}_1)$ which commutes with the other monoidal 
structure on $\cat{C}_2$. Let us spell out this last condition: given an 
object $a$ in $\cat{C}_2$, write
$
\begin{diagram}
x\ar@/^/[r]^-{a}&{y}  
\end{diagram}
$  to indicate $d_2(a) = x$ and $c_2(a) = y$.

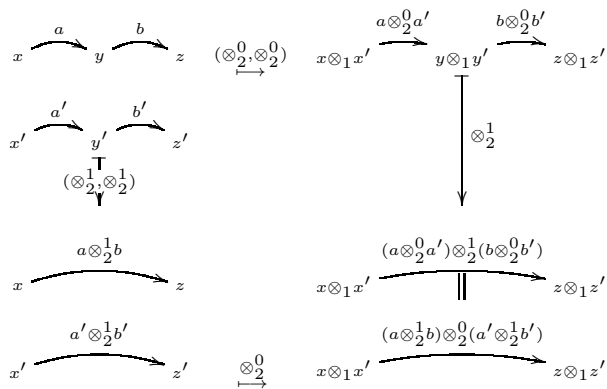
\[\xymatrixrowsep{1.5pc}
 \xymatrixcolsep{1.5pc}
\let\objectstyle=\scriptscriptstyle
\let\labelstyle=\scriptscriptstyle
\begin{diagram}
x\ar@/^/[r]^-{a}&{y}\ar@/^/[r]^-{b}&{z}\ar@{}[rr]|-{\stackrel{(\tensor_2^{0},\tensor_2^{0})}{\longmapsto}}&&
{x{\tensor_1}x'}\ar@/^/[r]^-{a{\tensor_2^{0}}a'}&{y{\tensor_1}y'}\ar@{|->}[dd]^{\tensor_2^{1}}\ar@/^/[r]^-{b{\tensor_2^{0}}b'}&{z{\tensor_1}z'}\\
x'\ar@/^/[r]^-{a'}&{y'}\ar@{|->}[d]|-{(\tensor_2^{1},\tensor_2^{1})}\ar@/^/[r]^-{b'}&{z'}& 
& &\\
& &{ } & & & & \\
x\ar@/^/[rr]^-{a{\tensor_2^{1}}b}&&{z}& &{x{\tensor_1}x'}\rrtwocell\omit{=}
\ar@/^/[rr]^-{(a{\tensor_2^{0}}a')\tensor_2^{1}(b{\tensor_2^{0}}b')}&&{z{\tensor_1}z'}\\
x'\ar@/^/[rr]^-{a'{\tensor_2^{1}}b'}&&{z'}\ar@{}[rr]|-{\stackrel{\tensor_2^{0}}{\longmapsto}}
&
&{x{\tensor_1}x'}\ar@/^/[rr]^-{(a{\tensor_2^{1}}b)\tensor_2^{0} 
(a'{\tensor_2^{1}}b')}&& {z{\tensor_1}z'}
 \end{diagram}
\]

\noindent Now to give a pseudo-$\Cat(T_2)$-algebra structure on the same 
graph we weaken all the monoids to pseudo-monoids (but $d_2$ and $c_2$ 
remain {\em strict\/} monoidal morphisms with respect to $\tensor_2^{0}$ 
and $\tensor_1$) and $\tensor_2^0$ and $\tensor_2^{1}$ pseudo-commute, in 
the sense that the equality in the above chasing-tensors diagram becomes 
an isomorphism 
$\isomorph{\gamma_{a,b}^{a',b'}}{(a{\tensor_2^{1}}b)\tensor_2^{0} 
(a'{\tensor_2^{1}}b')}
{(a{\tensor_2^{0}}a')\tensor_2^{1}(b{\tensor_2^{0}}b')}$ coherent with 
respect to the pseudo-monoidal structures involved. This amounts to the 
fact that $\tensor_2^{1}$ is a strong monoidal functor 
$\cat{C}_2\times_{\cat{C}_1}\cat{C}_2$ (which inherits a pairwise monoidal 
structure from 
$(\cat{C}_2,\tensor_2^{0},I_2^{0},\alpha,\lambda,\rho)$) to $\cat{C}_2$. 
Simliar considerations apply to the units. This is the point of view 
adopted in \cite{BalteanuFiedorowiczSchwanzlVogt98}, which gives a sound 
conceptual account of interchange constraints and their attendant axioms. 
Extrapolating to the colimit, we conclude

\vskip1em
\begin{center}
\framebox[1.1\width]{$\Alg{\Cat(T_\omega)}\equiv 
{\mbox{\textit{Strict-Monoidal-Globular-Cat}}}
$}
\end{center}
\vskip1em
\begin{center}
  \framebox[1.1\width]{$\pseudo{\Alg{\Cat(\bold{M})}}\equiv
    {\mbox{\textit{Monoidal-Globular-Cat}}} $}
\end{center}
\vskip5mm
where the 2-categories on the right give the `infinitary' presentations of the objects of the corresponding ones in \cite{Batanin98}.


\subsection{Representable multicategories in $\wgph$}
\label{sec:rep-mult-wgph}

Consider an object in $\SpnT{{\mathsf T}_2}{\ngph{2}}$: it consists of two 2-graphs $O$ and $A$:

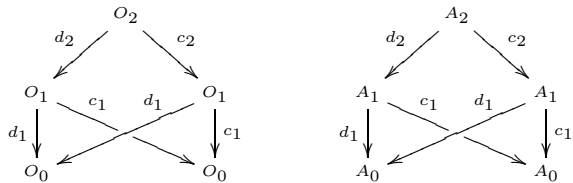
\[\xymatrixrowsep{1.5pc}
 \xymatrixcolsep{1.5pc}
\let\objectstyle=\scriptscriptstyle
\let\labelstyle=\scriptscriptstyle
\begin{diagram}
&\dlto_{d_2}{O_2}\drto^{c_2}&
& & &\dlto_{d_2}{A_2}\drto^{c_2}&
\\
 \dto_{d_1}{O_1}\drrto^(.3){c_1}|!{[d];[rr]}\hole& 
&\dllto_(.3){d_1}{O_1}\dto^{c_1}
&&  \dto_{d_1}{A_1}\drrto^(.3){c_1}|!{[d];[rr]}\hole& 
&\dllto_(.3){d_1}{A_1}\dto^{c_1}
\\
 {O_0}& & {O_0}
&&
{A_0}& & {A_0}
\end{diagram}
\]

\noindent and a span 

\[\xymatrixrowsep{1.5pc}
 \xymatrixcolsep{1.5pc}
\let\objectstyle=\scriptscriptstyle
\let\labelstyle=\scriptscriptstyle
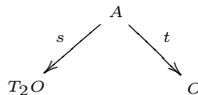
\begin{diagram}
&\dlto_{s}{A}\drto^{t}& \\
{T_2O}& & {O}
\end{diagram}
\]
in $\ngph{2}$ which in turn consists of 3 ordinary spans 

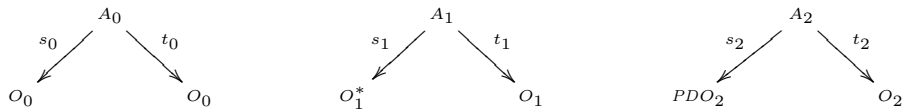
\[\xymatrixrowsep{1.5pc}
 \xymatrixcolsep{1.5pc}
\let\objectstyle=\scriptscriptstyle
\let\labelstyle=\scriptscriptstyle
\begin{diagram}
&\dlto_{s_0}{A_0}\drto^{t_0}&
& &
&\dlto_{s_1}{A_1}\drto^{t_1}&
& &&\dlto_{s_2}{A_2}\drto^{t_2}& \\
{O_0}& & {O_0}
&& 
{O_1^{*}}& & {O_1}
&&{{\mathit PD}O_2}& & {O_2}
\end{diagram}
\]

\noindent compatible with $d_i,c_i (i=0,1,2)$,  where $O_1^{*}$ is the set of composable 1-cells in $O$ and ${\mathit PD}O_2$ is the set of pastings of 2-cells in $O$. See \cite{Street96} for a detailed description of the  pasting diagrams of the related (and more complex) situation of computads and their role in the construction of free 2-categories.

A {\em multicategory\/} in $\ngph{2}$ is therefore such a span endowed with a monoid structure. This means that the first span above sets up a category with object of objects $O_0$ and object of morphisms $A_0$. The second span is then a multicategory as follows: let us denote the elements  of $A_i$ by symbols $a^{i}_x$ and adopt a similar convention for the elements of $O_i$. The combinatorial information of source, target, domain and codomain for an element of $A_1$ can be represented diagramatically as a 2-cell:
\[\xymatrixrowsep{1pc}
 \xymatrixcolsep{1.5pc}
\let\objectstyle=\scriptscriptstyle
\let\labelstyle=\scriptscriptstyle
\begin{diagram}
{o^{0}_{X_1}}\dto_-{o^1_{s_1}}
\rto^{a^{0}_d}&{o^0_X}\ar@{->}[ddd]^-{o^{1}_t}\\
{o^0_{Y_1}}\ar@{.}[d]
& & \\
{o^0_{X_{n}}}\dto_-{o^1_{s_n}}& & \\
{o^0_{Y_n}}\xtwocell[-3,0]{}\omit{<4>{a^1_\alpha}^{}}
\rto_-{a^0_c}&{o^{0}_Y}
\end{diagram}
\]
\noindent to indicate 

\begin{eqnarray*}
  d_1(a^1_\alpha) & = & a^0_d \\
c_1(a^1_\alpha) & = & a^0_c \\
s_1(a^1_\alpha) & = & \tuple{o^1_{s_1},\ldots,o^1_{s_n}} \\
t_1(a^1_\alpha) & = & o^1_t
\end{eqnarray*}
\noindent and the further evident information about $0$-domain and codomains.
These cells are equipped with a composition operation

\[\xymatrixrowsep{1.5pc}
 \xymatrixcolsep{2.5pc}
\let\objectstyle=\scriptstyle
\let\labelstyle=\scriptscriptstyle
\begin{diagram}
{o^0_{X_{11}}}\ar@{}[ddr]|-{a^1_{\alpha_1}}
\rto^-{{a^0_{d_1}}^{}}
\dto_-{o^1_{s_{11}}}&{o^{0}_{X_1}}\ddto_-{o^1_{s_1}}
\ar@{}[dddddr]|-{a^1_{\alpha}}
\rto^{a^{0}_d}&{o^0_X}
\ar@{->}[ddddd]^-{o^{1}_t}
& &{o^0_{X_{11}}}\dto_-{o^1_{s_{11}}}
\ar@{}[dddddr]|-{a^1_{\tiny \alpha(\alpha_1,\ldots,\alpha_n)}}
\rto^-{{a^0_d a^0_{d_1}}^{}}&{o^0_X}\ar@{->}[ddddd]^-{o^{1}_t}\\
{o^0_{X_{12}}}\ar@{.}[d]& & & & {o^0_{X_{12}}}\ar@{.}[ddd]& & \\
{o^0_{X_{1m}}}\ar@{.}[d]
\rto_-{a^0_{c_1}}&{o^0_{Y_1}}\ar@{.}[d]
& & 
{\longmapsto} & & & \\
{o^{0}_{X_{n1}}}
\ar@{.}[d]\ar@{}[ddr]|{a^1_{\alpha_n}}
\rto^-{a^0_{d_n}}&{o^{0}_{X_{n}}}\ddto_-{o^1_{s_n}}
 & & 
& &
& & && \\
{o^{0}_{X_{nm'}}}\dto_-{o^1_{s_{nm'}}}& & 
& &
{o^{0}_{X_{nm'}}}\dto_-{o^1_{s_{nm'}}}& & \\
{o^{0}_{X_{nm''}}}\rto_{a^0_{c_n}}&{o^0_{X_n}}\rto_-{a^0_c}&{o^0_Y}& &
{o^{0}_{X_{nm''}}}\rto_{{a^0_c a^0_{c_n}}^{}}&{o^0_Y}
\end{diagram}
\]
\noindent associative and unitary as expected. A similar `multicategory' structure is provided for the third span above, this time with the more complex pasting composition, which we will use below to account for interchange.

Let us now consider {\em representability\/} for these multicategories in $\ngph{2}$. A cell $\cell{a^1_\pi}{\tuple{o^1_{s_1},\ldots,o^1_{s_n}}}{o^1_{\tensor\tuple{s_1,\ldots,s_n}}}$  in $A_1$ is {\em universal\/} if precomposition with it sets up a bijection
\[\xymatrixrowsep{1pc}
 \xymatrixcolsep{2pc}
\let\objectstyle=\scriptscriptstyle
\let\labelstyle=\scriptscriptstyle
\begin{diagram}
{o^{0}_{X_1}}\dto_-{o^1_{s_1}}
\ar@{}[dddr]|{a^1_\alpha}
\rto^{a^{0}_d}&{o^0_X}\ar@{->}[ddd]^-{o^{1}_t}
& & & {o^0_S}\ar@{->}[ddd]_{{o^{1}_{\tensor\tuple{s_1,\ldots,s_n}}}^{}}
\ar@{}[dddr]|{\widehat{a}^1_{\alpha}}
\rto^{\widehat{a}^0_{\sigma}}&{o^0_X}\ar@{->}[ddd]^-{o^{1}_t}
\\
{o^0_{Y_1}}\ar@{.}[d]
& & {\longleftrightarrow} & &\\
{o^0_{X_n}}\dto_-{o^1_{s_n}}& && &  \\
{o^0_{Y_n}}\rto_-{a^0_c}&{o^{0}_T}
& & &{o^0_{Y_n}}\rto_-{\widehat{a}^0_\tau}&{o^0_Y}
\end{diagram}
\]
\noindent natural in $o^1_{t}$. Furthermore, such universal cells should be closed under composition
and (quite importantly) $\morph{d_i,c_i}{A_i}{A_{i-1}}$ should 
preserve universal cells.
 To see how the interchange isomorphism arises with this reformulation, consider again for simplicity that the $0$-level is trivial and look once again at the situation we considered before. We can realise the relevant tensor composites by universal cells with the appropriate sources.
Both ways around, the mutlicategory composite cells are universal (since these are closed under such composition) in $A_2$ for the same source element in ${\mathit PD}O_2$, hence their targets are canonically isomorphic.  

Thus, a \textbf{multicategory in} $\wgph$ consists of a multigraph in $\wgph$ so that the source of a $n$-cell (an element in $A_n$) is a pasting diagram of elements of $O_n$. Such a multigraph is equipped with unitary and associative multicategory composites as explained above. It is \textbf{representable} when for every pasting diagram $p$ in $O_n$ there is a universal $n$-cell whose source is $p$, and such universal cells are closed under composition and preserved by $\morph{d_i,c_i}{A_i}{A_{i-1}}$ for all $i > 1$. The morphisms between such are the evident ones, preserving the combinatorial source-target information, composites and identites and (for representables) preserving universal cells. 

By Theorems \ref{thm:monadicity} and
\ref{thm:bimod-spn} we have

\begin{corollary}[Universal version of monoidal globular categories]
\label{cor:rep-w-multicat}
\mbox{ }\hfill{ }
\begin{itemize}
\item \framebox[1.1\width]
{$\Repmulticat_s(\wgph) \equiv \Alg{T_{\omega}}
\equiv
   { \mbox{\textit{Strict-Monoidal-Globular-Cat}}}$}

\noindent where $\Repmulticat_s(\wgph)$ denotes the 2-category of strict representable multicategories and strict morphisms of such.
        
\item \framebox[1.1\width]{$\Repmulticat(\wgph) \equiv 
\pseudo{\Alg{T_{\omega}}}
    \equiv {\mbox{\textit{Monoidal-Globular-Cat}}}$}
\end{itemize}
\end{corollary}

As a further consequence of our general theory we can establish the following coherence results:

\begin{corollary}[Coherence for monoidal globular categories]
\label{coherence-monoidal-glob}
        
\mbox{ }\hfill{ }
\begin{itemize}

\item The inclusion
        
  $$\Repmulticat_s(\wgph)\hookrightarrow
\Repmulticat(\wgph) 
$$
        
        \noindent has a left biadjoint whose unit is a (pseudo-natural) equivalence
\item The inclusion
        
  $$
  \mathit{Strict Monoidal Globular Cat}\hookrightarrow
  \mathit{Monoidal Globular Cat}$$
        
        \noindent has a left biadjoint whose unit is a (pseudo-natural) equivalence
      \end{itemize}
\end{corollary}

We have therefore recovered the coherence result for monoidal globular categories by methods altogether different to those in
\cite[\S 4]{Batanin98}. Notice that Thm.~4.1 in \ibid. follows from the above and the fact that a strong morphism out of a free monoidal globular category can be strictified (\cf. \cite[Prop.10.3]{Hermida99a}).

\begin{remark}\label{rmk:operads}
A multicategory in $\wgph$ with object of objects the terminal $\omega$-graph $\bold{1}$ amounts to an {\em operad\/} in 
$\mathit{Span\/}$ in the sense \cite{Batanin98}
  
\end{remark}

\subsection{Globular monoids and their classifier}
\label{sec:globular-monoids}

The terminal monoidal globular category $\bold{1}$ has underlying $\omega$-graph the terminal such, that is, the one with only one cell at every dimension (denoted $U_\omega$ in \cite{Batanin98}). A lax morphism from $\bold{1}$ into 
a monoidal globular category $\cat{C}$ amounts to what Batanin and Street call a {\em globular monoid\/} in $\cat{C}$ (\cf. \cite{BataninStreet98}). In elementary terms, a globular monoid $(M,\iota,\mu)$ in $\cat{C}$ is given by:
\begin{itemize}
\item objects $M_n$ of $\cat{C}_n$ for every $n$, compatible with domain and codomain in $\cat{C}$ ($d_i(M_i) = c_i(M_i) = M_{i-1}$),
\item morphisms $\morph{\iota_n}{I_n(M_n)}{M_n}$ and $\morph{\mu_n}{M_n{\tensor_n}M_n}{M_n}$ which make $(M_n,\iota_n,\mu_n)$ a monoid in the monoidal category $(\cat{C}_n,I_n(M_n),\tensor_n)$,

\item the monoid structure commutes with the interchange isomorphisms: for every $i < j < n$,
  \begin{eqnarray*}
    \mu_j{\circ}(\mu_i{\tensor_j}\mu_i){\circ}\gamma^{M_n,M_n}_{M_n,M_n} & = & \mu_j{\circ}(\mu_j{\tensor_i}\mu_j) \\
\mu_i{\circ}(\iota_j{\tensor_i}\iota_j) & = & \isomorph{\lambda}{I_j(M_n){\tensor_i}I_j(M_n)}{I_j(M_n)}
  \end{eqnarray*}

\end{itemize}

According to our classification of lax morphisms of \sectref{sec:lax-morph}, or rather its more explicit  version in terms of multicategories in \sectref{sec:monoid-classifier}, such a globular monoid corresponds to a (strict) morphism of (strict) monoidal globular categories $FR({\bold{1}})\rightarrow\cat{C}$. 

\begin{corollary}[Classification of globular monoids]
\label{cor:globular-monoid-classifier}
  Given a (strict) monoidal globular category $\cat{C}$, there is an (isomorphism) equivalence of categories 

\begin{displaymath}
  \mathsf{GlobularMonoid}(\cat{C}) \equiv (\mathsf{Ps}-)\mathit{Strict Monoidal Globular Cat}(FR({\mathbf 1}),\cat{C})
\end{displaymath}
\noindent realised (from right to left) by precomposition with the \underline{generic monoid} 
$\morph{G}{{\mathbf 1}}{FR({\mathbf 1})}$.
\end{corollary}

To obtain an explicit description of the globular monoid classifier
$FR(\bold{1})$ we need to know the free monoidal globular category on $R\bold{1}$. To describe $R\bold{1}$ in turn we must know $T_{\omega}(\bold{1})$, the free $\omega$-category on the terminal $\omega$-graph. This object has been described by Batanin \cite{Batanin98,BataninStreet98} in terms of trees. The $n$-cells $\omega$-category $\mathsf{Tree}$  are finite trees of height $n$
(such a tree is formally described as a functor $\morph{\tau}{[n]^{\mathit{op}}}{\Delta}$ such that $\tau(0) = [0]$ in \ibid.) which give the relevant combinatorial information of the pasting diagrams they represent:

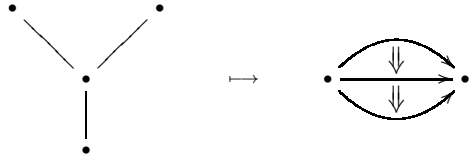
\[\xymatrixrowsep{1.5pc}
 \xymatrixcolsep{1.5pc}
\let\objectstyle=\scriptscriptstyle
\let\labelstyle=\scriptscriptstyle
\begin{diagram}
{\bullet}\ar@{-}[dr]& &\ar@{-}[dl]{\bullet}& &  & \\
&{\bullet}\ar@{-}[d]& &{\longmapsto}& {\bullet}\rruppertwocell<6>{<-1.5>}\rrto\rrlowertwocell<-6>{<1.5>}&&{\bullet} \\
&{\bullet}& & & & &
\end{diagram}
\]

\noindent There is in fact an explicit construction associating to a tree $\tau$ the $\omega$-graph (globular set) $\|\tau\|$
it 
represents  \cf. \cite{Batanin98,Street99b}.

The (strict) monoidal globular category $FR(\bold{1})$ has underlying graph

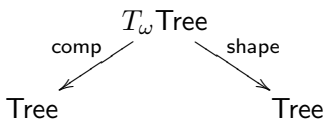
\[\xymatrixrowsep{1.5pc}
        \xymatrixcolsep{1.5pc}
        \begin{diagram}
        &\dlto_{\mathsf{comp}}{T_\omega{{\mathsf Tree}}}\drto^{\mathsf{shape}}& \\
 {\mathsf Tree}& &{{\mathsf Tree}}
               \end{diagram}
        \]

\noindent where the cells of $T_\omega{{\mathsf Tree}}$ are trees labelled by trees. A {\em labelling\/} of a tree is formally defined via its associated $\omega$-graph: a labelling of $\tau$ is a morphism in  $\morph{l}{\|\tau\|}{\mathsf{Tree}}$ in $\wgph$. The composition operation $\morph{{\mathsf comp}}{T_\omega{{\mathsf Tree}}}{\mathsf{Tree}}$ amounts to grafting the trees in the labelling according to the shape of the labelling tree, or equivalently  performing the pasting composite of the associated $\omega$-graphs. We thus obtain the monoidal globular category $\Omega$ of \cite{BataninStreet98}, together with its generic globular monoid $U_\omega$ (with its unique globular monoid structure).

\section{Pseudo-functors}
\label{sec:pseudo-functors}

Given a small category $\cat{C}$ with set of objects $C$, consider the 
functor 2-category $\twocat{K} = [C,\Cat]$. Clearly $\twocat{K}$ admits a 
calculus of bimodules and Kleisli objects for them, with all the relevant 
structure given pointwise. The category $\cat{C}$ acts on $\twocat{K}$ 
defining a 2-monad $\morph{\cat{C}\star\_}{\twocat{K}}{\twocat{K}}$
which can be described explicitly as follows:

\begin{description}
        \item[functor]  Given a functor $\morph{F}{C}{\Cat}$, the functor 
        $\morph{\cat{C}\star{F}}{C}{\Cat}$ has action
        \[  \cat{C}\star{F}(x) = \coprod_{\morph{f}{y}{x}} F(y)\]
        \noindent and we write  $\tuple{f,\varphi}$ for an object of 
        $\cat{C}\star{F}(x)$, with $\varphi \in F(\mathsf{dom} f)$.
        
        \item[unit]  $\cell{\eta_F}{F}{\cat{C}\star{F}}$ takes the object 
        $\varphi \in Fx$ to $\tuple{\mathit{id}_x,\varphi}$.
        
        \item[multiplication]  
        $\cell{\mu_F}{\cat{C}\star\cat{C}\star{F}(x)}{\cat{C}\star{F}(x)}$ 
        takes  the object $\tuple{f,\tuple{g,\varphi}}$ to 
        $\tuple{f\circ{g},\varphi}$.
\end{description}

\begin{remark}
The action of $\cat{C}$ on $\twocat{K}$ can be grasped more easily by 
identifying $\twocat{K}$ with $\Cat/C$, by regarding $\morph{F}{C}{\Cat}$ 
as the family $\coprod(F)\rightarrow C: (x,\varphi) \mapsto x$. Given the 
small category $\cat{C}$:
{\xymatrixrowsep{1.5pc}
 \xymatrixcolsep{1.5pc}
               \begin{diagram}
   {C}&\lto_{d}{C_1}\rto^{c} & {C}
\end{diagram}
}
 the family corresponding  to $\cat{C}\star{F}$ is obtained by pullback 
 against $d$, followed by composition with $c$:

 \[\xymatrixrowsep{1.5pc}
 \xymatrixcolsep{1.5pc}
 \let\objectstyle=\scriptstyle
               \begin{diagram}
              {\coprod(F)}\dto_{p}&\lto {d^{*}(\coprod(F))}\dto
              \drto^{`\cat{C}\star{F}'}&\\
   {C}&\lto_{d}{C_1}\rto^{c} & {C}
\end{diagram}
\]
\noindent With this point of view we can easily show that $\cat{C}\star$ is cartesian, along the lines of \cite{Benabou90}.
\end{remark}

\vskip1em
\noindent
\underline{$\cat{C}\star$-\textbf{algebras and pseudo-algebras}:} Let 
$\morph{\alpha}{\cat{C}\star{F}}{F}$ be a $\cat{C}\star$-algebra. This 
amounts to a $C$-collection of functors 
$\morph{\alpha_x}{\coprod_{\morph{f}{y}{x}} F(y)}{F(x)}$. It is clear 
that such $\alpha_x$'s satisfying the algebra axioms, give the action on 
morphisms to endow $F$ with the structure of a functor 
$\morph{F}{\cat{C}}{\Cat}$, namely $$\morph{F(\morph{f}{y}{x}) = 
\alpha_x(\tuple{f,\_})}{F(y)}{F(x)}$$

\noindent and a further easy identification of 1- and 2-cells yields

\begin{center}
  \framebox[1.1\width]{$\Alg{\cat{C}\star}\equiv [\cat{C},\Cat]$}
\end{center}

\noindent and similarly

\begin{center}
  \framebox[1.1\width]{$\pseudo{\Alg{\cat{C}\star}}\equiv
    \pseudo{[\cat{C},\Cat]}$}
\end{center}
\noindent where $\pseudo{[\cat{C},\Cat]}$ consists of pseudo-functors, 
pseudo-natural transformations and modifications.

\vskip1em
\noindent
\underline{\textbf{Lax} $\Bimod{\cat{C}\star}$-\textbf{algebras}:} 
To obtain the corresponding objects over which pseudo-functors become 
properties, we must analyse the bicategory $\BimodT{\cat{C}\star}{\Cat}$.
A bimodule $\rel{M}{\cat{C}\star{F}}{F}$ amounts to a $C$-indexed 
collection of bimodules in $\Cat$, $\rel{M_x}{\cat{C}\star{F}(x)}{F(x)}$. 
Such $M_x$ corresponds in turn to a  family of bimodules 
$\tuple{\rel{M^f}{F(y)}{F(x)}}_{\morph{f}{y}{x}}$ according to the 
following chain of identifications

\[
\prooftree{
\coprod_{\morph{f}{y}{x}}F(y)\not\rightarrow F(x)}
\Justifies
\[\[
{\left(\coprod_{\morph{f}{y}{x}}F(y)\right)^{\mathit{op}}\times F(x) 
\rightarrow \Set}
\Justifies
{\prod_{\morph{f}{y}{x}}[F(y)^{\mathit{op}}\times F(x) \rightarrow \Set]}
\]
\Justifies 
{\tuple{F(y)\not\rightarrow F(x)}_{\morph{f}{y}{x}}}
\]
\thickness=0.08em
\endprooftree
\]

A normal monad (or equivalently a normal lax 
$\Bimod{\cat{C}\star}$-algebra) on such $M$  consists of 

\begin{description}
        \item[unit]  $\cell{\iota}{\eta_F^{*}}{M}$ which induces the 
identity
        $$\morph{\Hom{Fx} = M^{\mathit{id}_x}}{Fx}{Fx} $$
        
        \item[multiplication]  
        $\cell{m}{(\cat{C}\star{M}){\bullet}M}{\im{(\mu_F)}{\bullet}M}$ 
amounts to 
        a collection of bimodule morphisms
        
        $$ 
        \tuple{\cell{m^{f,g}}{M^{f}{\bullet}M^{g}}{M^{f{\circ}g}}}
_{f:y{\rightarrow}x, g:z{\rightarrow}y}$$
\end{description}

\noindent subject to the associativity and unit axioms which make $M$ a 
\textbf{(normal) lax functor} $\morph{M}{\cat{C}}{\Bimod{\Cat}}$ such that $Mx = 
Fx$ on objects. A morphism $h$ between two such (normal) monads 
$\rel{M}{\cat{C}\star{F}}{F}$ and $\rel{N}{\cat{C}\star{G}}{G}$ amounts 
to a $C$-collection of functors $\morph{h_x}{Fx}{Gx}$ and bimodule 
morphisms $$\cell{\theta_f}{M^f\bullet\im{(h_x)}}{\im{(h_y)}{\bullet}N^f}$$
\noindent for 
each $\morph{f}{y}{x}$ in $\cat{C}$; the compatibility with the unit and 
multiplication of the monads makes such collection into a lax 
transformation $\morph{\cell{\im{h}}{M}{N}}{\cat{C}}{\Bimod{\Cat}}$. 
Similarly, a 2-cell $\cell{\alpha}{h}{k}$ between two such morphisms 
corresponds to a modification $\cell{\im{\alpha}}{\im{h}}{\im{k}}$. Hence,

\begin{center}
  \framebox[1.1\width]{$\lax{\Alg{\Bimod{\cat{C}\star}}} \equiv
    \laxrep{[\cat{C},\Bimod{\Cat}]}$}
\end{center}
\noindent where $\laxrep{[\cat{C},\Bimod{\Cat}]}$ denotes the 2-category 
of (normal) lax 
functors, representable lax transformations (those whose component 
bimodules are representable) and modifications.


\begin{remark}
\label{rmk:fibrations}
        Here we could proceed further from the general theory and identify
         $\laxrep{[\cat{C}.\Bimod{\Cat}]}$  with 
$\Cat/\cat{C}$ (a 
         result usually credited to B\'{e}nabou). Now, adjoint pseudo-algebras on 
$\Cat/\cat{C}$ are 
         Grothendieck (co)fibrations, as expected: 
$$ {\mathit Strict\/}\laxrep{[\cat{C},\Bimod{\Cat}]} \simeq {\mathit Split Cofibrations\/}/\cat{C}$$
$$ {\mathit Representable\/}\laxrep{[\cat{C},\Bimod{\Cat}]} \simeq {\mathit Cofibrations\/}/\cat{C}$$
\end{remark}

\begin{remark}
\label{rmk:representable-pseudo-functors}
  In this situation it is fairly straightforward to see Theorem \ref{thm:representable-characterisation} at work: 
a lax functor $\morph{M}{\cat{C}}{\Bimod{\Cat}}$ is representable iff each bimodule $\rel{M_f}{Mx}{My}$ is so, thus 
$M_f = \im{Ff}$ for some functor $\morph{Ff}{Fx}{Fy}$. The structural 2-cell $m$ is an isomorphism iff
each $\tuple{\cell{m^{f,g}}{M^{g}{\bullet}M^{f}}{M^{f{\circ}g}}}
_{f:y{\rightarrow}x, g:z{\rightarrow}y}$ is invertible. Clearly such isomorphisms give the required
$\isomorph{m_{f,g}}{Ff{\circ}Fg}{F(f{\circ}g)}$ which make $\morph{F}{\cat{C}}{\Cat}$ a pseudo-functor.
So ${\mathit Representable-\/}\laxrep{[\cat{C},\Bimod{\Cat}]} \simeq \pseudo{[\cat{C},\Cat]}$.
\end{remark}

Once again Theorem \ref{thm:classification-strong-morphisms} applies to obtain the standard coherence results.

\begin{corollary}[Coherence for pseudo-functors]
\label{cor:coherence-pseudo-functors}
\mbox{ }\hfill{ }
\begin{itemize}
\item 

The inclusion
        
  $${\mathit Strict\/}\laxrep{[\cat{C},\Bimod{\Cat}]}
\hookrightarrow
{\mathit Representable\/}\laxrep{[\cat{C},\Bimod{\Cat}]}
$$
        
        \noindent has a left biadjoint whose unit is a (pseudo-natural) equivalence
\item The inclusion
        
  $$[\cat{C},\Cat]\hookrightarrow\pseudo{[\cat{C},\Cat]}
$$
        
        \noindent has a left biadjoint whose unit is a (pseudo-natural) equivalence
\end{itemize}
  
\end{corollary}

{\small
\bibliographystyle{alpha}
\bibliography{references}
}

\end{document}